\pdfoutput=1
\RequirePackage{ifpdf}
\ifpdf 
\documentclass[pdftex]{sigma}
\else
\documentclass{sigma}
\fi
\usepackage{mathabx}
\usepackage{algorithmic,algorithm}

\newtheorem{thm}{Theorem}[section]
\newtheorem{cor}[thm]{Corollary}
\newtheorem{lem}[thm]{Lemma}
\newtheorem{prop}[thm]{Proposition}

{ \theoremstyle{definition}

\newtheorem{defn}[thm]{Definition}}

\numberwithin{equation}{section} \numberwithin{figure}{section}
\numberwithin{table}{section} 

\theoremstyle{remark}

\def\a{{\alpha}}
\def\b{{\beta}}

\def\kb{{\mathbf \it k}}

\def\jb{{\mathbf j}}
\def\kb{{\mathbf k}}

\def\sb{{\mathbf s}}
\def\tb{{\mathbf t}}

 \def\NN{{\mathbb N}}
 \def\CA{{\mathcal A}}
 
 \def\CP{{\mathcal P}}

 \def\CL{{\mathcal L}}

 \def\CV{{\mathcal V}}
 \def\CH{{\mathcal H}}
 
 \def\HH{{\mathbb H}}

 \def\RR{{\mathbb R}}
 \def\ZZ{{\mathbb Z}}
 \def\A{{\mathcal A}}
 \def\G{{\mathcal G}}
 \def\S{{\mathcal S}}

\newcommand{\e}{\mathrm{e}}

\newcommand{\tr}{{\mathsf {tr}}}
\newcommand{\TC}{{\mathsf {C}}}
\newcommand{\TS}{{\mathsf {S}}}
 \def\TSS{{\mathsf{SS}}}
 \def\TCC{{\mathsf{CC}}}
 \def\TCS{{\mathsf{CS}}}
 \def\TSC{{\mathsf{SC}}}
 \def\cc{\operatorname{cc}}
 \def\sc{\operatorname{sc}}
 \def\cs{\operatorname{cs}}
 \def\ss{\operatorname{ss}}

 \def\sspan{\operatorname{span}}
\def \la {\langle}
\def \ra {\rangle}
\newcommand{\wt}{\widetilde}
\newcommand{\wh}{\widehat}

\begin{document}

\allowdisplaybreaks

\renewcommand{\PaperNumber}{067}

\FirstPageHeading

\ShortArticleName{Discrete Fourier Analysis and Chebyshev Polynomials with $G_2$ Group}

\ArticleName{Discrete Fourier Analysis and Chebyshev Polynomials with $\boldsymbol{G_2}$ Group}

\Author{Huiyuan LI~$^\dag$, Jiachang SUN~$^\dag$ and Yuan XU~$^\ddag$}

\AuthorNameForHeading{H.~Li, J.~Sun and Y.~Xu}

\Address{$^\dag$~Institute of Software, Chinese Academy of Sciences, Beijing 100190, China}
\EmailD{\href{mailto:huiyuan@iscas.ac.cn}{huiyuan@iscas.ac.cn}, \href{mailto:sun@mail.rdcps.ac.cn}{sun@mail.rdcps.ac.cn}}

\Address{$^\ddag$~Department of Mathematics, University of Oregon,
    Eugene, Oregon 97403-1222, USA}
\EmailD{\href{mailto:yuan@uoregon.edu}{yuan@uoregon.edu}}
\URLaddressD{\url{http://uoregon.edu/~yuan/}}

\ArticleDates{Received May 04, 2012, in f\/inal form September 06, 2012; Published online October 03, 2012}

\Abstract{The discrete Fourier analysis on the $30^{\degree}$--$60^{\degree}$--$90^{\degree}$ triangle
is deduced from the corresponding results on the regular hexagon by considering
functions invariant under the group $G_2$, which leads to the def\/inition of four
families generalized Chebyshev polynomials. The study of these polynomials
leads to a Sturm--Liouville eigenvalue problem that contains two parameters, whose
solutions are analogues of the Jacobi polynomials. Under a concept of $m$-degree
and by introducing a new ordering among monomials, these polynomials are
shown to share properties of the ordinary orthogonal polynomials. In
particular, their common zeros generate cubature rules of Gauss type.}

\Keywords{discrete Fourier series; trigonometric; group $G_2$; PDE; orthogonal polynomials}
\Classification{41A05; 41A10}


\section{Introduction}

In our recent works \cite{LSX, LSX10, LX09} we studied discrete Fourier analysis
associated with translation lattices. In the case of two dimension, our results include
discrete Fourier analysis of exponential functions on the regular hexagon and, by
restricting to symmetric and antisymmetric exponentials on the hexagon under the
ref\/lection group $\CA_2$ (the group of symmetry of the regu\-lar hexagon), the generalized
cosine and sine functions on the equilateral triangle, which can also be transformed into
the generalized Chebyshev polynomials on a domain bounded by the hypocycloid.
These polynomials possess maximal number of common zeros, which implies
the existence of Gaussian cubature rules, a rarity that is only the second example
ever found. The f\/irst example of Gaussian cubature rules is connected with the
trigonometric functions on the $45^{\degree}$--$45^{\degree}$--$90^{\degree}$ triangle. The richness
of these results prompts us to look into similar results on the $30^{\degree}$--$60^{\degree}$--$90^{\degree}$
triangle in the present work. This case is also considered recently in \cite{Patera} as an example
under a general framework of cubature rules and orthogonal polynomials for the
compact simple Lie groups, for which the group is $G_2$.

It turns out that much of the discrete Fourier analysis on the $30^{\degree}$--$60^{\degree}$--$90^{\degree}$
triangle can be obtained, perhaps not surprisingly,  though symmetry from our results
on the hexagonal domain. The most direct way of deduction, however, is not through our
results on the equilateral triangle. The reason lies in the underline group $G_2$, which is
a composition of $\CA_2$ and its dual $\CA_2^*$, the symmetric group of the regular
hexagon and its rotation. Our framework of discrete Fourier analysis incorporates two
lattices, one determines the domain and the other determines the space of exponentials.
Our results on the equilateral triangle are obtained from the situation when both lattices
are taken to be the same hexagonal lattices~\cite{LSX}. Another choice is to take one
lattice as the hexagonal lattice and the other as the rotation of the same lattice by $90^{\degree}$
degree~\cite{LSX10}, with the symmetric groups $\CA_2$ and $\CA_2^*$, respectively.
As we shall see, it is from this set up that our results on the $30^{\degree}$--$60^{\degree}$--$90^{\degree}$
triangle can be deduced directly via symmetry. The results include cubature rules and
orthogonal trigonometric functions that are analogues of cosine and sine functions. There
are four families of such functions and they have also been studied recently in~\cite{Patera, Sz}.
While the results in these two papers concern mainly with orthogonal polynomials, our emphasis
is on the discrete Fourier analysis and cubature rules, and on the connection to the results in
the hexagonal domain.

The generalized cosine and sine functions on the $30^{\degree}$--$60^{\degree}$--$90^{\degree}$
triangle are also eigenfunctions of the Laplace operator with suitable boundary conditions. There
are four families of such functions. Under proper change of variables, they become orthogonal
polynomials on a domain bounded by two curves. However, unlike the equilateral triangle,
these polynomials do not form a complete orthogonal basis in the usual sense of total order of
monomials. To understand the structure of these polynomials, we consider the Sturm--Liouville
problem for a general pair of parameters~$\alpha$,~$\beta$, with the four families that correspond to
the generalized cosine and sine functions as $\alpha = \pm \frac 12$, $\beta = \pm \frac12$.
The dif\/ferential operator of this eigenvalue problem has the form
\begin{gather*}
 \CL_{\a,\b} := -A_{11}(x,y) \partial_x^2 -2A_{12} (x,y)\partial_x \partial_y  - A_{22}(x,y)\partial_y^2
   + B_1 (x,y)\partial_x + B_2 (x,y) \partial_y.
\end{gather*}
Such operators have long been studied in association with orthogonal polynomials in two
va\-riab\-les; see for example \cite{K, K2, KS, Su}, as well as \cite{Beer} and the references therein.
Our opera\-tor~$\CL_{\a,\b}$, however, is dif\/ferent in the sense that the coef\/f\/icient functions~$A_{i,j}$
are usually assumed to be of quadratic polynomials to ensure that the operator has $n+1$
polynomials of degree $n$ as eigenfunctions, whereas $A_{2,2}$ in our $L_{\a,\b}$ is a polynomial
of degree $3$ for which it is no longer obvious that a full set of eigenfunctions exists. Nevertheless,
we shall prove that the eigenvalue problem $\CL_{\a,\b} u = \lambda u$ has a complete set of
polynomial solutions, which are also orthogonal polynomials, analogue of the Jacobi polynomials.
Upon introducing a new or\-dering among monomials, these
polynomials  can be shown to be uniquely determined by their highest term in the new
ordering.  As a matter of fact, this ordering def\/ines the region of  inf\/luence and dependence
in the polynomial space for each solution. Furthermore, it
preserves the $m$-degree of polynomials, a~concept introduced
in \cite{Patera}, rather than the total degree. In the case of $\a = \pm \frac12$ and
$\b = \pm \frac12$, the common zeros of these polynomials determine the Gauss,
Gauss--Lobatto and Gauss--Radau cubature rules, respectively, all in the sense of $m$-degree.
It is known that the cubature rule of degree $2n-1$ exists if and only if its nodes form a
variety of an ideal generated by certain ortho\-go\-nal polynomials. It is somewhat surprising
that this relation is preserved when the $m$-degree is used in place of the ordinary degree.

The paper is organized as follows. The following section contains what we need from the
discrete Fourier analysis on the hexagonal domain. The results on the   $30^{\degree}$--$60^{\degree}$--$90^{\degree}$
triangle is developed in Section~\ref{section3}, which are translated into generalized Chebyshev polynomials
in Section~\ref{section4}. The Sturm--Liouville problem is def\/ined and studied in Section~\ref{section5} and the
cubature rules are presented in Section~\ref{section6}.

\section{Discrete Fourier analysis on hexagonal domain}\label{section2}

Before stating the results on the hexagonal domain, we give a short narrative of
the necessary back\-ground on the discrete Fourial analysis with lattice as developed
in \cite{LSX, LX09}. We refer to~\mbox{\cite{CS,DM, Ma,Mun}} for some applications of discrete
Fourier analysis in several variables.

A lattice $L$ in $\RR^d$ is a discrete subgroup $L =  L_A := A\ZZ^d$, where
$A$, called a generator matrix, is nonsingular. A bounded set $\Omega$
of $\RR^d$, called the fundamental domain of $L$, is said to tile
$\RR^d$ with the lattice $L$ if $\Omega + L = \RR^d$, that is,
\[
  \sum_{\alpha \in L} \chi_\Omega (x + \alpha) = 1, \qquad
      \hbox{for almost all $x \in \RR^d$},
\]
where $\chi_\Omega$ denotes the characteristic function of $\Omega$. For a given
lattice $L_A$, the dual lattice $L_A^\perp$ is given by $L_A^\perp =  A^{- \tr}\ZZ^d$.
A result of Fuglede \cite{F} states that a bounded open set $\Omega$ tiles
$\RR^d$ with the lattice $L$ if, and only if, $\{\e^{2 \pi i \alpha \cdot x}: \alpha
\in L^\perp\}$ is an orthonormal basis with respect to the inner product
\begin{equation} \label{ipOmega}
  \langle f, g \rangle_\Omega = \frac{1}{|\Omega|}
      \int_\Omega f(x) \overline{g(x)} dx.
\end{equation}
Since $L^\perp_A = A^{-\tr}\ZZ^d$, we can write $\alpha = A^{-\tr} k$ for $\alpha \in
L_A^\perp$ and $k \in \ZZ^d$, so that $\e^{2\pi i \alpha \cdot x}
= \e^{2 \pi i k^{\tr} A^{-1} x}$.

For our discrete Fourier analysis, the boundary of $\Omega$ matters. We shall
f\/ix an $\Omega$ such that $0 \in \Omega$ and  $\Omega + A\ZZ^d =  \RR^d$
holds {\it pointwisely} and {\it without overlapping}.

\begin{defn} \label{def:N}
Let  $\Omega_A$ and $\Omega_B$ be the fundamental domains of $A\ZZ^d$
and $B\ZZ^d$, respectively. Assume {\it all entries of the matrix
$N:=B^\tr A$ are integers}.
Def\/ine
\[
\Lambda_N :=   \big\{ k \in \ZZ^d: B^{-\tr}k \in \Omega_A\big\} \qquad \hbox{and} \qquad
\Lambda_{N}^\dag:=   \big\{ k \in \ZZ^d: A^{-\tr}k \in \Omega_B\big\}.
\]
Furthermore, def\/ine the f\/inite-dimensional subspace of exponential functions
\[
\CV_N := \sspan \big\{\e^{2\pi i\, k^\tr A^{-1} x}, \, k \in \Lambda_{N}^\dag \big\}.
\]
\end{defn}

A function $f$ def\/ined on $\RR^d$ is called a periodic function with respect
to the lattice $A \ZZ^d$ if
\[
      f (x + A k) = f(x) \qquad \hbox{for all $k \in \ZZ^d$}.
\]
The function $x \mapsto e^{2\pi i k^\tr A^{-1} x}$ is periodic with respect to the lattice $A\ZZ^d$
and $\CV_N$ is a space of periodic exponential functions. We can now state the central result
in the discrete Fourier analysis.

\begin{thm} \label{thm:df}
Let $A$, $B$ and $N$ be as in Definition~{\rm \ref{def:N}}.  Define
\[
    \langle f, g \rangle_N =  \frac{1}{|\det (N)|}
        \sum_{j \in \Lambda_N } f(B^{-\tr} j ) \overline{g(B^{-\tr} j )}
\]
for $f$, $g$ in $C(\Omega_A)$, the space of continuous functions on~$\Omega_A$.
Then
\begin{equation}\label{c-d-inner}
    \langle f, g \rangle_{\Omega_A} =  \langle f, g \rangle_{N}, \qquad f, g \in \CV_N.
\end{equation}
\end{thm}

It follows readily that~\eqref{c-d-inner} gives a cubature formula exact for functions
in~$\CV_N$. Furthermore, it implies an explicit Lagrange interpolation by exponential
functions, which we shall not state since it will not be needed in the present work.

In the following, we shall call the lattice $L_A$ as the lattice for the physical space, as it determines the
domain on which our analysis lies, and the lattice~$L_B$  as the lattice for the frequency space,
as it determines the points that def\/ines the inner product.

The classical discrete Fourier analysis of two variables is the tensor product of the
results in one variable, which corresponds to $A =  B = I$, the identity matrix. We
are interested in choosing $A$ as the generating matrix $H$ of the hexagonal domain,
\[
   H =  \begin{pmatrix} \sqrt{3} & 0\\ -1 & 2\end{pmatrix} \qquad \hbox{with} \qquad
   \Omega_H=\left\{x\in \RR^2:\  -1\leq x_2, \tfrac{\sqrt{3} x_1}{2} \pm
\tfrac{x_2}{2}  < 1 \right\}.
\]
If we choose $B  =  \frac{n}{2} H$, so that $N = B^\tr A$ has all integer entries, we
are back to the  situation studied in~\cite{LSX}, which is the one that leads to the
discrete Fourier analysis on the equilateral triangle. The other choices are considered
in~\cite{LSX10}.

For the case that we are interested in, we choose $A = H$, the matrix for the hexagonal
lattice in the physical space, and $B  =  n H^{-\tr}$ with $n\in \ZZ$, the matrix for the
hexagonal lattice in the frequency space. Then $N=B^{\tr}A=n I$ has all integer entries.
This case was studied in \cite{LSX10}, which will be used to deduce the case that we are
interested in by an additional symmetry. As shown in \cite{LSX,Sun}, it is more
convenient to use homogeneous coordinates  $(t_1,t_2,t_3)$ def\/ined by
\begin{gather}\label{coordinates}
\begin{pmatrix}t_1 \\ t_2 \\ t_3 \end{pmatrix} =
 \begin{pmatrix}  \frac{\sqrt{3}}{2} & -\frac12 \\
                          0 &  1\\
                  -\frac{\sqrt{3}}{2} & -\frac12
 \end{pmatrix}  \begin{pmatrix}x_1 \\ x_2 \end{pmatrix} := E x,
\end{gather}
which satisfy $t_1 + t_2 +t_3 =0$. We adopt the convention of using bold letters,
such as $\tb$ to denote points in homogeneous coordinates. We def\/ine by
\[
    \RR_H^3 : = \big\{\tb = (t_1,t_2,t_3)\in \RR^3: t_1+t_2 +t_3 =0\big\} \qquad \hbox{and} \qquad
     \HH^{\dag} := \ZZ^3 \cap \RR^3_H
\]
the spaces of points and integers in homogeneous coordinates, respectively.
In such coordinates,
the fundamental  domains of the lattices $L_A$ and $L_B$ are then given by
\begin{gather*}
  \Omega:=  \Omega_A=\left\{\tb \in \RR_H^3:\  -1< t_1,t_2,-t_3\le 1 \right\},\\
\hphantom{\Omega:=}{} \ \Omega_B=\left\{\tb \in \RR_H^3:\  -n < t_1-t_2, t_1-t_3, t_2-t_3 \le n  \right\},
\end{gather*}
where $\Omega_A$ can be viewed as the intersection of the plane $t_1+t_2+t_3=0$ with the cube $[-1,1]^3$.
Def\/ine the index sets in homogeneous coordinates
\begin{gather*}
   \HH_n:  = \big\{\jb \in \HH^{\dag}: -n \le j_1,j_2,j_3\le n, \   \jb \equiv 0\ (\bmod\ 3) \big\},\\
   \HH_n^{\dag}: = \big\{\kb \in \HH^{\dag}: -n \le k_3- k_2, k_1-k_3, k_2-k_1 \le n \big\},
\end{gather*}
where $\tb \equiv 0 \pmod m$ means, by def\/inition, $t_1\equiv t_2 \equiv t_3
 \pmod m$.
We note that $\HH_n$ and $\HH_n^{\dag}$ serve as the symmetric counterparts of $\Lambda_N$ and $\Lambda_N^{\dag}$,
respectively, so that $\HH_n$ determines the points in the discrete inner product and $\HH_n^{\dag}$ determines
the space of exponentials. Moreover, the index set $\HH_n$ can be obtained from a rotation of $\HH_n^{\dag}$,
as shown in the following proposition.

\begin{figure}[htb]
\hfill%
\includegraphics[width=0.4\textwidth]{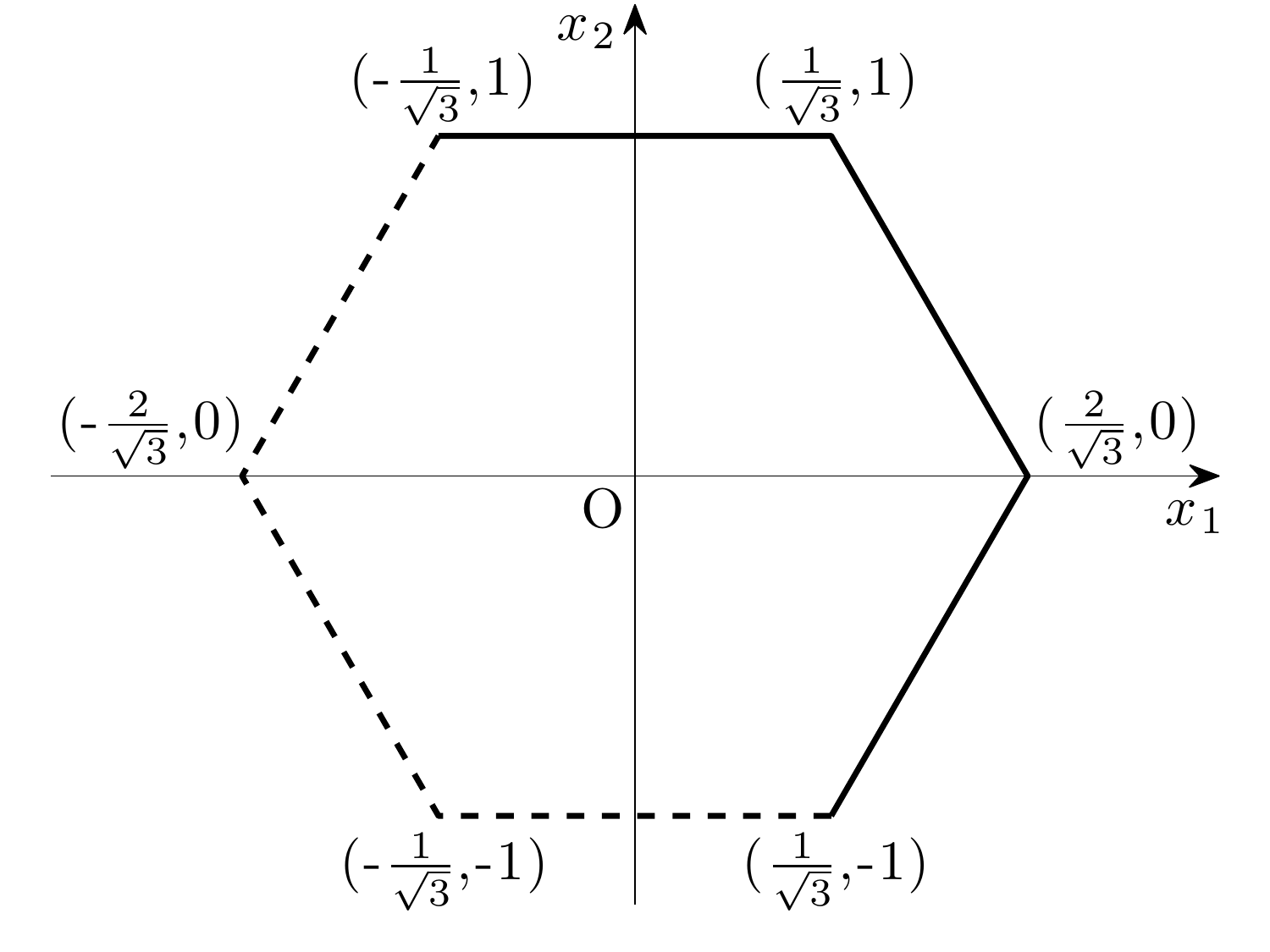}%
\hfill%
\includegraphics[width=0.4\textwidth]{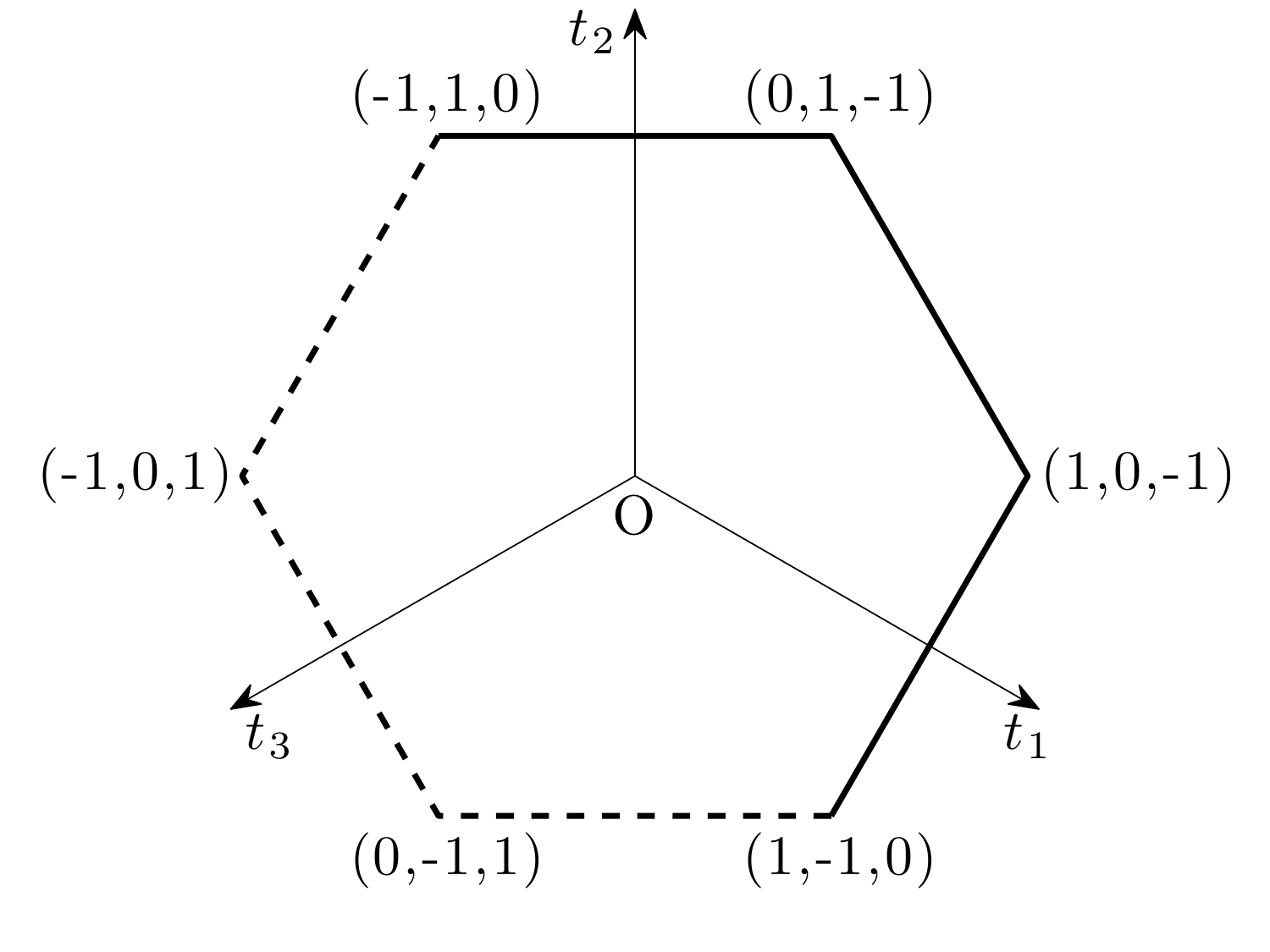}%
\hspace*{\fill}%
\caption{ $\Omega_A$ in Cartesian coordinates (left)  and homogeneous coordinates (right).}
\label{hp}
\end{figure}
\begin{figure}[htb]
\hfill%
\includegraphics[width=0.4\textwidth]{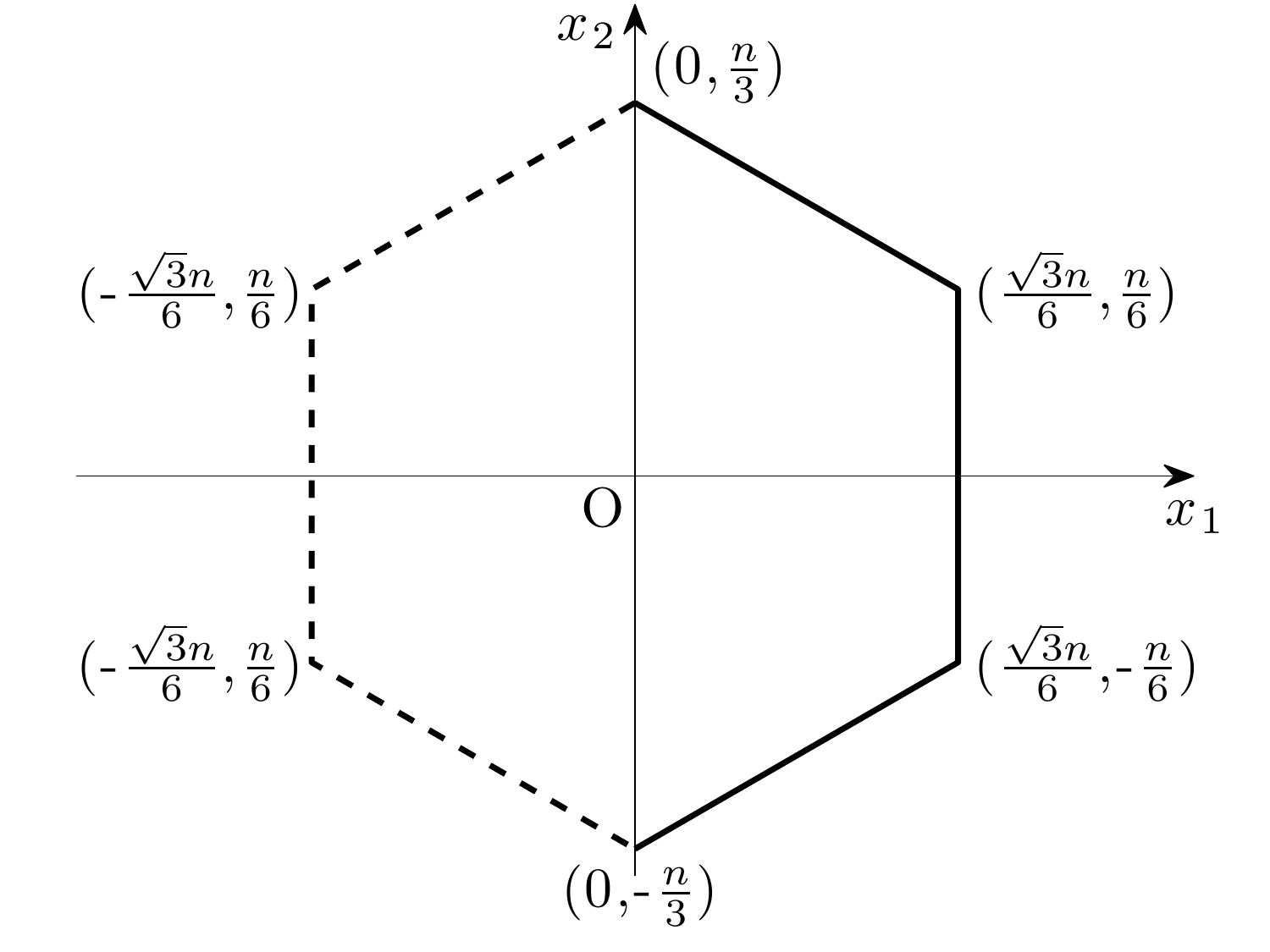}%
\hfill%
\includegraphics[width=0.4\textwidth]{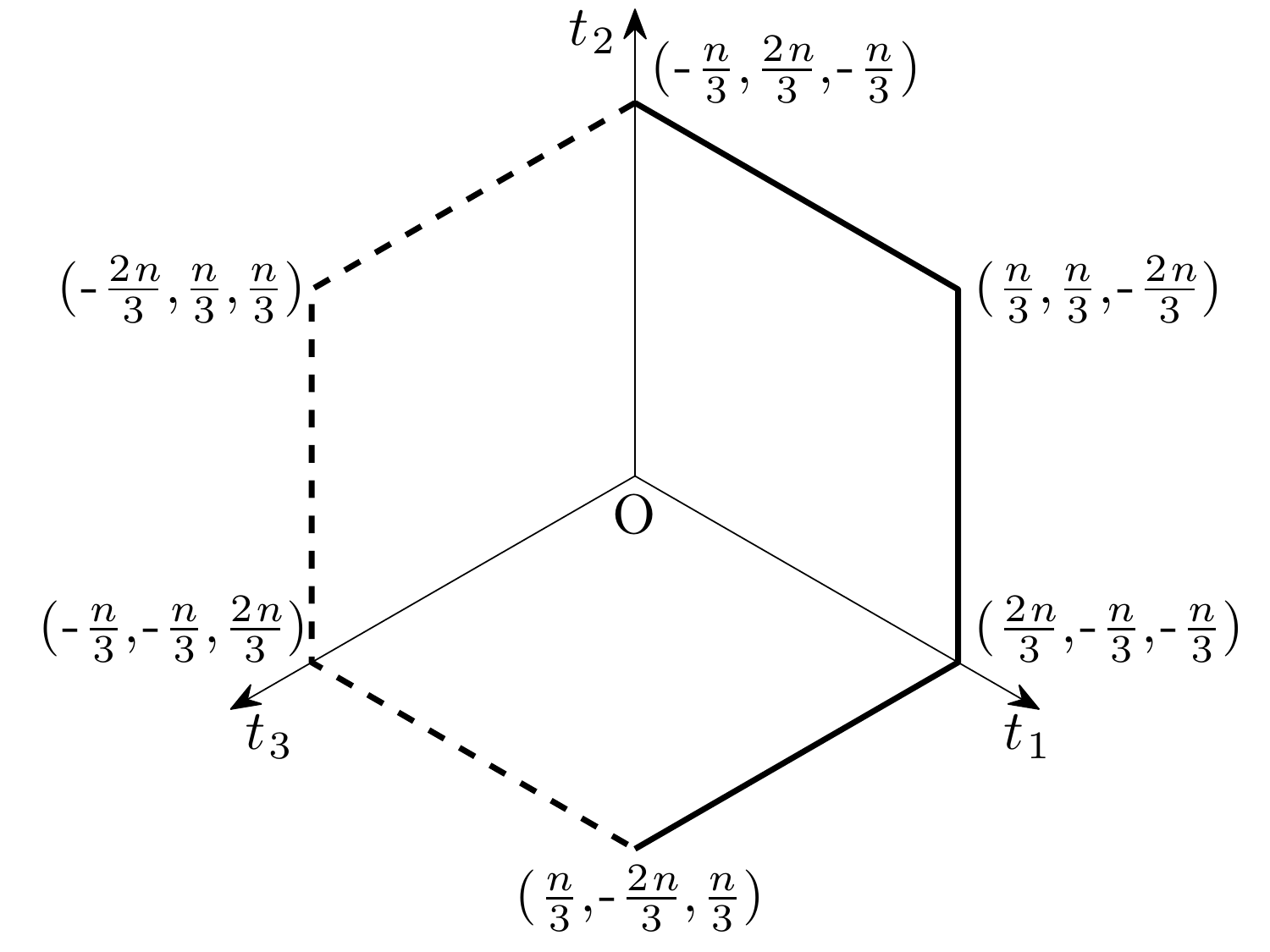}%
\hspace*{\fill}%
\caption{ $\Omega_B$ in Cartesian coordinates (left)  and homogeneous coordinates (right).}
\label{hf}
\end{figure}

\begin{prop}[\cite{LSX10}] \label{prop:H*H}
For  $\tb=(t_1,t_2,t_3)\in \RR^3_H$, define $\wh{\tb}:=(t_3-t_2,t_1-t_3,t_2-t_1)$.
Then $\frac{\wh\kb}{3}\in \HH_n^{\dagger}$ if $\kb\in \HH_n$ and
 $\wh\kb\in \HH_n$ if $\kb\in \HH_n^{\dagger}$.
\end{prop}

Proposition~\ref{prop:H*H} states that $\HH_n=\wh{\HH}_n^{\dag}:=\big\{\wh{\kb}:\;  \kb\in \HH^{\dag}_n \big\}$.
Similarly, we can def\/ine $\HH=\wh{\HH}^{\dag}:=\big\{\wh{\kb}:\;  \kb\in \HH^{\dag} \big\}=
\big\{\jb\in \HH^{\dag}:\;  \jb\equiv 0 \pmod{3} \big\}$. The set $\HH_n^{\dag}$ is the index set for the
space of expo\-nen\-tials.
Def\/ine the f\/inite-dimensional space $\CH_n^{\dagger}$ of exponential functions
\[
   \CH_n^{\dagger}: =   \sspan \left \{ \phi_\kb=\e^{\frac{2i\pi}{3} \kb \cdot \tb }: \kb \in  \HH_n^{\dagger} \right \}.
\]
By induction, it is not dif\/f\/icult to verify that
\[
\dim \CH_n^{\dagger} = |\HH_n^{\dagger}|=|\HH_n|  = \begin{cases} n^2 + n +1, \quad \hbox{if $n\not \equiv 1 \pmod 3$},\\
    n^2 + n -1, \quad \hbox{if $n\equiv 1 \pmod 3$}.
    \end{cases}
\]
\begin{figure}[htb]
\hfill%
\includegraphics[width=0.31\textwidth]{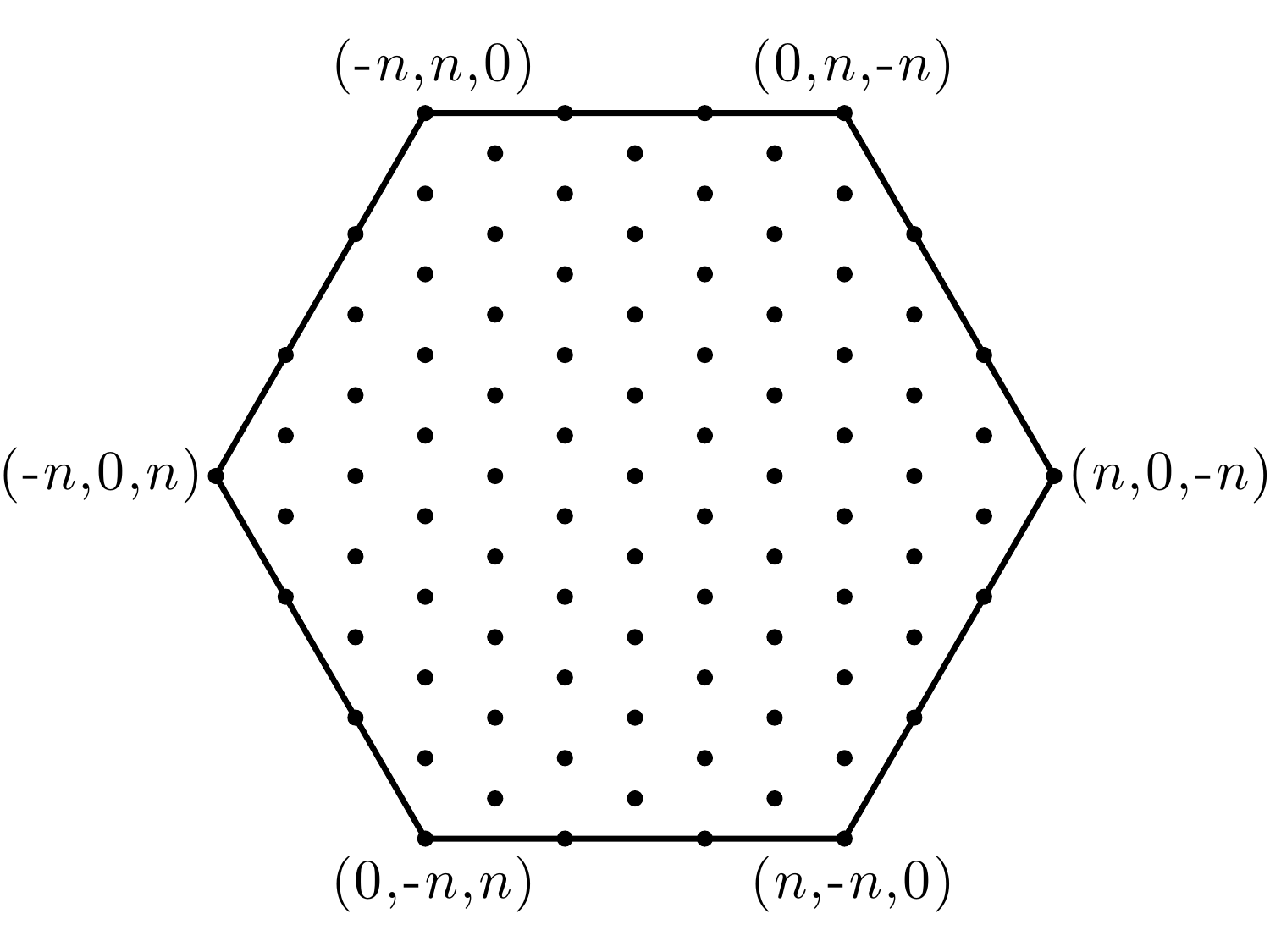}%
\hfill%
\includegraphics[width=0.31\textwidth]{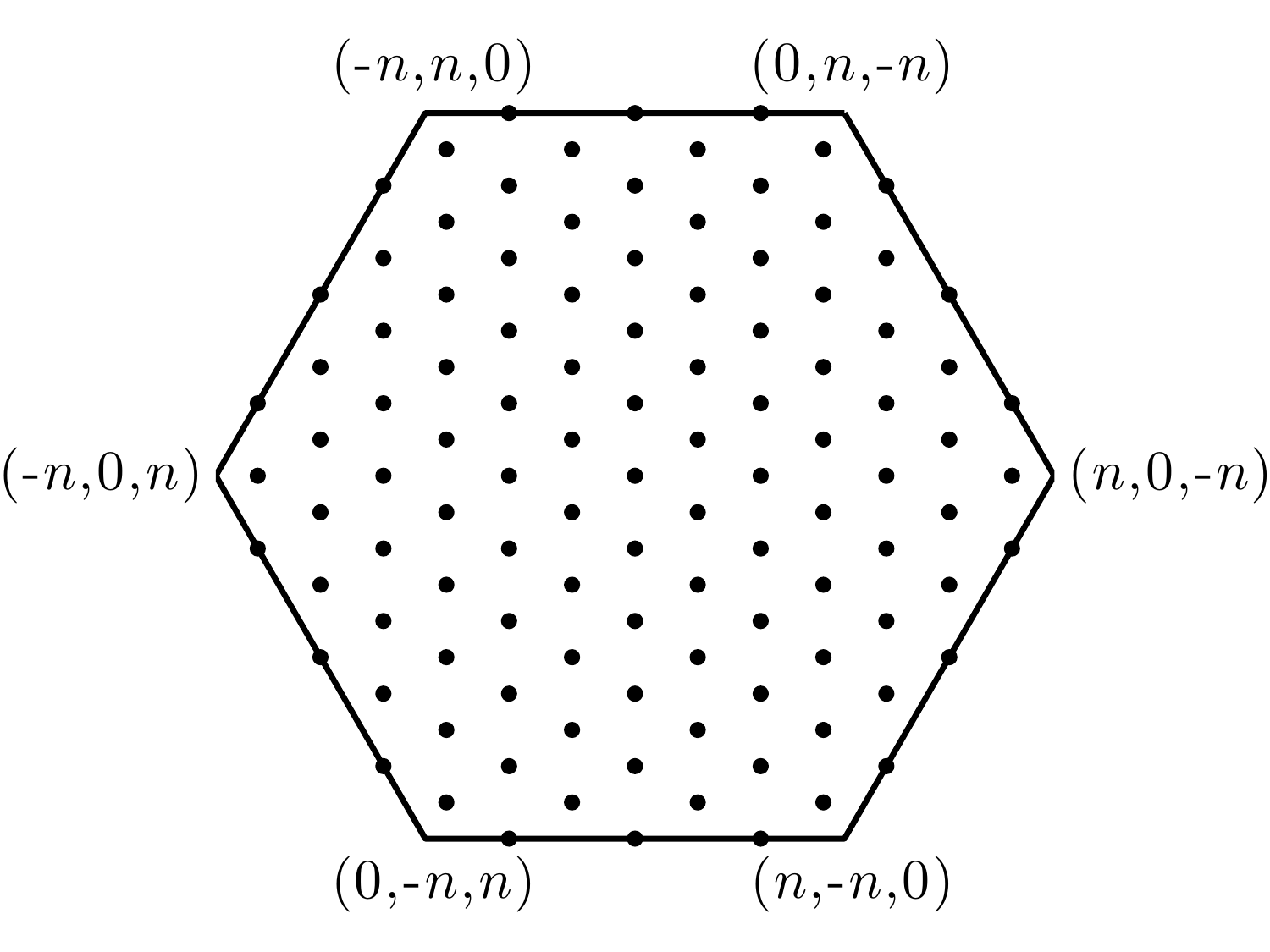}%
\hfill%
\includegraphics[width=0.31\textwidth]{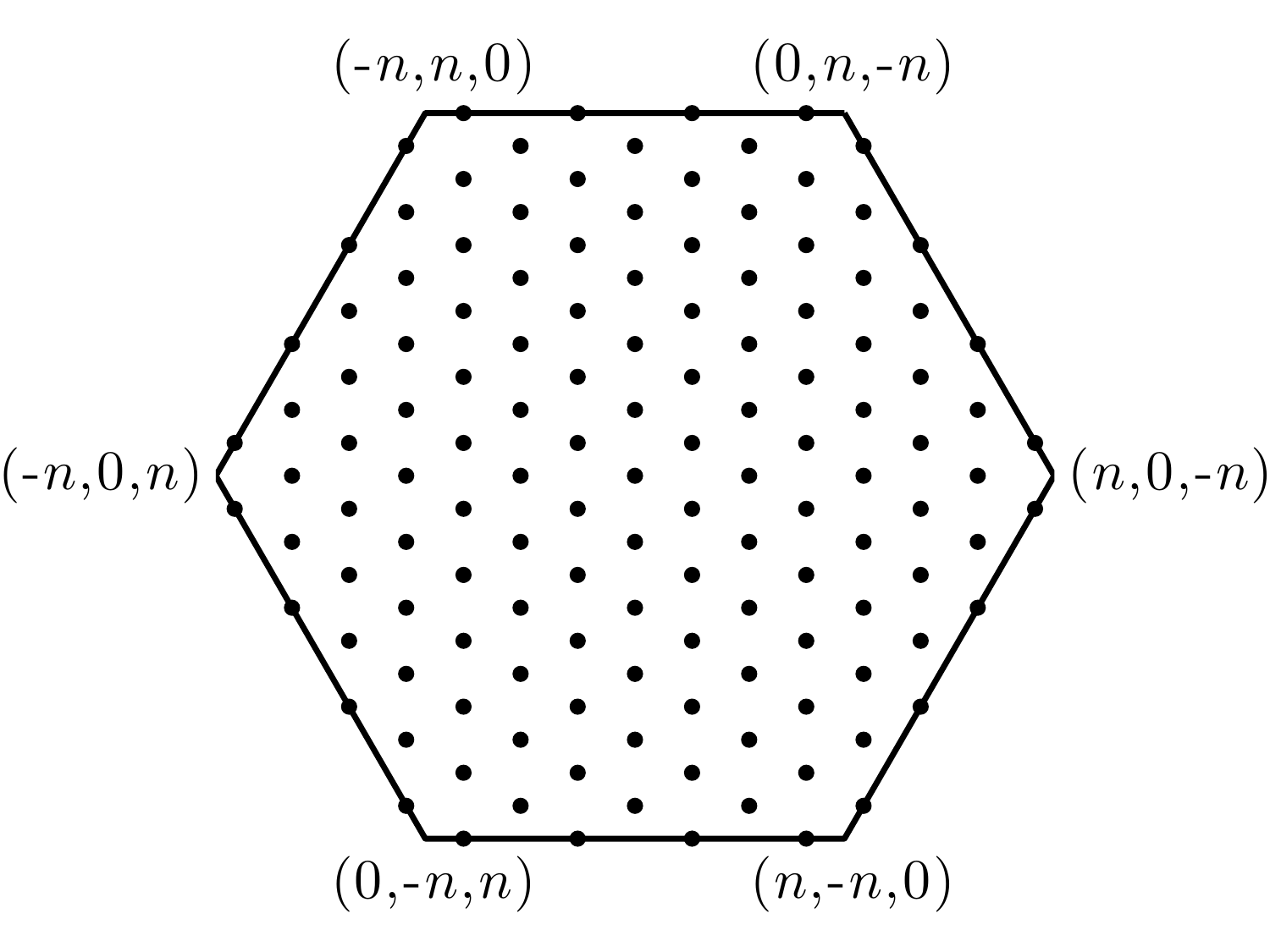}%
\hspace*{\fill}%
\caption{ $\HH_n$ for $n=9$ (left), $n=10$ (center) and $n=11$ (right).}
\label{hp306090}
\end{figure}
\begin{figure}[htb]
\hfill%
\includegraphics[width=0.31\textwidth]{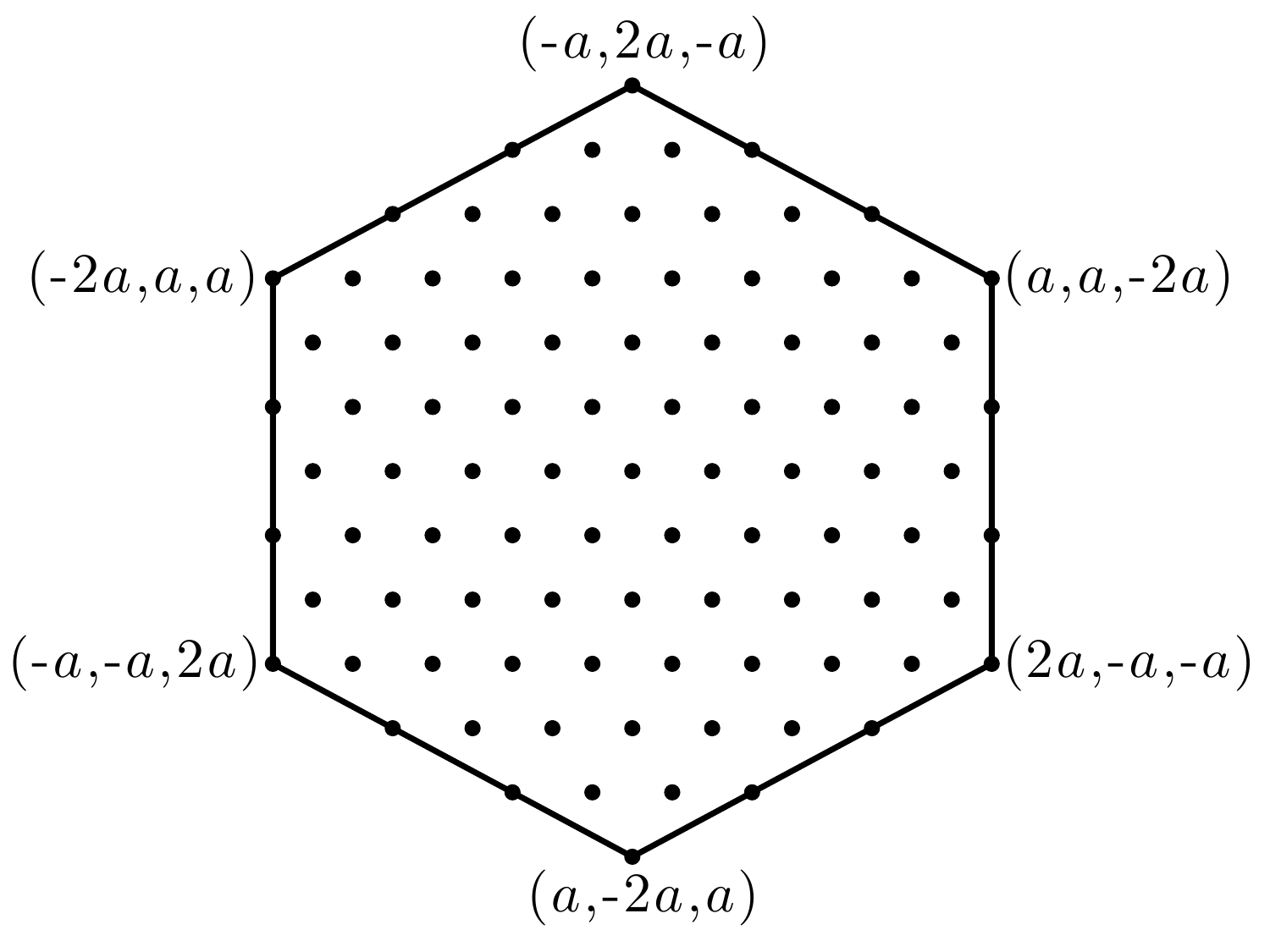}%
\hfill%
\includegraphics[width=0.31\textwidth]{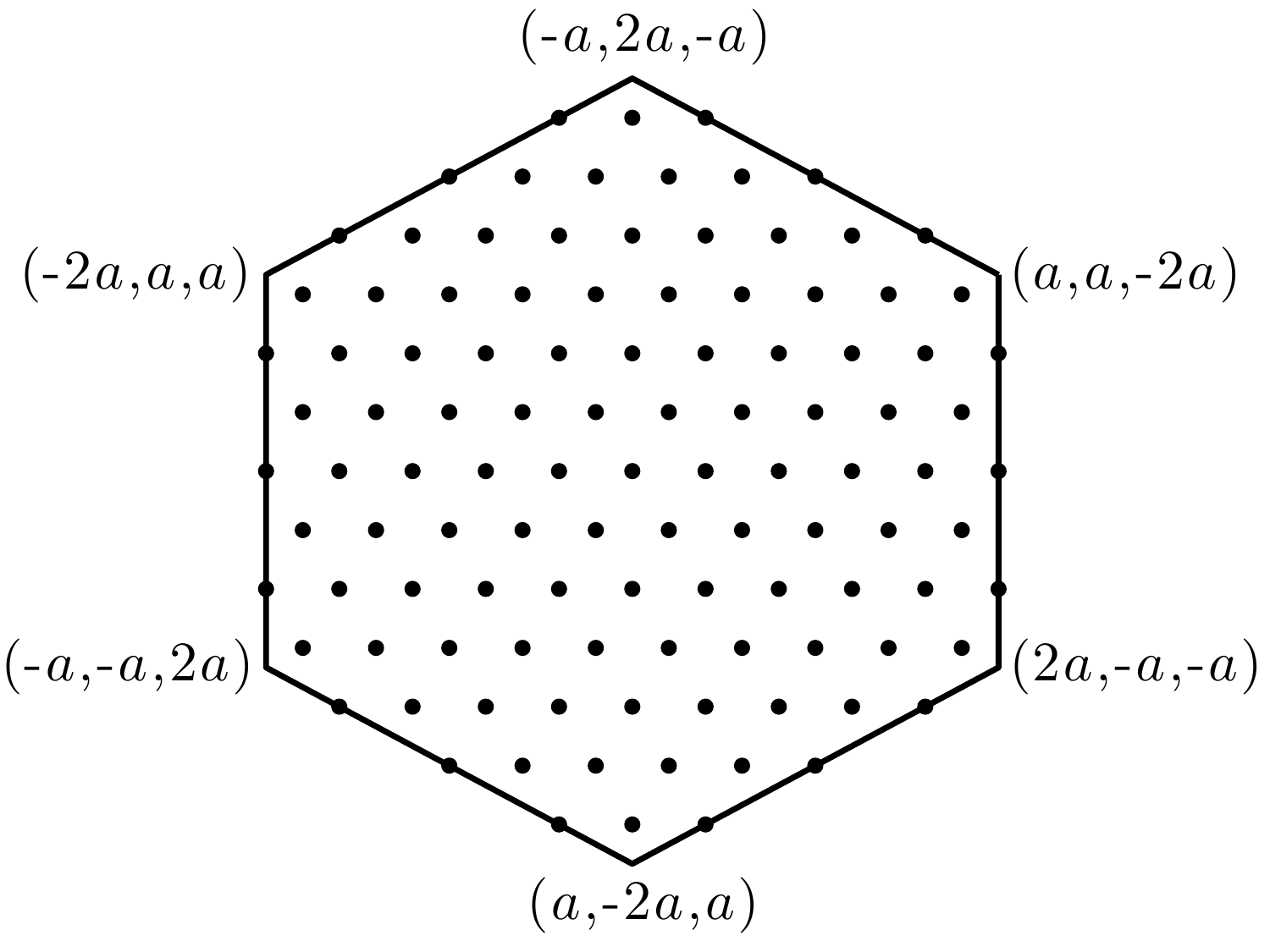}%
\hfill%
\includegraphics[width=0.31\textwidth]{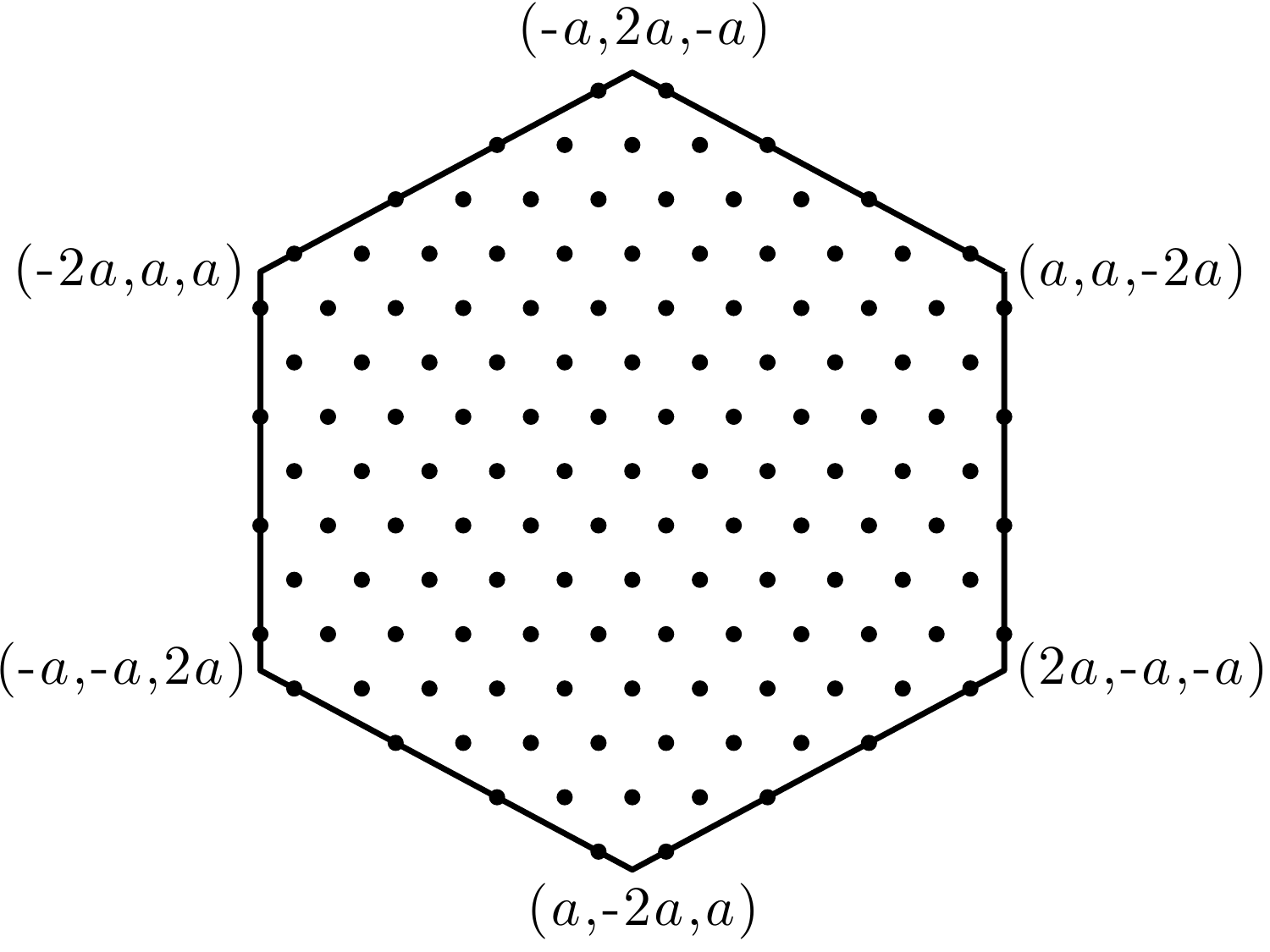}%
\hspace*{\fill}%
\caption{ $\HH^{\dagger}_n$ for $n=9$ (left), $n=10$ (center) and $n=11$ (right), where $a=\frac{n}{3}$.}
\label{hf306090}
\end{figure}
Under the homogeneous coordinates \eqref{coordinates}, $x\equiv y \pmod{H}$ becomes  $\tb  \equiv \sb \pmod 3$.
We call a function $f$ H-periodic if $f(\tb ) = f(\tb + \jb)$ whenever  $\jb  \equiv 0\, (\bmod 3)$.  Since $\jb, \kb \in \HH$
implies that $2 \jb \cdot \kb = (j_1-j_2)(k_1-k_2) + 3 j_3k_3$, we see that $\phi_\jb$ is H-periodic.

\begin{thm}[\cite{LSX10}]\label{thm:HH2ip}
The following cubature rule holds for any $f\in \CH_{2n-1}^{\dagger}$,
\begin{gather}  \label{cuba-HHD}
   \frac{1}{|\Omega|} \int_{\Omega} f(\tb)d\tb  = \frac{1}{n^2} \sum_{\jb\in \HH_n} c_\jb^{(n)}
   f\big(\tfrac{\jb}{n}\big), \qquad      c_{\jb}^{(n)} = \begin{cases} 1, & \jb \in \HH_n^{\degree},\\
       \frac{1}{2},  & \jb \in \HH_n^e,\\ \frac{1}{3},    & \jb \in \HH_n^v, \end{cases}
\end{gather}
where $\HH_n^{{\degree}}$, $\HH_n^{v}$ and $\HH_n^{e}$ denote the set of points in interior, set of vertices,
and set of points on the edges but not on the vertices; more precisely,
$\HH_n^{{\degree}}=\left\{\jb\in \HH:\,  -n<j_1,j_2,j_3 <n \right\}$, $\HH_n^{v}=\left\{(n,0,-n)\sigma \in \HH:\,
\sigma\in \A_2 \right\}$ and $\HH_n^{e}=\HH_n\setminus (\HH_n^{{\degree}}\cup \HH_n^{v})=
\left\{ (j,n-j,-n)\sigma\in \HH:\,   1\le j\le n-1 \right\}$.
In particular, let $Q_n f$ denote the right hand side of \eqref{cuba-HHD}; then for any $\kb\in\HH^{\dag}$,
$Q_n\phi_\kb = 1$ if  $\hat \kb \equiv 0 \pmod{3n}$ and $Q_n\phi_\kb = 0$ otherwise.
\end{thm}

Here we state the main result in terms of the cubature rule \eqref{cuba-HHD}, from which the
discrete inner product can be easily deduced. For further results in this regard, including
interpolation, we refer to \cite{LSX10}.

\section[Discrete Fourier analysis on the $30^{{\degree}}$--$60^{{\degree}}$--$90^{{\degree}}$ triangle]{Discrete Fourier analysis on the $\boldsymbol{30^{{\degree}}}$--$\boldsymbol{60^{{\degree}}}$--$\boldsymbol{90^{{\degree}}}$ triangle}\label{section3}

In this section we deduce a discrete Fourier analysis on the $30^{{\degree}}$--$60^{{\degree}}$--$90^{{\degree}}$ triangle
from the analysis on the hexagon by working with invariant functions.

\subsection{Generalized trigonometric functions}

The group $\A_2$  is generated by the ref\/lections in the edges of the equilateral triangles
inside the regular hexagon~$\Omega$. In homogeneous coordinates, the three ref\/lections
$\sigma_1$, $\sigma_2$, $\sigma_3$ are def\/ined by
\[
  \tb  \sigma_1 :=  -(t_1,t_3,t_2),  \qquad \tb \sigma_2 := -(t_2,t_1,t_3),
     \qquad \tb\sigma_3:= -(t_3,t_2,t_1).
\]
Because of the relations $\sigma_3=\sigma_1\sigma_2\sigma_1=\sigma_2\sigma_1\sigma_2$,
the group is given by
\[
\mathcal{A}_2 = \left\{ 1, \sigma_1, \sigma_2, \sigma_3, \sigma_1\sigma_2
, \sigma_2\sigma_1\right\}.
\]

The group $\A_2^*$ of isometries of the hexagonal lattice is generated by the ref\/lections
in the median of the equilateral triangles inside it, which can be derived from  the ref\/lection group
$\A_2$ by a rotation of $90^{\degree}$ and is exactly the permutation group of three elements.
To describe the elements in $\A_2^*$, we def\/ine the ref\/lection $-\sigma$ for any $\sigma\in \A_2$ by
\[
            \tb(-\sigma) : = -\tb\sigma, \qquad \forall\, \tb \in \RR_H^3.
\]
With this notation, the group $\A_2^*$ is given by
\[
\mathcal{A}_2^* =\left\{ 1, -\sigma_1, -\sigma_2, -\sigma_3, \sigma_1\sigma_2, \sigma_2\sigma_1\right\},
\]
in which   $-\sigma_1$, $-\sigma_2$, $-\sigma_3$ serve as the three basic ref\/lections. The group~$A_2^*$ is
the same as the permutation group $\S_3$ with three elements.

The group $G_2$ is exactly the composition of $\A_2$ and $\A_2^*$,
\[
   G_2 = \left\{\sigma\sigma^*: \sigma\in \A_2, \sigma^*\in \A_2^*\right\}
 =\left\{ \pm 1, \pm \sigma_1, \pm \sigma_2, \pm \sigma_3, \pm \sigma_1\sigma_2, \pm \sigma_2\sigma_1\right\}.
\]

Let $\G$ denote the group of $\A_2$ or $\A_2^*$ or $G_2$. For a function $f$ in homogeneous coordinates,
the action of the group $\G$ on $f$ is def\/ined by $\sigma f(\tb) = f (\tb \sigma)$, $\sigma \in \G$. A function $f$
is called {\it invariant} under $\G$ if $\sigma f =f$ for all $\sigma \in \G$, and called {\it anti-invariant} under $\G$
if $\sigma f = (-1)^{|\sigma|}f$ for all $\sigma \in \G$, where $|\sigma|$ denotes the inversion of $\sigma$ and
$(-1)^{|\sigma|} = 1$ if $\sigma =\pm 1, \pm \sigma_1 \sigma 2, \pm \sigma_2\sigma_1$,
and $(-1)^{|\sigma|} = -1$ if $\sigma = \pm \sigma_1,\pm \sigma_2,\pm \sigma_3$.
The following proposition is easy to verify (see \cite{K}).

\begin{prop}
Define the operators $\CP^+$ and $\CP^-$ acting on $f(\tb)$ by
\begin{equation} \label{CP^+}
\CP^\pm f(\tb) =  \frac{1}{6} \left[f(\tb) + f(\tb \sigma_1\sigma_2)+f(\tb \sigma_2\sigma_1)
        \pm  f(\tb \sigma_1) \pm  f(\tb \sigma_2)  \pm  f(\tb \sigma_3) \right].
\end{equation}
Then the operators $\CP^+$ and $\CP^-$ are projections from the class of
H-periodic functions onto the class of invariant, respectively anti-invariant,
functions under $\A_2$.
Furthermore,  define the operators $\CP_{\!*}^+$ and $\CP_{\!*}^-$ acting on $f(\tb)$ by
\begin{equation} \label{CP*^+}
\CP^\pm_{\!*} f(\tb) =  \frac{1}{6}  \left[f(\tb)+ f(\tb\sigma_1\sigma_2)  + f(\tb\sigma_2\sigma_1)
 \pm f(-\tb\sigma_1) \pm f(-\tb\sigma_2) \pm f(-\tb\sigma_3) \right].
\end{equation}
Then the operators $\CP_{\!*}^+$ and $\CP_{\!*}^-$ are projections from the class of
H-periodic functions onto the class of invariant, respectively anti-invariant
functions under $\A_2^*$.
\end{prop}
\begin{figure}[htb]
\hfill%
\includegraphics[width=0.31\textwidth]{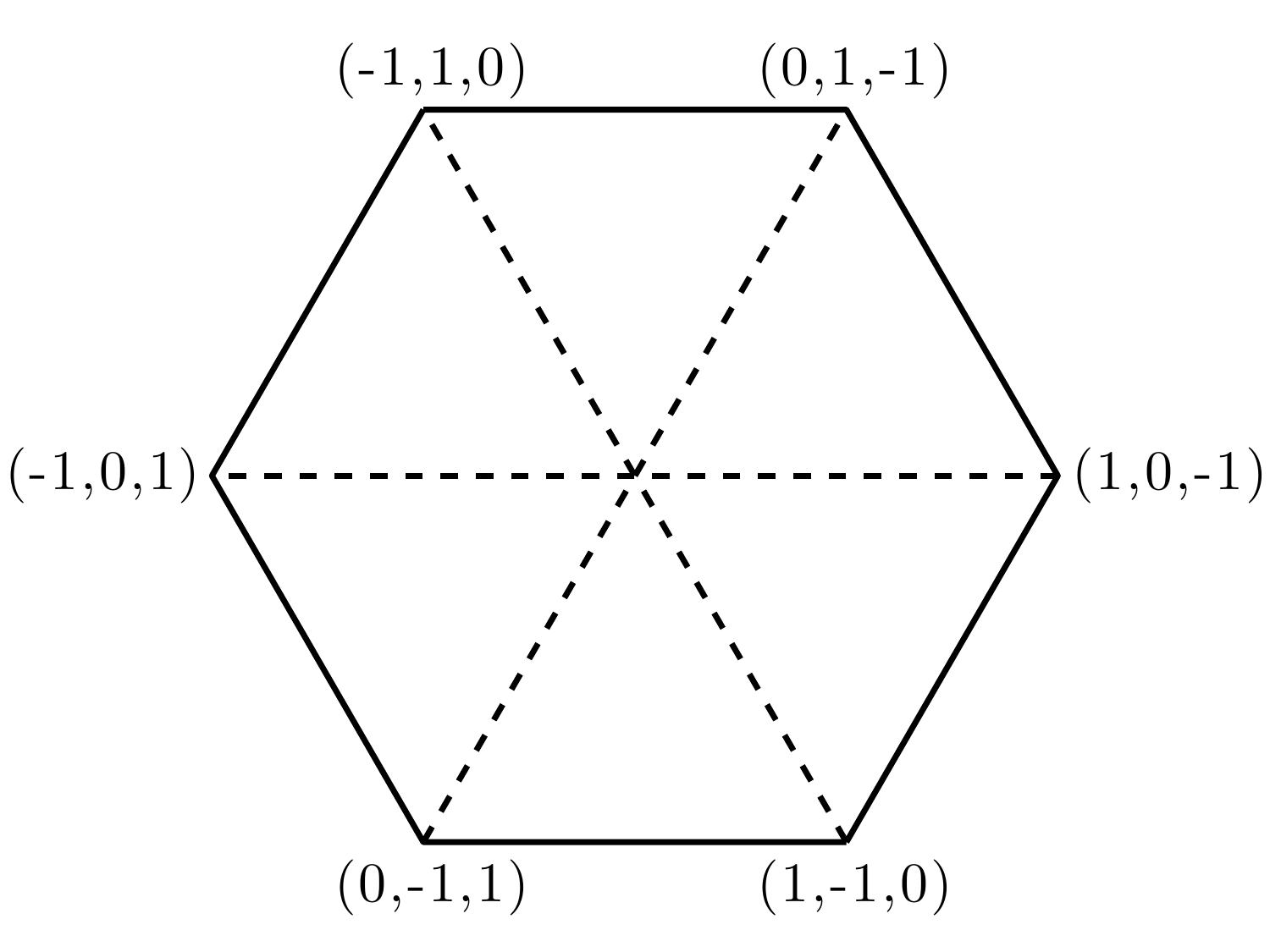}%
\hfill%
\includegraphics[width=0.31\textwidth]{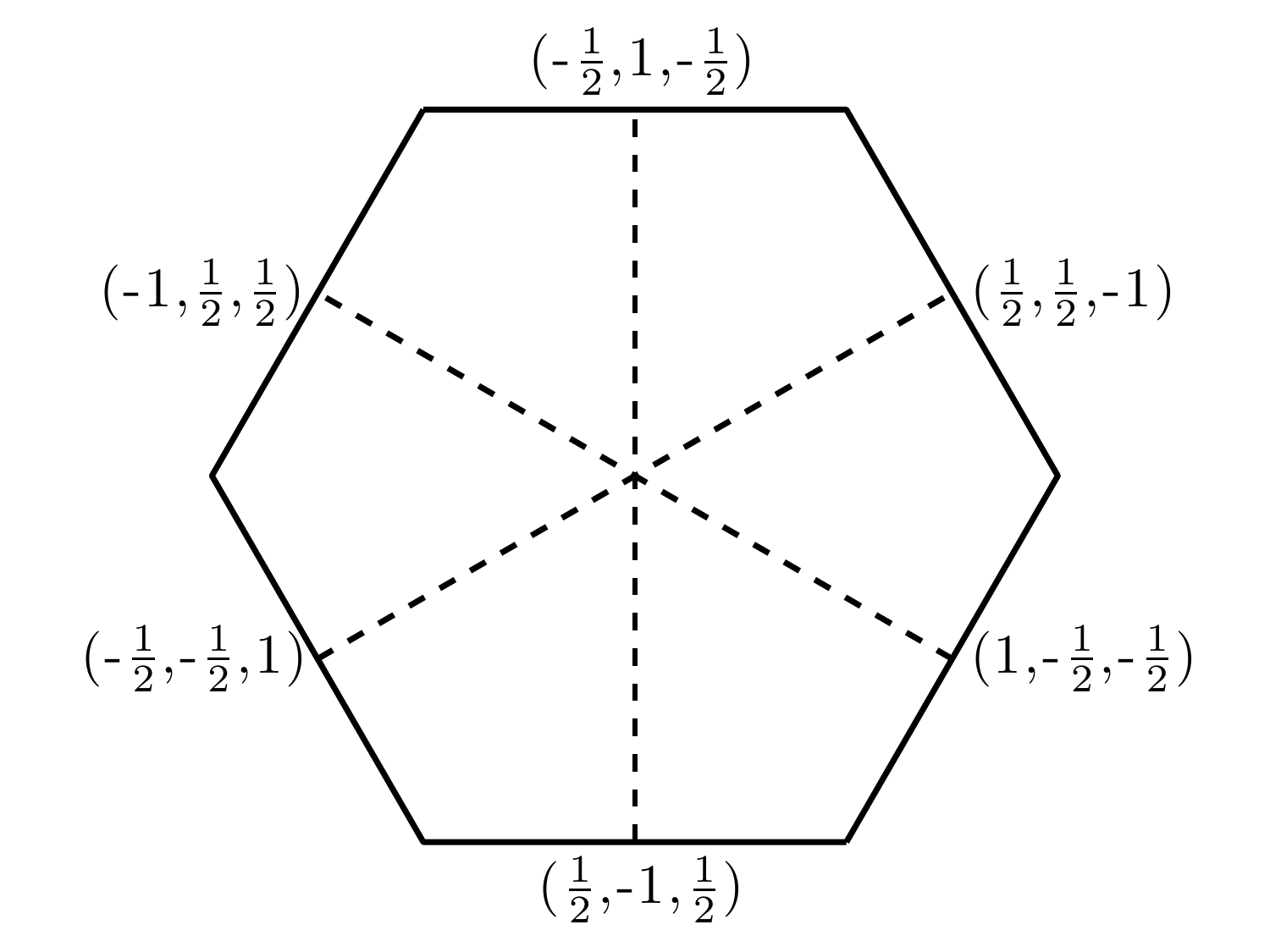}%
\hfill%
\includegraphics[width=0.31\textwidth]{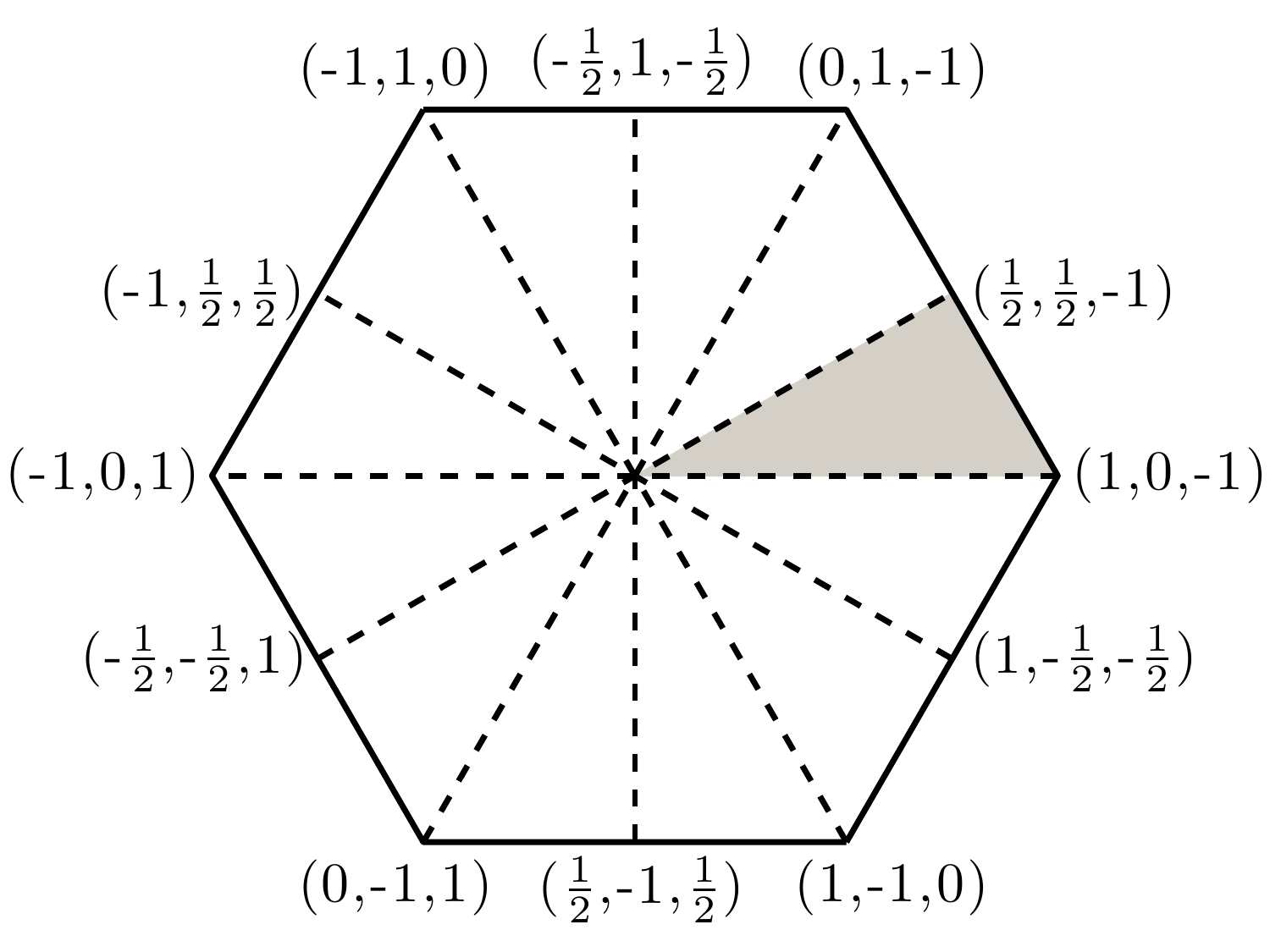}%
\hspace*{\fill}%
\caption{Symmetry under $\A_2$ (left), $\A_2^*$ (center) and $G_2$ (right) in the physical space. The shaded area is the fundamental triangle of $\Omega_A$ under $G_2$.}
\label{hpsym}
\end{figure}
\begin{figure}[htb]
\hfill%
\includegraphics[width=0.31\textwidth]{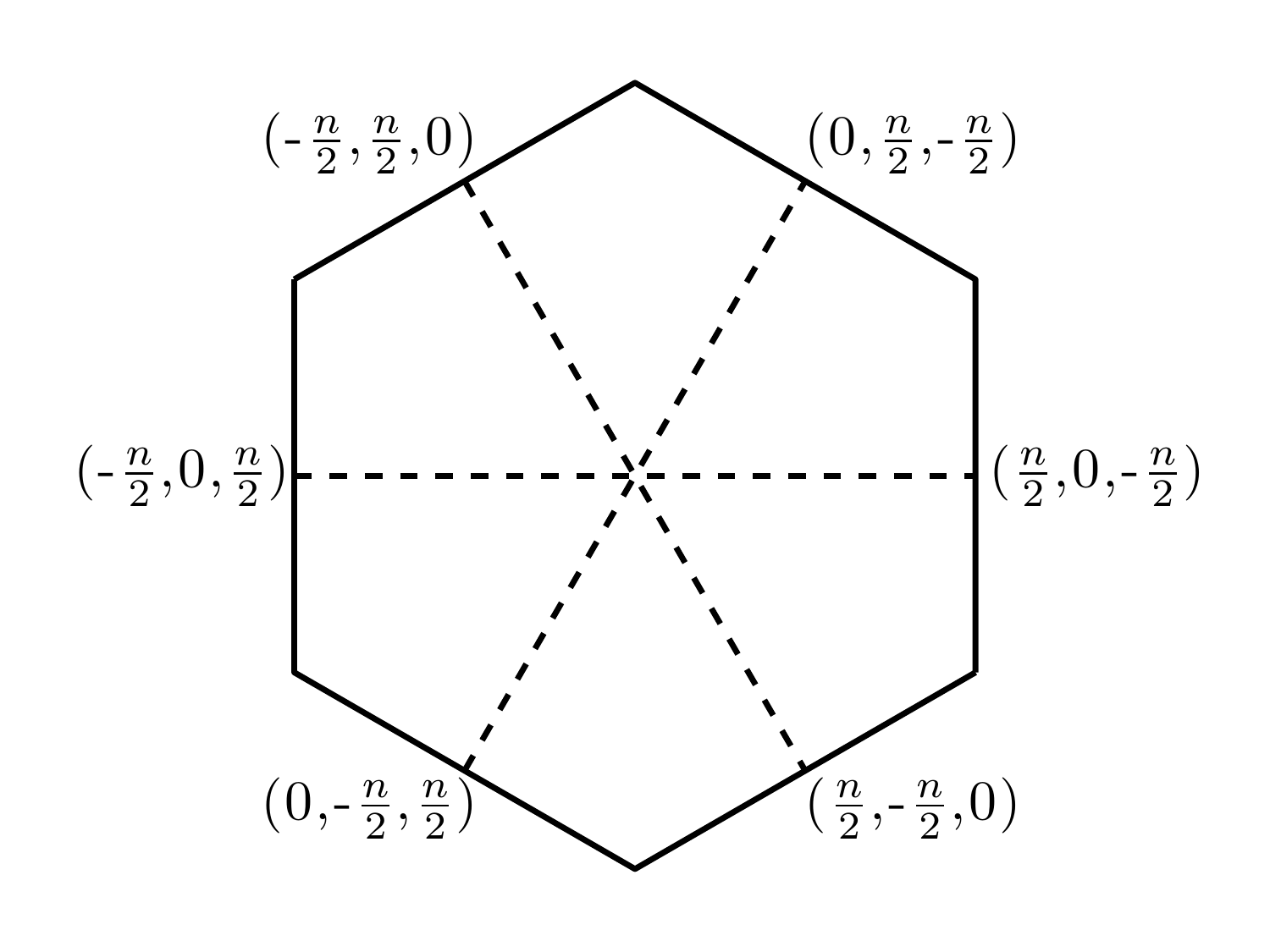}%
\hfill%
\includegraphics[width=0.31\textwidth]{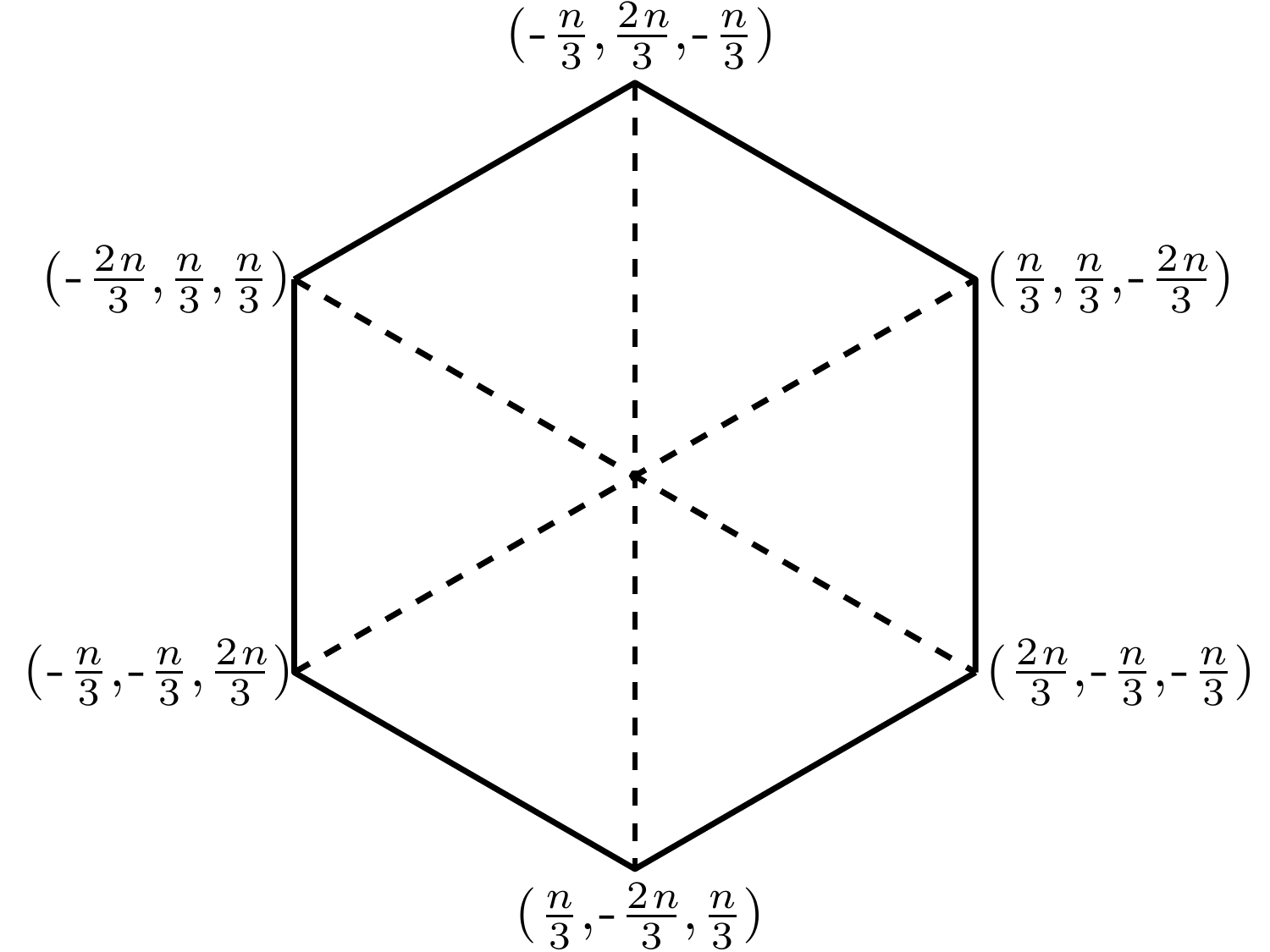}%
\hfill%
\includegraphics[width=0.31\textwidth]{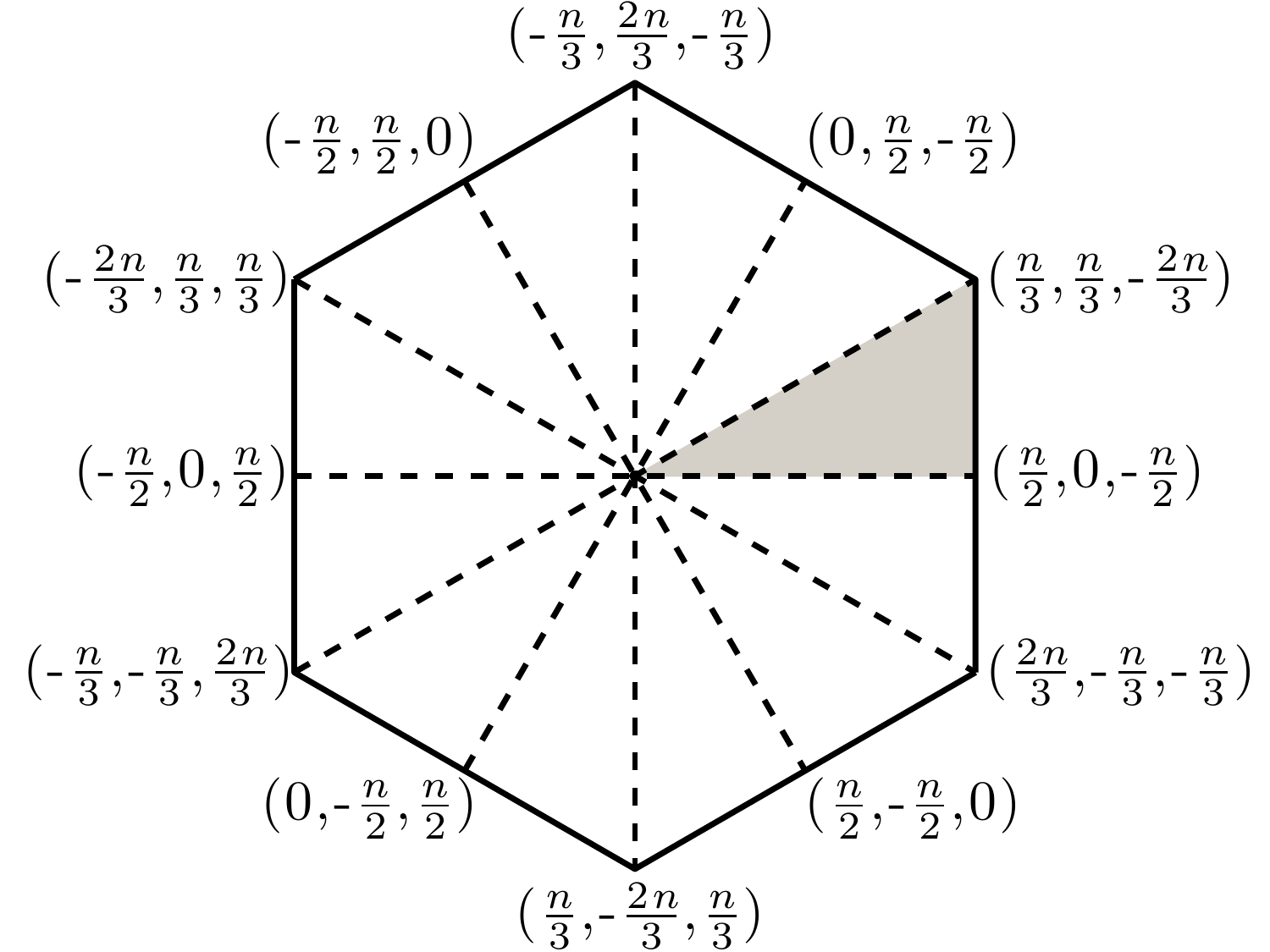}%
\hspace*{\fill}%
\caption{Symmetry under $\A_2$ (left), $\A_2^*$ (center) and $G_2$ (right) in the frequency space. The shaded area is the fundamental triangle of $\Omega_B$ under $G_2$.}
\label{hfsym}
\end{figure}

For $\sigma \in G_2$, the number of
inversion $|\sigma|$ satisf\/ies $|-\sigma|=|\sigma|$. The following lemma can be easily verif\/ied (writing down the
table of $\sigma \sigma^*$ for $\sigma \in A_2$ and $\sigma^* \in A_2^*$ if necessary).

\begin{lem}
Let $f$ be a generic H-periodic function. Then
\begin{gather*}
 \CP^{+}_{\!*} \CP^+ f(\tb)  = \frac{1}{12} \sum_{\sigma\in \mathcal{A}_2}\left(  f(\tb\sigma) +  f(-\tb\sigma)\right)= \frac{1}{12} \sum_{\sigma\in \mathcal{G}_2}  f(\tb\sigma),
\\  \CP^{-}_{\!*} \CP^+ f(\tb)  = \frac{1}{12} \sum_{\sigma\in \mathcal{A}_2}\left(  f(\tb\sigma)-  f(-\tb\sigma)\right)=\frac{1}{12} \sum_{\sigma\in \mathcal{A}_2^*}(-1)^{|\sigma|}\left(  f(\tb\sigma)-  f(-\tb\sigma)\right),
\\  \CP^{+}_{\!*} \CP^- f(\tb)  = \frac{1}{12} \sum_{\sigma\in \mathcal{A}_2} (-1)^{|\sigma|}\left(  f(\tb\sigma)- f(-\tb\sigma)\right)=\frac{1}{12} \sum_{\sigma\in \mathcal{A}_2^*}\left(  f(\tb\sigma)-  f(-\tb\sigma)\right),
\\  \CP^{-}_{\!*} \CP^- f(\tb)  = \frac{1}{12} \sum_{\sigma\in \mathcal{A}_2}(-1)^{|\sigma|}\left(  f(\tb\sigma)+ f(-\tb\sigma)\right)= \frac{1}{12} \sum_{\sigma\in \mathcal{G}_2} (-1)^{|\sigma|} f(\tb\sigma).
\end{gather*}
\end{lem}

For $\phi_\kb(\tb) = \e^{\frac{2 \pi i \kb \cdot \tb}{3}}$, the action of $\CP^+$ and $\CP^-$ on $\phi_\kb$ are called
the {\it generalized cosine} and {\it generalized sine} functions in \cite{LSX}, which are trigonometric functions
given by
\begin{gather}
 \TC_\kb(\tb) :=  \CP^+ \phi_\kb(\tb)  =
 \frac{1}{3} \left[ \e^{\frac{i\pi}{3}(k_1-k_3)(t_1-t_3)}\cos k_2\pi t_2  \right.\notag\\
       \left.
       \hphantom{\TC_\kb(\tb) :=  \CP^+ \phi_\kb(\tb)  =}{}
       +    \e^{\frac{i\pi}{3}(k_1-k_3)(t_2-t_1)}\cos k_2\pi t_3
   +\e^{\frac{i\pi}{3}(k_1-k_3)(t_3-t_2)}\cos k_2\pi t_1\right],  \label{TC_cos}\\
 \TS_\kb(\tb): = \frac{1}{i} \CP^- \phi_\kb(\tb) =
\frac{1}{3} \left[ \e^{\frac{i\pi}{3}(k_1-k_3)(t_1-t_3)}\sin k_2\pi t_2  \right.\notag\\
     \left.
     \hphantom{\TS_\kb(\tb): = \frac{1}{i} \CP^- \phi_\kb(\tb) =}{}
      +    \e^{\frac{i\pi}{3}(k_1-k_3)(t_2-t_1)}\sin k_2\pi t_3
   +\e^{\frac{i\pi}{3}(k_1-k_3)(t_3-t_2)}\sin k_2\pi t_1\right].\label{TC_sin}
   \end{gather}
Because of the symmetry, we only need to consider these functions on the fundamental domain
of the group $\CA_2$, which is one of the equilateral triangles of the regular hexagon. These functions
form a complete orthogonal basis on the equilateral triangle and they are the analogues of the
cosine and sine functions on the equilateral triangle. These generalized cosine and sine functions
are the building blocks of the discrete Fourier analysis on the equilateral triangle and subsequent
analysis of generalized Chebyshev polynomials in~\cite{LSX}.

We now def\/ine the analogue of such functions on~$G_2$. Since the fundamental domain of the
group $G_2$ is the $30^{\degree}$--$60^{\degree}$--$90^{\degree}$ triangle, which is half of the equilateral
triangle, we can relate the new functions to the generalized cosine and sine functions on the latter
domain. There are, however, four families of such functions,  def\/ined as follows:
\begin{gather*}
  \TCC_\kb(\tb) :=  \CP_{\!*}^+ \CP^+\phi_\kb(\tb)  =  \frac{1}{12} \sum_{\sigma\in \mathcal{A}_2}\left(  \phi_{\kb\sigma}(\tb) + \phi_{-\kb\sigma}(\tb)\right)
    =\frac{1}{2} \big(\TC_\kb(\tb)+\TC_{-\kb}(\tb)\big), \\
  \TSC_\kb(\tb) :=\frac1i\CP_{\!*}^- \CP^+\phi_\kb(\tb)= \frac{1}{12i} \sum_{\sigma\in \mathcal{A}_2}\left(  \phi_{\kb\sigma}(\tb) -\phi_{-\kb\sigma}(\tb)\right)  =
    \frac{1}{2i} \big(\TC_\kb(\tb)-\TC_{-\kb}(\tb)\big), \\
  \TCS_\kb(\tb) := \frac1i \CP_{\!*}^+ \CP^-\phi_\kb(\tb)  = \frac{1}{12i} \sum_{\sigma\in \mathcal{A}_2}(-1)^{|\sigma|}\left(  \phi_{\kb\sigma}(\tb) -\phi_{-\kb\sigma}(\tb)\right)
   =\frac{1}{2} \big(\TS_\kb(\tb)-\TS_{-\kb}(\tb)\big), \\
  \TSS_\kb(\tb): = -\CP_{\!*}^- \CP^-\phi_\kb(\tb) = -\frac{1}{12} \sum_{\sigma\in \mathcal{A}_2}(-1)^{|\sigma|}\left(  \phi_{\kb\sigma}(\tb) + \phi_{-\kb\sigma}(\tb)\right)
    =\frac{1}{2i } \big(\TS_\kb(\tb)+\TS_{-\kb}(\tb)\big),
\end{gather*}
where the second and the third equalities follow directly from the def\/inition. We call these functions {\it generalized trigonometric functions}.
As their names indicate, they are of the mixed type of cosine and sine functions.

From \eqref{TC_cos} and \eqref{TC_sin}, we can derive explicit formulas for these functions, which are
\begin{gather}
\TCC_\kb (\tb) = \frac{1}{3} \Big[  \cos {\tfrac{\pi (k_1-k_3)(t_1-t_3)}{3}}\cos \pi k_2 t_2
              + \cos {\tfrac{\pi (k_1-k_3)(t_2-t_1)}{3}}\cos \pi k_2 t_3 \nonumber \\
\hphantom{\TCC_\kb (\tb) =}{}   + \cos {\tfrac{\pi (k_1-k_3)(t_3-t_2)}{3}}\cos \pi k_2 t_1\Big],\label{TCC}
\\
\TSC_\kb(\tb) =  \frac{1}{3} \Big[ \sin {\tfrac{\pi (k_1-k_3)(t_1-t_3)}{3}}\cos \pi k_2 t_2
              + \sin {\tfrac{\pi (k_1-k_3)(t_2-t_1)}{3}}\cos \pi k_2 t_3 \nonumber\\
\hphantom{\TSC_\kb(\tb) =}{}
  + \sin {\tfrac{\pi (k_1-k_3)(t_3-t_2)}{3}}\cos \pi k_2 t_1\Big], \label{TSC}
\\
\TCS_\kb(\tb) =  \frac{1}{3} \Big[ \cos {\tfrac{\pi (k_1-k_3)(t_1-t_3)}{3}}\sin \pi k_2 t_2
              + \cos {\tfrac{\pi (k_1-k_3)(t_2-t_1)}{3}}\sin \pi k_2 t_3 \nonumber\\
\hphantom{\TCS_\kb(\tb) =}{}
 + \cos {\tfrac{\pi (k_1-k_3)(t_3-t_2)}{3}}\sin \pi k_2 t_1\Big],\label{TCS}
\\
\TSS_\kb(\tb) =   \frac{1}{3} \Big[ \sin {\tfrac{\pi (k_1-k_3)(t_1-t_3)}{3}}\sin \pi k_2 t_2
              + \sin {\tfrac{\pi (k_1-k_3)(t_2-t_1)}{3}}\sin \pi k_2 t_3 \nonumber\\
\hphantom{\TSS_\kb(\tb) =}{}
 + \sin {\tfrac{\pi (k_1-k_3)(t_3-t_2)}{3}}\sin \pi k_2 t_1\Big]. \label{TSS}
\end{gather}
In particular, it follows from \eqref{TSC}--\eqref{TSS}  that $\TCS_\kb(\tb)\equiv \TSS_\kb(\tb) \equiv 0$
whenever $\kb$ contains zero component and $\TSC_\kb(\tb)\equiv \TSS_\kb(\tb) \equiv 0$ whenever
$\kb$ contains equal elements.
Similar formulas can be derived from the permutations of $t_1$, $t_2$, $t_3$. In fact, the functions~$\TCC_\kb$ and~$\TSS_\kb$ are invariant
and anti-invariant under~$G_2$, respectively, whereas the functions~$\TCS_\kb$ and~$\TSC_\kb$ are of the mixed type, with the
f\/irst one invariant under~$\CA_2$ and anti-invariant under~$\CA_2^*$ and the second one invariant under~$\CA_2^*$ and
anti-invariant under~$\CA_2$. More precisely, these invariant properties lead to the following identities:
\begin{alignat}{3}
\label{SG2}
&\TCC_\kb(\tb\sigma) = \TCC_\kb(\tb), \qquad \TSS_\kb(\tb\sigma) = (-1)^{|\sigma|}\TSS_\kb(\tb),\qquad &&
\sigma \in \G_2, &\\
\label{SA2-1}
&\TSC_\kb(\tb\sigma) = -\TSC_\kb(-\tb\sigma) = \TSC_\kb(\tb),  && \sigma \in \A_2,&
\\
\label{SA2-2}
&\TCS_\kb(\tb\sigma) = -\TCS_\kb(-\tb\sigma) = (-1)^{|\sigma|}\TCS_\kb(\tb),  &&\sigma \in \A_2, &\\
\label{SA2*-1}
&\TSC_\kb(\tb\sigma) = -\TSC_\kb(-\tb\sigma) = (-1)^{|\sigma|}\TSC_\kb(\tb),  && \sigma \in \A_2^*,& \\
\label{SA2*-2}
&\TCS_\kb(\tb\sigma) = -\TCS_\kb(-\tb\sigma) = \TCS_\kb(\tb),  && \sigma \in \A_2^*.&
\end{alignat}
In particular, it follows from \eqref{TSC}--\eqref{TSS}  that $\TCS_\kb(\tb)\equiv \TSS_\kb(\tb) \equiv 0$
whenever $\kb$ contains zero component and $\TSC_\kb(\tb)\equiv \TSS_\kb(\tb) \equiv 0$ whenever
$\kb$ contains equal elements. Moreover, for any $\kb\in \HH^{\dag}$,  $\TCS_\kb(\tb)=\TSS_\kb(\tb) = 0$
whenever $\tb$ contains zero component and $\TSC_\kb(\tb)=\TSS_\kb(\tb) = 0$ whenever $\tb$
contains equal elements. 

\begin{figure}[htb]
\hfill%
\begin{minipage}{0.4\textwidth}\includegraphics[width=1\textwidth]{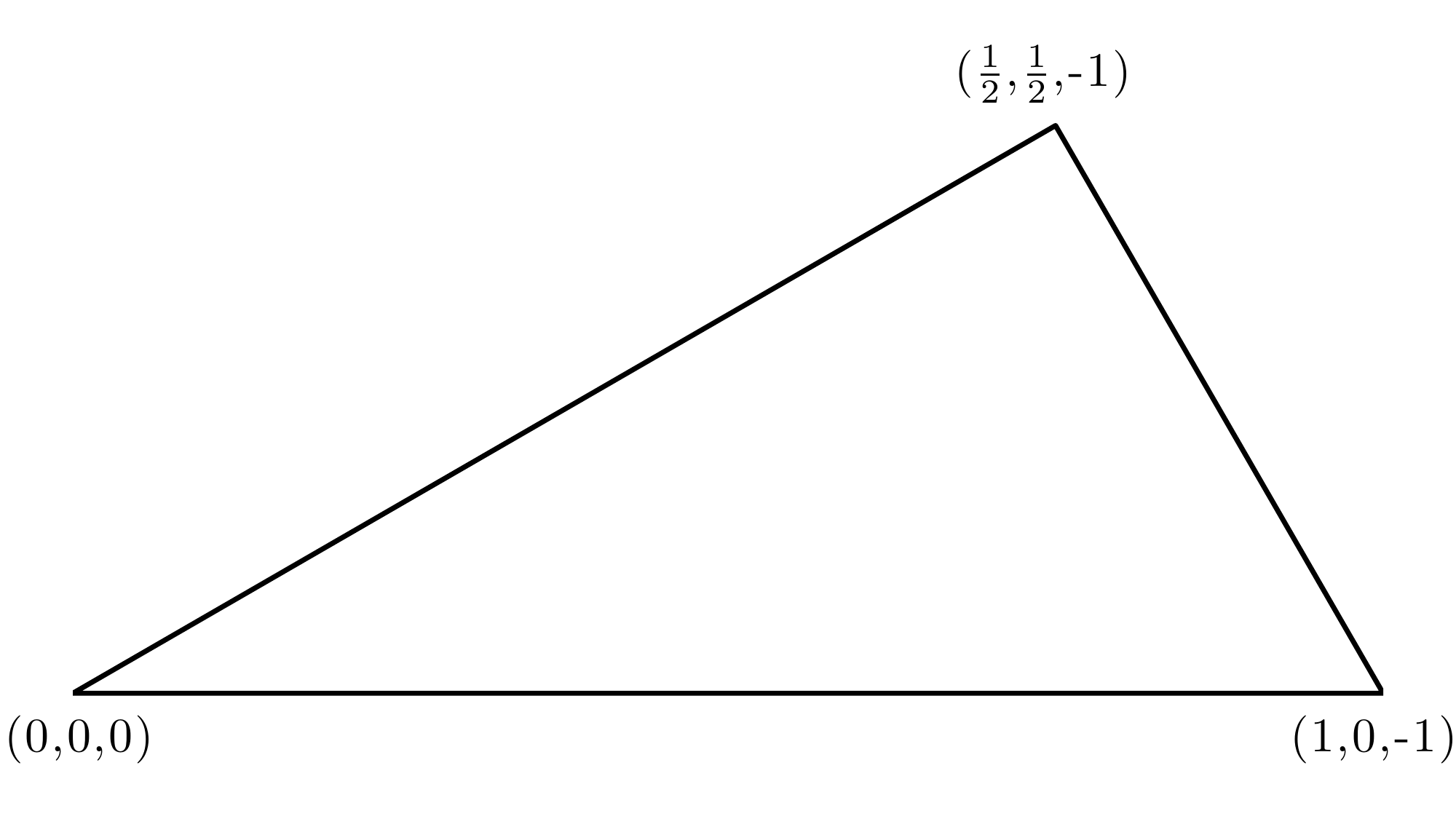}\end{minipage}%
\hfill%
\begin{minipage}{0.4\textwidth}\includegraphics[width=1\textwidth]{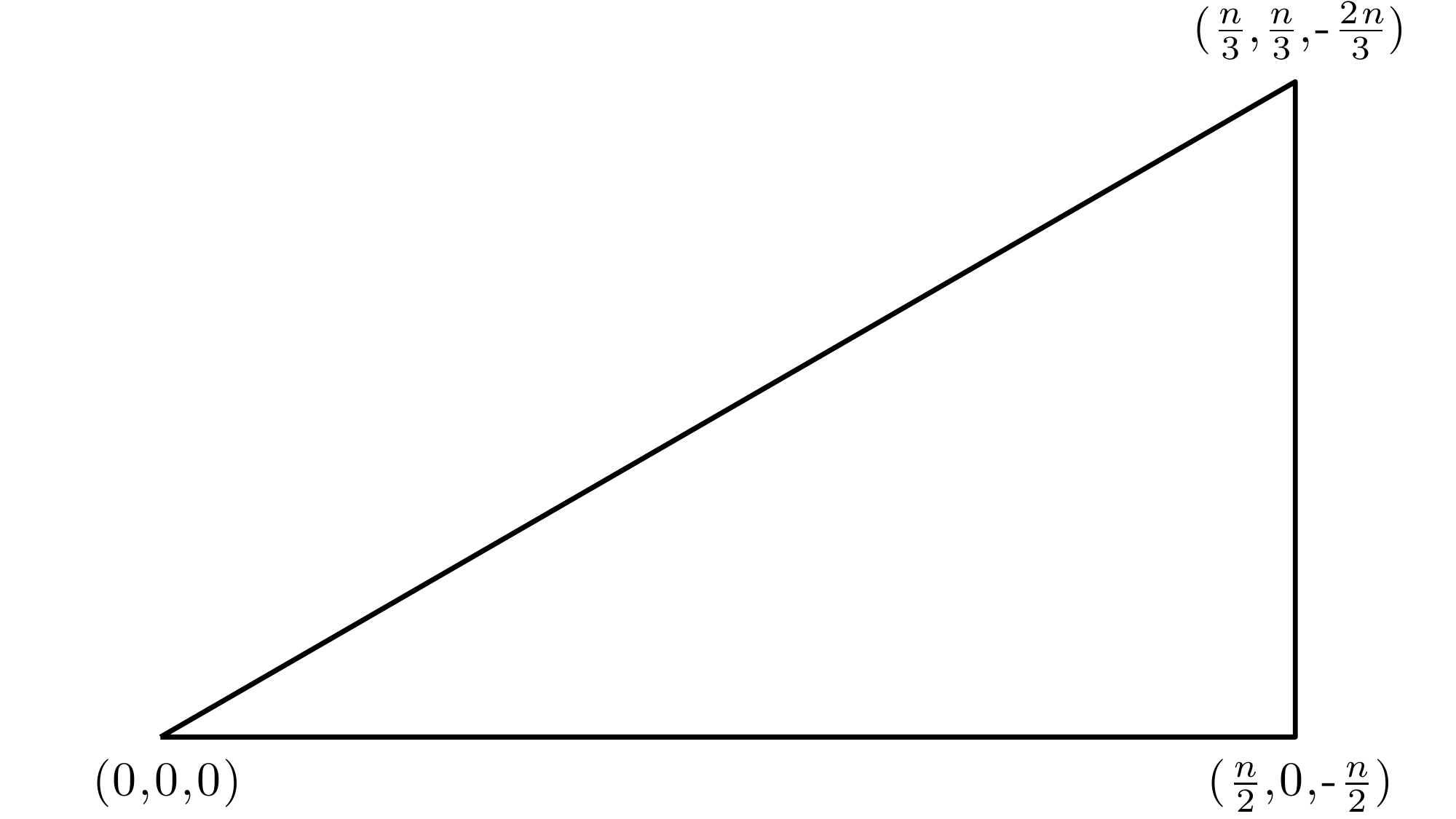}\end{minipage}%
\hspace*{\fill}%
\caption{The fundamental triangles in $\Omega_A$ (left) and $\Omega_B$ (right).}
\label{t306090}
\end{figure}

Because of their invariant properties, we only need to consider these functions on one of the
twelve $30^{{\degree}}$--$60^{{\degree}}$--$90^{{\degree}}$  triangles in the hexagon $\Omega$. We shall
choose the triangle as
\begin{gather} \label{Delta}
   \triangle :=   \{\tb\in \RR_H^3 :    0 \le t_2 \le t_1 \le -t_3\le 1\}.
\end{gather}
The region $\triangle$ and its relative position in the hexagon  are depicted
in Figs.~\ref{t306090} and~\ref{hpsym}.

When $\TCC_\kb$, $\TSC_\kb$, $\TCS_\kb$, $\TSS_\kb$ are restricted to the triangle $\triangle$, we only need to consider
a subset of $\kb \in \HH^{\dag}$ as can be seen by the relations in \eqref{SG2}--\eqref{SA2*-2}. Indeed, we
can restrict $\kb$ to the index sets
\begin{gather} \label{Lambda*}
 \Gamma = \Gamma^{\cc}: = \big\{\kb \in \HH^{\dag}:\; 0 \le k_2 \le  k_1 \big\},   \qquad \Gamma^{\sc} := \big\{\kb \in \HH^{\dag}:\; 0 \le k_2 <  k_1 \big\}, \\
\hphantom{\Gamma =}{} \  \Gamma^{\cs}: = \big\{\kb \in \HH^{\dag}:\; 0 < k_2 \le  k_1 \big\} , \qquad  \Gamma^{\ss} : = \big\{\kb \in \HH^{\dag}:\; 0 < k_2 <  k_1 \big\},
 \end{gather}
respectively, where the notation is self-explanatory; for example, $\Gamma^{\cc}$ is the index set for $\TCC_\kb$.
We def\/ine an inner product on $\triangle$ by
\[
   \langle f, g \rangle_{\triangle} :=
       \frac{1}{|\triangle|}\int_{\triangle} f(\tb)\overline{g(\tb)} d\tb
      = 4 \int_{0}^{\frac12} dt_2 \int_{t_2}^{1-t_2} f(\tb) \overline{g(\tb)} dt_1.
\]
If $f \bar g$ is invariant under the group $G_2$, then it is easy to see that $\la f,g\ra_{\Omega} = \la f, g \ra_{\triangle}$.
Consequently, we can deduce the orthogonality of $\TCC_\kb$, $\TSC_\kb$, $\TCS_\kb$, $\TSS_\kb$ from that of
$\phi_\kb$ on $\Omega$.

\begin{prop}\label{prop:trig-ortho*}
It holds that
\begin{alignat}{3}\label{TCC-ortho}
     &\la\TCC_\kb, \TCC_{\jb} \ra_{\triangle}= \frac{\triangle_{\kb, \jb}}{|\kb G_2|}  =  \triangle_{\kb, \jb} \begin{cases}
       1, & \kb = 0,\\
       \frac{1}{6}, &  k_2 (k_1 -k_2) =0,\, k_1> 0,\\
       \frac{1}{12}, & k_1 >  k_2 > 0,   \end{cases}  \qquad & &\jb,\kb\in \Gamma^{\cc},& \\
     &\la\TSC_\kb, \TSC_{\jb} \ra_{\triangle} = \frac{\triangle_{\kb, \jb}}{|\kb G_2|} =  \triangle_{\kb, \jb} \begin{cases}
       \frac{1}{6}, &  k_2 =0,\\
       \frac{1}{12}, & k_1 >  k_2 > 0,  \end{cases}  && \jb,\kb\in \Gamma^{\sc}, &\\
     &\la\TCS_\kb, \TCS_{\jb} \ra_{\triangle}  = \frac{\triangle_{\kb, \jb}}{|\kb G_2|}=  \triangle_{\kb, \jb} \begin{cases}
       \frac{1}{6}, &  k_1 =k_2 >0,\\
       \frac{1}{12}, & k_1 >  k_2 > 0,   \end{cases}  && \jb,\kb\in \Gamma^{\cs}, & \\
     &\la\TSS_\kb, \TSS_{\jb} \ra_{\triangle} = \frac{\triangle_{\kb, \jb}}{|\kb G_2|} = \tfrac{1}{12} \triangle_{\kb, \jb}
       , && \jb,\kb\in \Gamma^{\ss},&
\end{alignat}
where $\kb G_2=\left\{\kb\sigma : \sigma\in G_2 \right\}$ denotes the orbit of $\kb$ under $G_2$.
\end{prop}

 \subsection[Discrete Fourier analysis on the $30^{{\degree}}$--$60^{{\degree}}$--$90^{{\degree}}$ triangle]{Discrete Fourier analysis on the $\boldsymbol{30^{{\degree}}}$--$\boldsymbol{60^{{\degree}}}$--$\boldsymbol{90^{{\degree}}}$ triangle}

Using the fact that $\TCC_\kb$, $\TSC_{\kb}$ and $\TCS_\kb$, $\TSS_\kb$ are invariant and anti-invariant
under $\A_2$ and that $\TCC_\kb$, $\TCS_{\kb}$ and $\TSC_\kb$, $\TSS_\kb$ are invariant and anti-invariant
under $\A_2^*$, we can deduce a discrete orthogonality for the generalized trignometric functions.
Again, we state the main result in terms of cubature rules. The index set for the nodes of the cubature rule
is given by
\[
\Upsilon_n := \left\{ \jb\in \HH:\; 0\le j_2 \le j_1 \le -j_3 \le n \right\},
\]
which are located inside $n\triangle$ as seen by \eqref{Delta}. The space of invariant functions being
integrated exactly by the cubature rule are indexed by
\begin{gather*}
 \Gamma_n=   \Gamma_n^{\cc} :=  \Gamma \cup \HH^{\dag}_n
      = \big\{\kb\in \HH^{\dag}: \,   0 \le k_2 \le k_1\le k_3+ n\big\},\\
\hphantom{\Gamma_n=}{} \ \Gamma_n^{\sc} := \Gamma^{\sc} \cup \HH^{\dag}_n
      = \big\{\kb\in \HH^{\dag}: \,   0 \le k_2 <  k_1 < k_3+n\big\},\\
\hphantom{\Gamma_n=}{} \ \Gamma_n^{\cs} := \Gamma^{\cs} \cup \HH^{\dag}_n
      = \big\{\kb\in \HH^{\dag}: \,   0 <  k_2 \le k_1  \le k_3 + n\big\},\\
\hphantom{\Gamma_n=}{} \ \Gamma_n^{\ss} : = \Gamma^{\ss} \cup \HH^{\dag}_n
      = \big\{\kb\in \HH^{\dag}: \,   0 < k_2 < k_1< k_3+ n\big\}.
\end{gather*}
Correspondingly, we def\/ine the following subspaces of $\CH_n^{\dag}$,
\begin{gather*}
     \CH^{\cc}_n   := \sspan\{ \TCC_\kb:  \kb \in \Gamma_n^{\cc}\},\qquad
        \CH^{\sc}_n : = \sspan\{ \TSC_\kb:  \kb \in \Gamma_n^{\sc}\},\\
     \CH^{\cs}_n   : = \sspan\{ \TCS_\kb:  \kb \in \Gamma_n^{\cs} \},\qquad
        \CH^{\ss}_n  := \sspan\{ \TSS_\kb:  \kb \in \Gamma_n^{\ss}\}.
\end{gather*}
It is easy to verify that
\begin{gather}
\dim \CH^{\cc}_n = |\Gamma_n^{\cc}| = \tfrac12 \big(3\lfloor\tfrac{n}{3}\rfloor-2n\big)  \big(\lfloor\tfrac{n}{3}\rfloor+1\big)
  -\big(\lfloor\tfrac{n}{2} \rfloor-n-1\big) \big(\lfloor\tfrac{n}{2} \rfloor+1\big),\nonumber\\
 \dim \CH^{\ss}_n=|\Gamma_n^{\ss}| = |\Gamma_{n-6}|,\qquad
 \dim \CH^{\sc}_n= |\Gamma_n^{\sc}|=\dim \CH^{\cs}_n =|\Gamma_n^{\cs}| = |\Gamma_{n-3}|.\label{dimCT}
 \end{gather}

\begin{figure}[htb]
\hfill%
\includegraphics[width=0.31\textwidth]{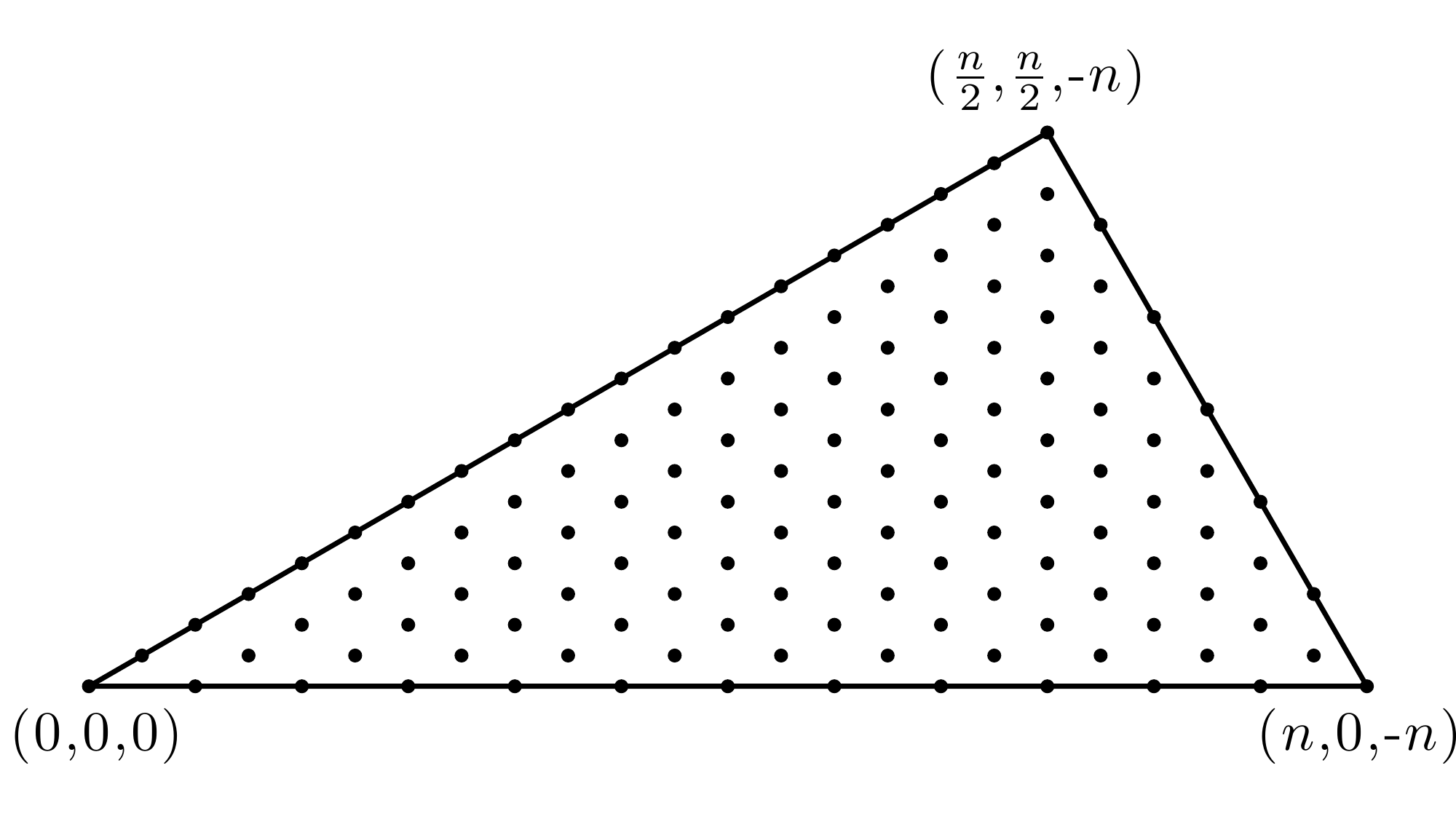}%
\hfill%
\includegraphics[width=0.31\textwidth]{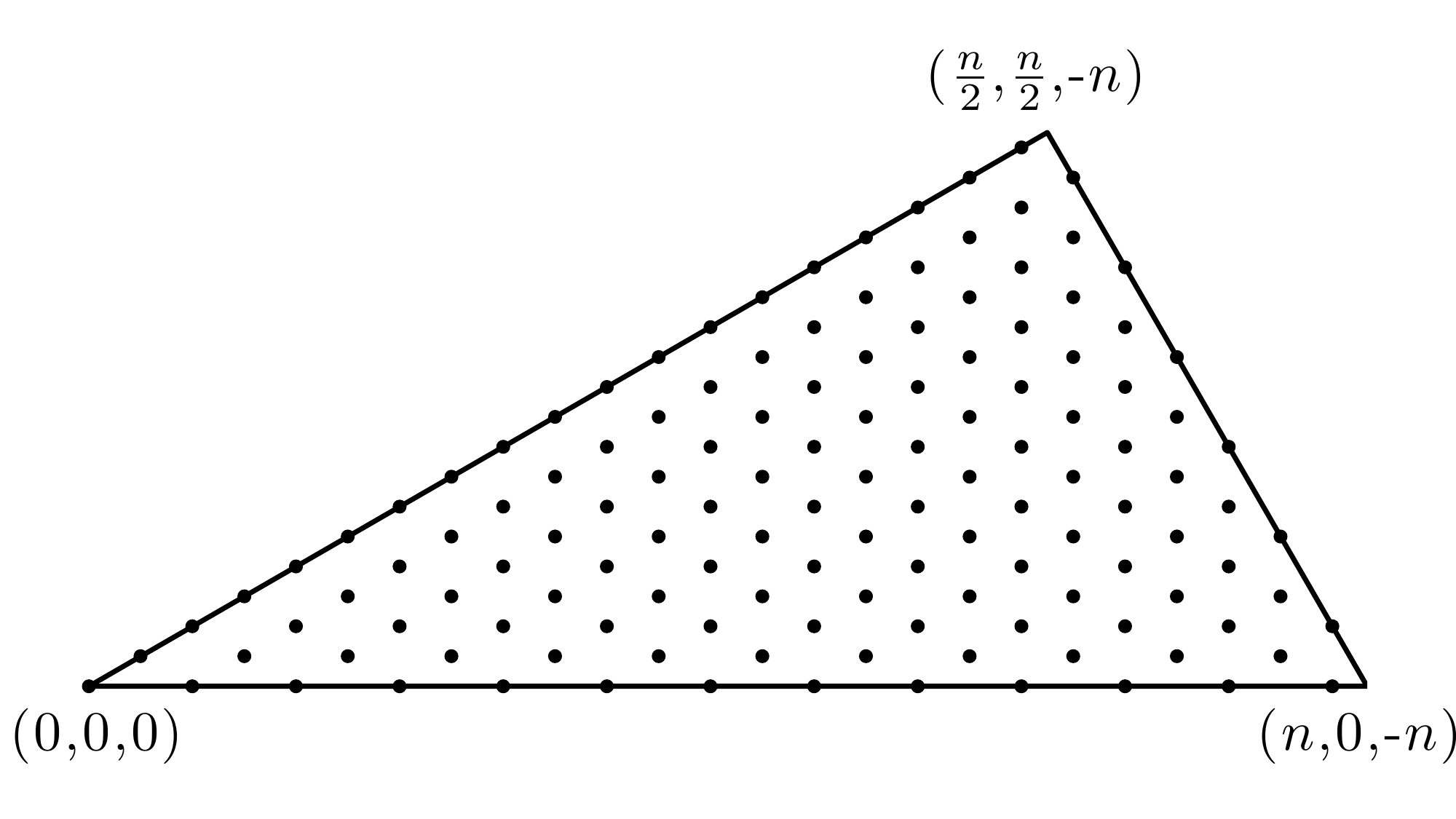}%
\hfill%
\includegraphics[width=0.31\textwidth]{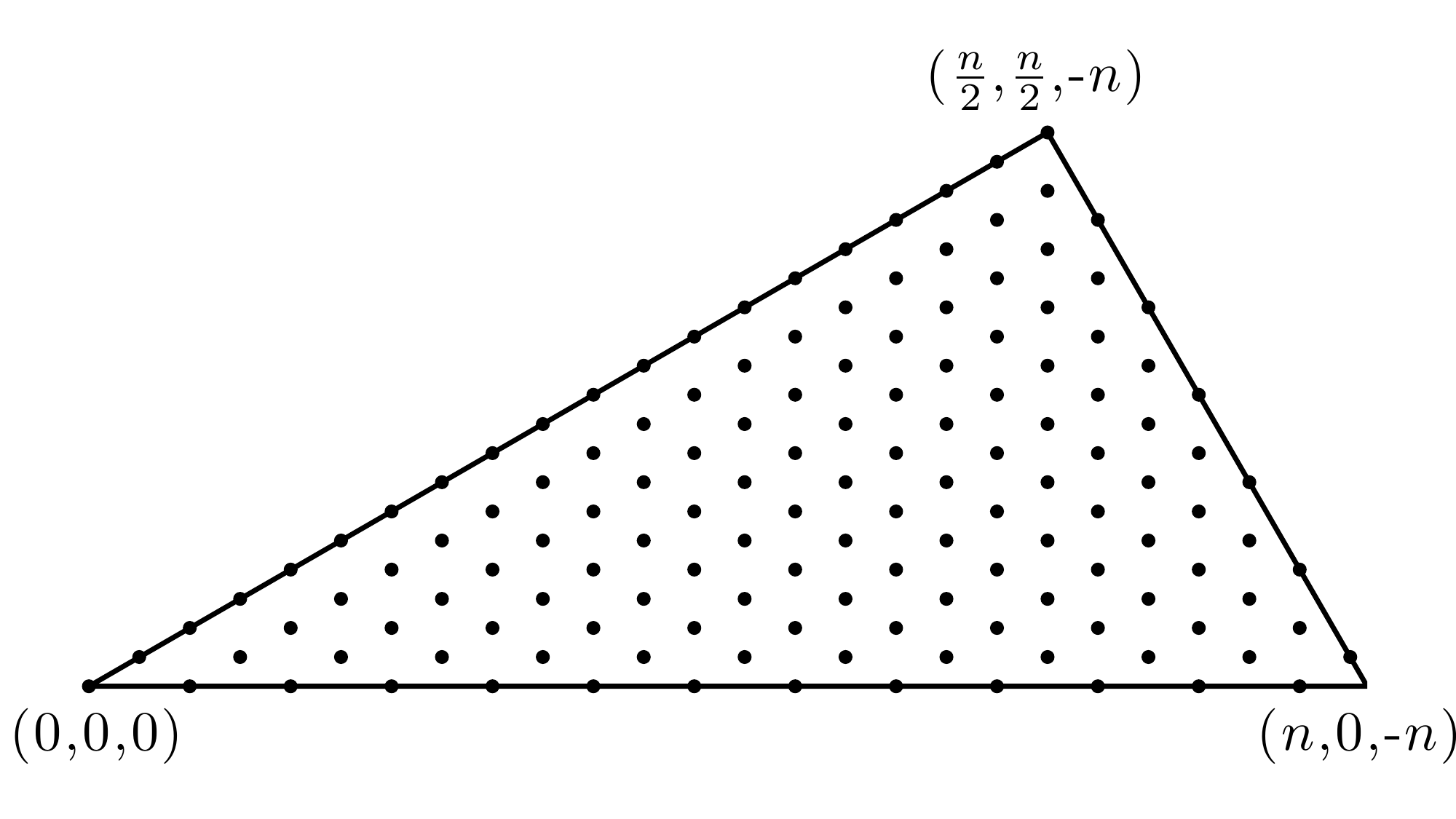}%
\hspace*{\fill}%
\caption{The index set $\Upsilon_n$. $n\equiv 0\pmod{3}$ (left), $n\equiv 1\pmod{3}$ (center) and
$n\equiv 2\pmod{3}$ (right).}
\label{Gp}
\end{figure}

\begin{figure}[htb]
\hfill%
\includegraphics[width=0.4\textwidth]{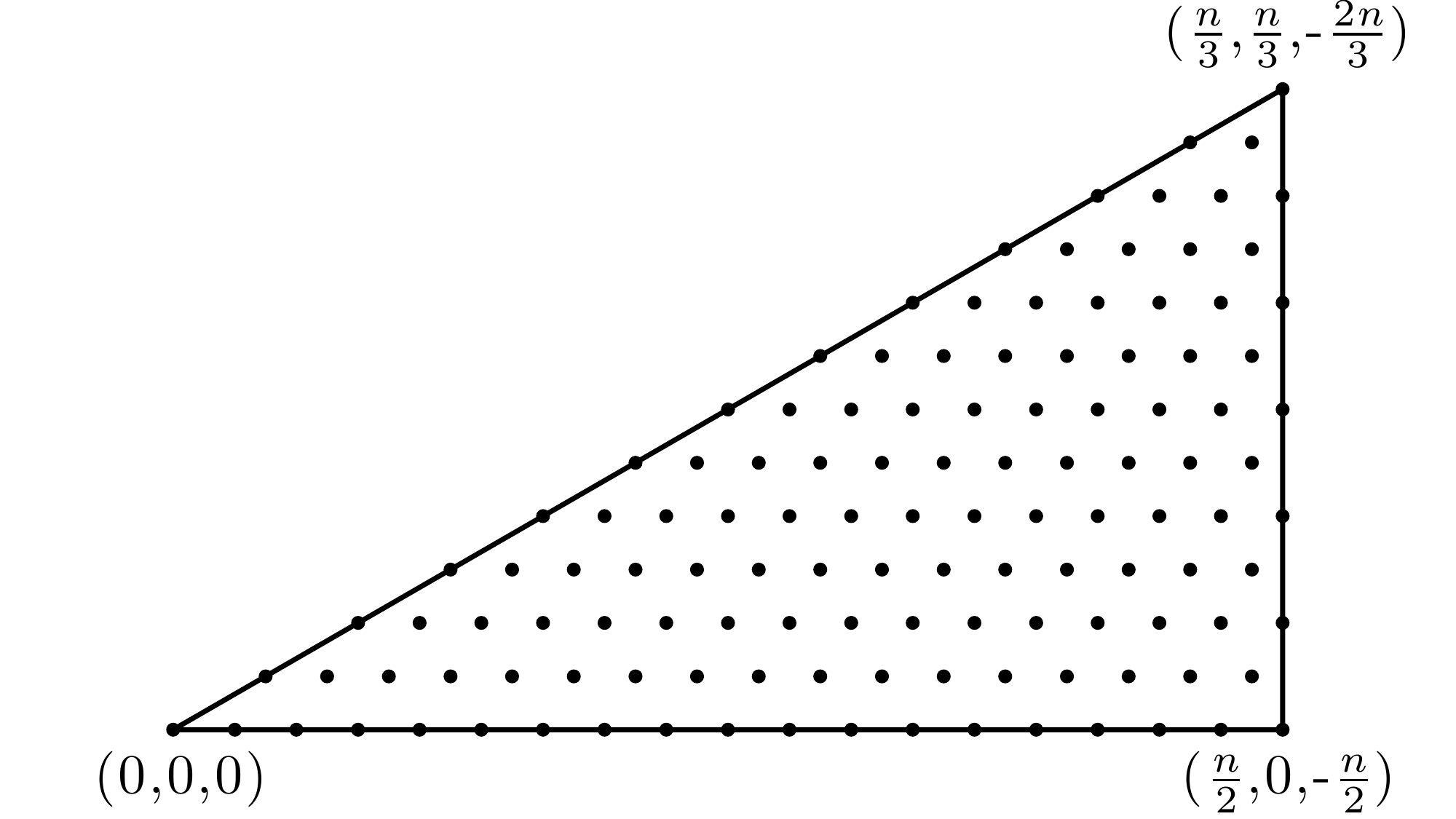}%
\hfill%
\includegraphics[width=0.4\textwidth]{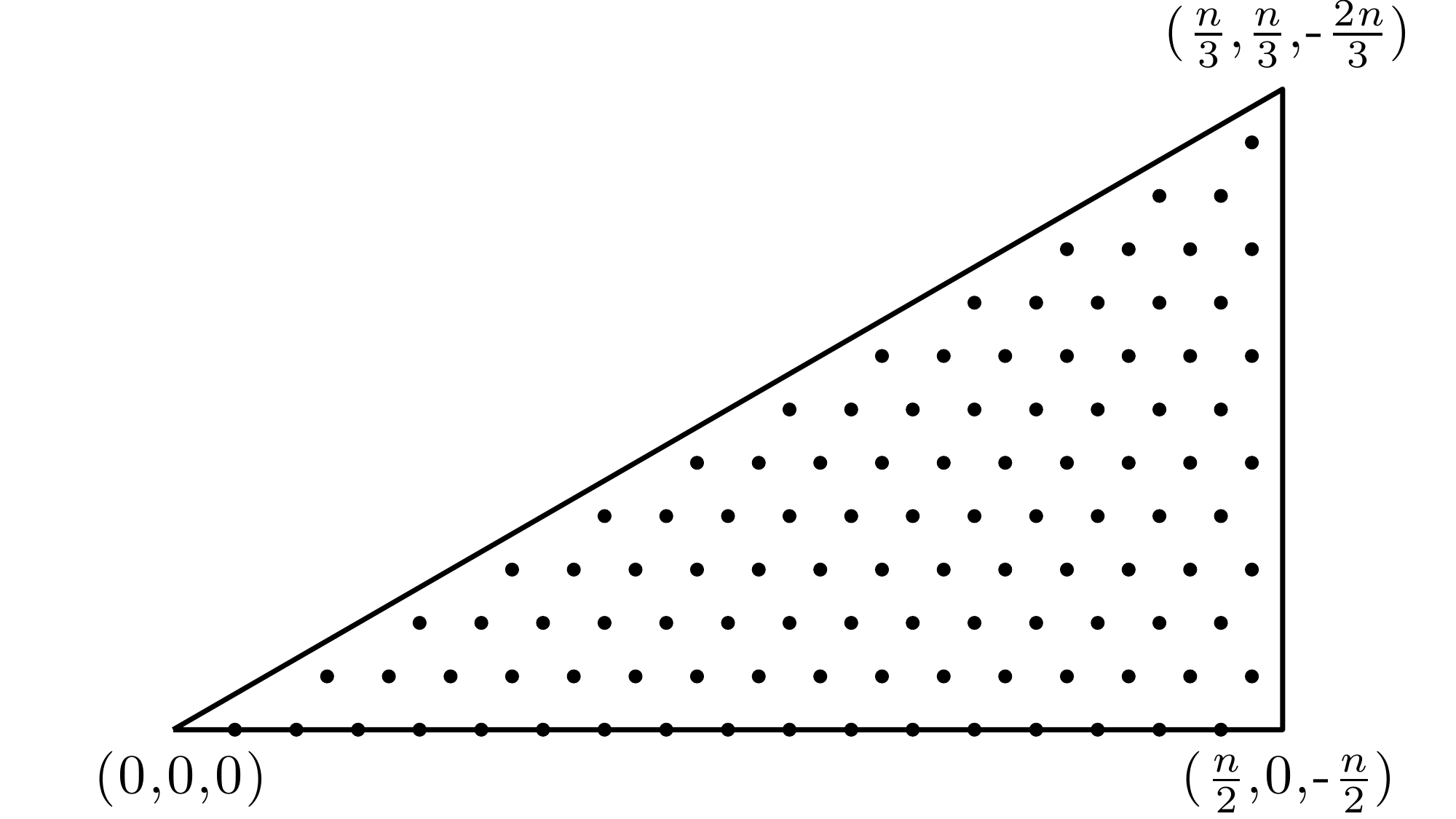}%
\hspace*{\fill}%

\hfill%
(cc)
\hfill%
(sc)
\hspace*{\fill}%

\hfill%
\includegraphics[width=0.4\textwidth]{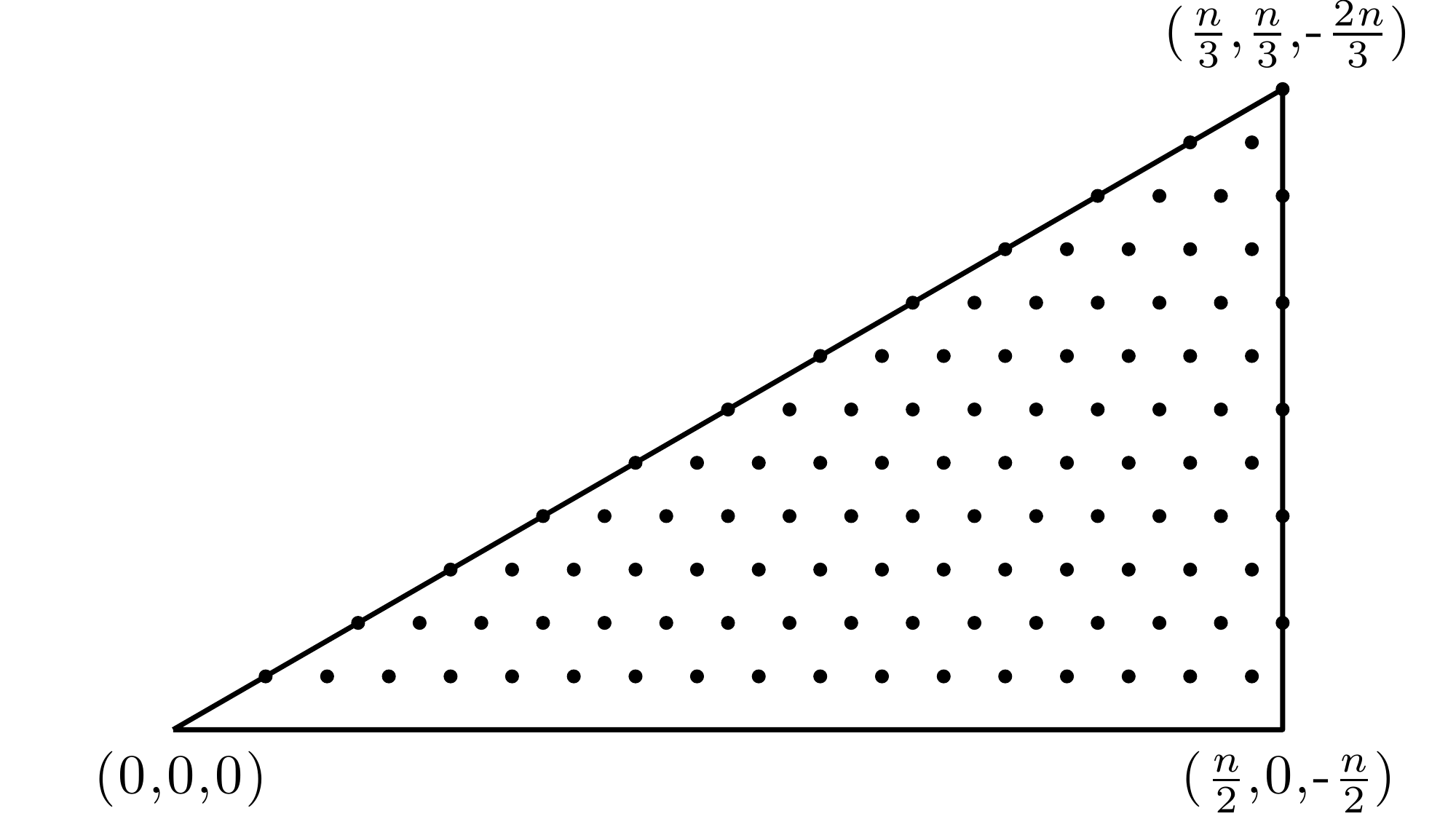}%
\hfill%
\includegraphics[width=0.4\textwidth]{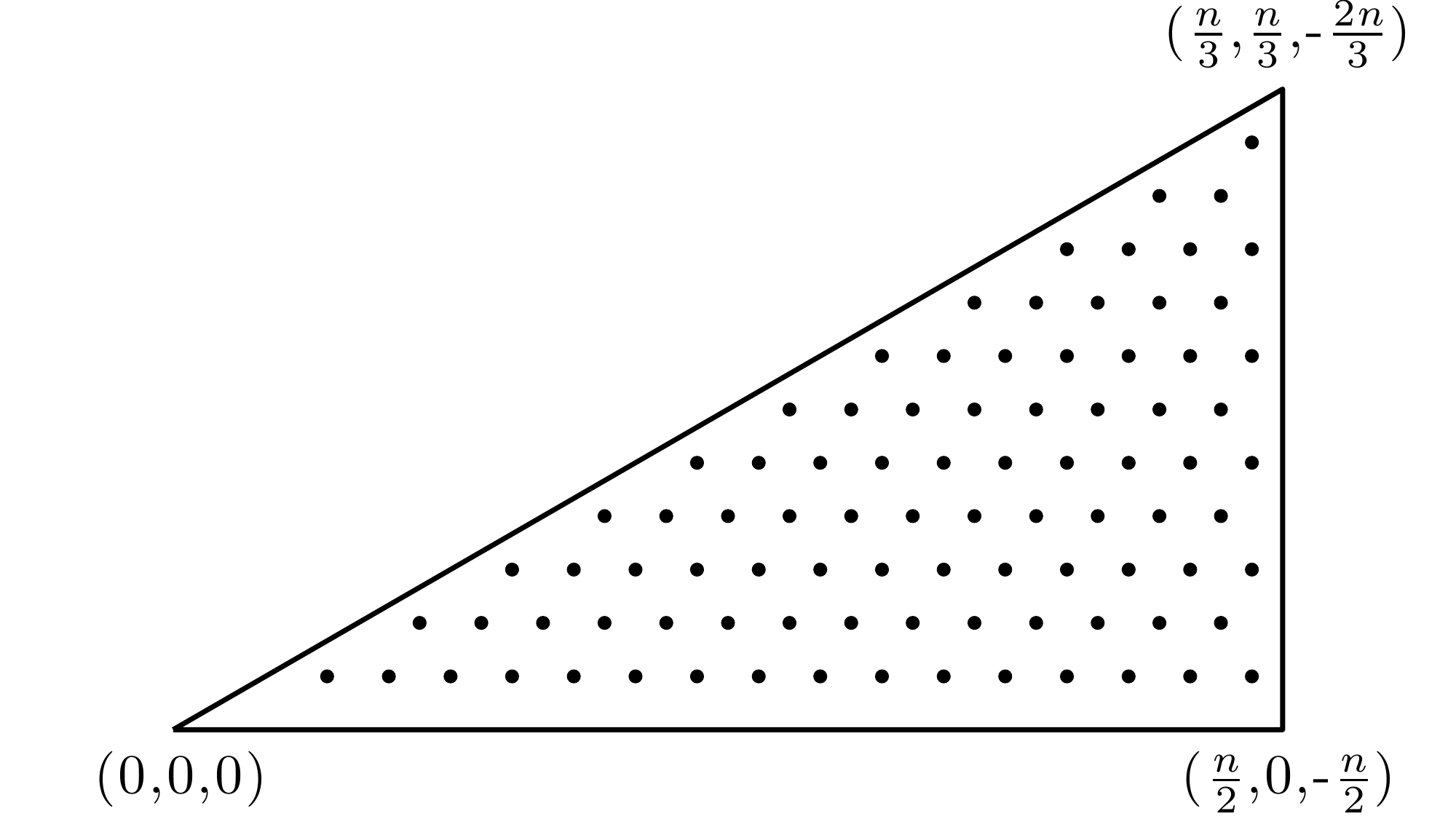}%
\hspace*{\fill}%

\hfill%
(cs)
\hfill%
(ss)
\hspace*{\fill}%

\caption{The index set $\Gamma_n$.}
\label{GfG*}
\end{figure}

%
%
%
%

\begin{thm} The following cubature is exact for all $f\in \CH^{\cc}_{2n-1}$
\begin{gather}
\label{cuba-HHT2}
\frac{1}{|\triangle|} \int_{\triangle} f(\tb)d\tb =\frac{1}{n^2} \sum_{\jb\in\Upsilon_n } \omega_{\jb}^{(n)}
  f\left(\frac{\jb}{n}\right),
\end{gather}
where
\begin{gather*}
\omega_\jb^{(n)} : =  c^{(n)}_{\jb} |\jb G_2|= \left\{\begin{array}{lll}
       12, &  \jb\in \Upsilon^{{\degree}}_n,
        & (\text{interior}),\\[0.2em]
       1, & \jb =\mathbf{0}, & (30^{{\degree}} \text{-vertex}),\\[0.2em]
       2, & \jb =(n,0,-n),   & (60^{{\degree}} \text{-vertex}), \\[0.2em]
       3, & \jb =(\frac{n}{2}, \frac{n}{2}, -n),  &(90^{{\degree}}\text{-vertex}),\\[0.2em]
       6 , &  \text{otherwise} , & (\text{boundaries}).
       \end{array}\right.
\end{gather*}
Moreover, if we define the discrete inner product $\la f,g \ra_{\triangle,n}=\frac1{n^2} \sum\limits_{\jb\in\Upsilon_n } \omega_{\jb}^{(n)} f(\frac{\jb}{n}) \overline{g(\frac{\jb}{n})}$, then
\begin{alignat*}{3}
&\la \TCC_{\jb}, \TCC_{\kb}\ra_{\triangle,n} = \frac{\triangle_{\jb,\kb}}{c^{(n)}_{\wh{\kb}} |\kb G_2|}
= \frac{\triangle_{\jb,\kb}}{\omega^{(n)}_{\wh{\kb}}},  \qquad  && \jb,\kb\in \Gamma_n, & \\
&\la \TSC_{\jb}, \TSC_{\kb}\ra_{\triangle,n} = \frac{\triangle_{\jb,\kb}}{c^{(n)}_{\wh{\kb}} |\kb G_2|}
= \frac{\triangle_{\jb,\kb}}{\omega^{(n)}_{\wh{\kb}}}, \qquad &&  \jb,\kb\in \Gamma^{\sc}_n, & \\
&\la \TCS_{\jb}, \TCS_{\kb}\ra_{\triangle,n} = \frac{\triangle_{\jb,\kb}}{c^{(n)}_{\wh{\kb}} |\kb G_2|}
= \frac{\triangle_{\jb,\kb}}{\omega^{(n)}_{\wh{\kb}}}, \qquad && \jb,\kb\in \Gamma^{\cs}_n, & \\
&\la \TSS_{\jb}, \TSS_{\kb}\ra_{\triangle,n} = \frac{\triangle_{\jb,\kb}}{c^{(n)}_{\wh{\kb}} |\kb G_2|} = \frac{\triangle_{\jb,\kb}}{12}, \qquad && \jb,\kb\in \Gamma^{\ss}_n, &
\end{alignat*}
where $\wh \kb = (k_3-k_2,k_1-k_3,k_2-k_1)$.
\end{thm}

The formula \eqref{cuba-HHT2} is derived from \eqref{cuba-HHD} by using the invariance
of the functions in $\CH^{\cc}_{2n-1}$ and upon writing $\Omega = \big(\bigcup_{\sigma\in G_2}
\{ \tb\sigma :  \tb \in \triangle^{{\degree}}\} \big) \bigcup \big( \bigcup_{\sigma\in G_2}
\{ \tb\sigma :  \tb \in \partial\triangle \} \big)$. The reason that $\wh \kb$ appears goes back to
Proposition \ref{prop:H*H}. As the proof is similar to that  in \cite{LSX}, we shall omit the details.

One may note that the formulation of the result resembles a Gaussian quadrature.
The connection will be discussed in Section~\ref{section6}.

\subsection[Sturm--Liouville eigenvalue problem for the Laplace operator]{Sturm--Liouville eigenvalue problem for the Laplace operator}

Recall the relation \eqref{coordinates} between the coordinates $(x_1,x_2)$ and the homogeneous coordinates $(t_1,t_2,t_3)$.
A quick calculation gives the expression of the Laplace operator in homogeneous coordinates,
\[
     \Delta:=  \frac{\partial^2}{\partial x_1^2} +  \frac{\partial^2}{\partial x_2^2} =
          \frac{1}{2} \left[ \left(\frac{\partial}{\partial t_1}-\frac{\partial}{\partial t_2}\right)^2
        +\left(\frac{\partial}{\partial t_2}-\frac{\partial}{\partial t_3}\right)^2
        +\left(\frac{\partial}{\partial t_3}-\frac{\partial}{\partial t_1}\right)^2  \right].
\]
A further computation shows that $\phi_\kb(\tb) = \e^{\frac{2 \pi i}{3} \kb \cdot \tb}$
are the eigenfunctions of the Laplace operator: for $\kb \in \HH$,
\begin{equation}\label{SL-Laplace}
   \Delta \phi_\kb = - \lambda_\kb \phi_\kb, \qquad \lambda_\kb := \frac{2\pi^2}{9} \left[ (k_1-k_2)^2+(k_2-k_3)^2+(k_3-k_1)^2 \right].
\end{equation}
As a consequence, our generalized trigonometric functions are the solutions of the Sturm--Liouville eigenvalue problem for the
Laplace operator with certain boundary conditions on the $30^{\degree}$--$60^{\degree}$--$90^{\degree}$ triangle. To be more precise, we
denote the three linear segments that are the boundary of this triangle by $B_1$, $B_2$, $B_3$,
\[
    B_1 := \{\tb \in \triangle: t_3 =-1\}, \qquad B_2 := \{\tb \in \triangle: t_2 =0\}, \qquad  B_3 := \{\tb \in \triangle: t_1 =t_2\}.
\]
Let $\frac{\partial}{\partial n}$ denote the partial derivative in the direction of the exterior norm of $\triangle$. Then
\[
    \frac{\partial}{\partial n} \Big \vert_{B_1} = - \frac{\partial }{\partial t_3}, \qquad
     \frac{\partial}{\partial n} \Big \vert_{B_2} =  - \frac{\partial }{\partial t_2}, \qquad
     \frac{\partial}{\partial n} \Big \vert_{B_1} = \frac{\partial }{\partial t_2} -  \frac{\partial }{\partial t_1}.
\]

\begin{thm} \label{thm:SL-Laplace}
The generalized trigonometric functions $\TCC_\kb$, $\TSC_\kb$, $\TCS_\kb$, $\TSS_\kb$
are the eigenfunctions of the Laplace operator, $\Delta u = - \lambda_\kb u$, that satisfy the boundary conditions:
\begin{alignat*}{3}
  & \TCC_\kb:    \  \frac{\partial u}{\partial n} \Big \vert_{B_1 \cup B_2 \cup B_3}  =0, \qquad
  && \TSC_\kb:    \  \frac{\partial u}{\partial n} \Big \vert_{B_1 \cup B_2} =0, \qquad u  \vert_{B_3} = 0, & \\
  &  \TCS_\kb:    \  \frac{\partial u}{\partial n} \Big \vert_{B_3} =0, \qquad u \vert_{B_1 \cup B_2} = 0,     \qquad
      && \TSS_\kb:    \    u \vert_{B_1 \cup B_2 \cup B_3} = 0. &
\end{alignat*}
\end{thm}

\begin{proof}
Since $\lambda_\kb$ is invariant under $G_2$, that is, $\lambda_\kb = \lambda_{\kb \sigma}$, $\forall \, \sigma \in G_2$,
that these functions satisfy $\Delta u = - \lambda_\kb u$ follows directly from their def\/initions. The boundary conditions
can be verif\/ied directly via the equations \eqref{TCC}, \eqref{TSC}, \eqref{TCS} and \eqref{TSS}.
\end{proof}

In particular, $\TCC_\kb$ satisf\/ies the Neumann boundary conditions and $\TSS_\kb$ satisf\/ies the Dirichlet
type boundary conditions.

\subsection{Product formulas for the generalized trigonometric functions}

Below we give a list of identities on the product of the generalized trigonometric functions,
which will be needed in the following section.

\begin{lem}  \label{recurH}
The generalized trigonometric functions satisfy the relations,
\begin{gather}
\label{TCCTCC}
   \TCC_{\jb} \TCC_{\kb} = \frac{1}{12} \sum_{\sigma \in G_2} \TCC_{\kb+\jb\sigma}
   =  \frac{1}{12} \sum_{\sigma\in G_2} \TCC_{\jb+\kb\sigma},\\
\label{TCCTSC}
   \TCC_{\jb} \TSC_{\kb} = \frac{1}{12} \sum_{\sigma \in G_2} \TSC_{\kb+\jb\sigma}
  = \frac1{12} \sum_{\tau \in \A_2^*} (-1)^{\tau} \big(\TSC_{\jb+\kb\tau}- \TSC_{\jb-\kb\tau}\big),\\
\label{TCCTCS}
    \TCC_{\jb} \TCS_{\kb} = \frac{1}{12} \sum_{\sigma \in G_2} \TCS_{\kb+\jb\sigma}
    = \frac1{12} \sum_{\tau \in \A_2^*} \big(\TCS_{\jb+\kb\tau}- \TCS_{\jb-\kb\tau}\big),\\
\label{TCCTSS}
    \TCC_{\jb} \TSS_{\kb} = \frac{1}{12}  \sum_{\sigma \in G_2}  \TSS_{\kb+\jb\sigma}
 =\frac{1}{12}   \sum_{\sigma\in G_2}   (-1)^{|\tau|}\TSS_{\jb+\kb\sigma},\\
    \TSC_{\jb} \TSC_{\kb}
   = -\frac{1}{12}  \sum_{\tau \in \A_2^*}   (-1)^{|\tau|}\big(\TCC_{\kb+\jb\tau} - \TCC_{\kb-\jb\tau}\big)\nonumber\\
 \hphantom{\TSC_{\jb} \TSC_{\kb}}{}
  = -\frac{1}{12}   \sum_{\tau \in \A_2^*}   (-1)^{|\tau|}\big(\TCC_{\jb+\kb\tau} - \TCC_{\jb-\kb\tau}\big), \label{TSCTSC}
\\ \label{TSCTCS}
\TSC_{\jb} \TCS_{\kb}
   = \frac{1}{12} \sum_{\tau \in \A_2^*} (-1)^{|\tau|}\big(\TSS_{\kb+\jb\tau} - \TSS_{\kb-\jb\tau}\big)
   = \frac{1}{12} \sum_{\tau \in \A_2^*}  \big(\TSS_{\jb+\kb\tau} - \TSS_{\jb-\kb\tau}\big), \\
\label{TCSTCS}
   \TCS_{\jb} \TCS_{\kb}
   = -\frac{1}{12} \sum_{\tau \in \A_2^*} \big(\TCC_{\kb+\jb\tau} - \TCC_{\kb-\jb\tau}\big)
   = -\frac{1}{12} \sum_{\tau \in \A_2^*} \big(\TCC_{\jb+\kb\tau} - \TCC_{\jb-\kb\tau}\big),
\\
\label{TSSTSS}
  \TSS_{\jb} \TSS_{\kb}
   = \frac{1}{12} \sum_{\sigma \in G_2}(-1)^{|\sigma|} \TCC_{\kb+\jb\sigma}
   = \frac{1}{12} \sum_{\sigma \in G_2} (-1)^{|\sigma|}\TCC_{\jb+\kb\sigma}.
\end{gather}
Furthermore, the following formulas hold:
\begin{gather}
\label{WT}
    3\TSC_{1,0,-1}(\tb) \TCS_{1,1,-2}(\tb) = \TSS_{2,1,-3}(\tb),\\
\label{WT1}
    [\TSC_{1,0,-1}(\tb)]^2 =  \frac{1}{3} \big[ 1 + 2\TCC_{1,1,-2}  \big]-[\TCC_{1,0,-1}]^2,\\
\label{WT2}
     [\TCS_{1,1,-2}]^2 +[\TCC_{1,1,-2}]^2 =  \frac{1}{3} \big[ 1 + 2\TCC_{3,0,-3}  \big],\\
\label{WT3}
    [\TCC_{1,0,-1}]^3  =  \frac{1}{36} \TCC_{3,0,-3} + \frac{1}{4} \TCC_{1,0,-1}  +  \frac{1}{6} \TCC_{1,1,-2} + \frac1{18} +    \frac1{2}\TCC_{1,1,-2} \TCC_{1,0,-1}.
 \end{gather}
\end{lem}

\begin{proof}
For \eqref{TCCTCC}--\eqref{TSSTSS}, we only prove \eqref{TSCTCS}. Other identities can be proved
similarly.  By the def\/inition of the generalized trigonometric functions,
\begin{gather*}
 \TSC_{\jb}  \TCS_{\kb} = \frac1{12i} \sum_{\sigma\in \A_2^*}(-1)^{|\sigma|}
   \big(\phi_{\jb \sigma} -\phi_{-\jb \sigma} \big)
\times  \frac1{12i} \sum_{\tau \in \A_2^*} \big(\phi_{\kb \tau} - \phi_{-\kb \tau} \big)  \\
\hphantom{\TSC_{\jb}  \TCS_{\kb}}{}
=  -\frac{1}{12^2}  \sum_{\tau \in \A_2^*}  (-1)^{|\tau|}   \sum_{\sigma\in \A_2^*}   (-1)^{|\sigma\tau^{-1}|} \big[ \phi_{(\kb +\jb \sigma \tau^{-1}) \tau}  + \phi_{-(\kb +\jb \sigma \tau^{-1}) \tau}\\
\hphantom{\TSC_{\jb}  \TCS_{\kb} =}{}
   - \phi_{(\kb -\jb \sigma \tau^{-1}) \tau}  - \phi_{-(\kb -\jb \sigma \tau^{-1}) \tau} \big]
\end{gather*}
upon using the relation $(-1)^{|\tau| + |\sigma \tau^{-1}|} = (-1)^{|\sigma|}$, consequently,
\begin{gather*}
 \TSC_{\jb} \TCS_{\kb} =  -\frac{1}{12^2} \sum_{\sigma\in \A_2^*} (-1)^{|\sigma|} \sum_{\tau \in \A_2^*} (-1)^{|\tau|} \big[ \phi_{(\kb +\jb \sigma) \tau}
  + \phi_{-(\kb +\jb \sigma) \tau} - \phi_{(\kb -\jb \sigma) \tau}  - \phi_{-(\kb -\jb \sigma) \tau} \big] \\
\hphantom{\TSC_{\jb} \TCS_{\kb}}{}
=  \frac{1}{12} \sum_{\sigma\in \A_2^*} (-1)^{|\sigma|} \big( \TSS_{\kb +\jb \sigma} -\TSS_{\kb -\jb \sigma}\big),
\end{gather*}
proving the f\/irst equality in \eqref{TSCTCS}. Further by \eqref{SG2},{\samepage
\begin{gather*}
\frac{1}{12} \sum_{\sigma\in \A_2^*}  (-1)^{|\sigma|} \big( \TSS_{\kb +\jb \sigma} -\TSS_{\kb -\jb \sigma}\big)
=\frac{1}{12} \sum_{\sigma\in \A_2^*} \big( \TSS_{\kb \sigma^{-1} +\jb} -\TSS_{\kb \sigma^{-1} -\jb}\big)
\\
\qquad{} =  \frac{1}{12} \sum_{\sigma\in \A_2^*} \big( \TSS_{\jb+\kb \sigma } -\TSS_{\kb\sigma-\jb }
\big) =\frac{1}{12} \sum_{\sigma\in \A_2^*} \big( \TSS_{\jb+\kb \sigma } -\TSS_{\jb-\kb \sigma}\big),
\end{gather*}
since $\TSS_{\jb} = \TSS_{- \jb}$ by \eqref{TCC}. This completes the proof of \eqref{TSCTCS}.}

We now prove the relations \eqref{WT}--\eqref{WT3}. By \eqref{TSCTCS},
\begin{gather*}
\TCS_{1,1,-2}(\tb) \TSC_{1,0,-1}(\tb) =
\frac{1}{6} \Big[ \big(\TSS_{2,1,-3}(\tb) - \TSS_{0,-1,1}(\tb) \big)+ \big(\TSS_{-1,1,0}(\tb) - \TSS_{3,-1,-2}(\tb)\big) \\
\hphantom{\TCS_{1,1,-2}(\tb) \TSC_{1,0,-1}(\tb) =}{}
 + \big(\TSS_{2,-2,0}(\tb) - \TSS_{0,2,-2}(\tb)\big) \Big] = \frac1{3}\TSS_{2,1,-3}(\tb),
\end{gather*}
which proves \eqref{WT}.  By \eqref{TSCTSC} and \eqref{TCCTCC}, we have
\begin{gather*}
  [\TSC_{1,0,-1}]^2+[\TCC_{1,0,-1}]^2 =   -\frac{1}{6} \big[ \TCC_{2,0,-2} - 1
      + 2\TCC_{1,-1,0} - 2\TCC_{1,1,-2}  \big]\\
\hphantom{[\TSC_{1,0,-1}]^2}{}   + \frac{1}{6} \big[ \TCC_{2,0,-2} + 1+ 2\TCC_{1,-1,0} + 2\TCC_{1,1,-2}  \big]
 =   \frac{1}{3} \big[ 1 + 2\TCC_{1,1,-2}  \big],
\end{gather*}
which is \eqref{WT1}. Next, from \eqref{TCSTCS} and \eqref{TCCTCC} we deduce that
\begin{gather*}
 [\TCS_{1,1,-2}]^2 +[\TCC_{1,1,-2}]^2 = -\frac{1}{6} \big[ \TCC_{2,2,-4} - 1+ 2\TCC_{1,1,-2} - 2\TCC_{3,0,-3} \big] \\
 \hphantom{[\TCS_{1,1,-2}]^2}{} +\frac{1}{6} \big[ \TCC_{2,2,-4} + 1 + 2\TCC_{1,1,-2} + 2\TCC_{3,0,-3}  \big]
 =  \frac{1}{3} \big[ 1 + 2\TCC_{3,0,-3}  \big],
 \end{gather*}
which is \eqref{WT2}. Finally, the identity \eqref{WT3} follows from a successive use of \eqref{TCCTCC}.
The proof is completed.
\end{proof}

\section{Generalized Chebyshev polynomials}\label{section4}

In \cite{LSX}, the generalized cosine and sine functions $\TC_\kb$ and $\TS_\kb$ are shown to be
polynomials under a~change of variables, which are analogues of Chebyshev polynomials of the f\/irst
and the second kind, respectively, in two variables. These polynomials, f\/irst studied in \cite{K, K2},
are orthogonal polynomials on the region bounded by the hypocycloid and they enjoy a remarkable
property on its common zeros, which yields a rare example of the Gaussian cubature rule.

In this section, we consider analogous polynomials related to our new generalized trigonometric
functions, which has a structure dif\/ferent from those related to $\TC_\kb$ and $\TS_\kb$.

The classical Chebyshev polynomials, $T_n(x)$, are obtained from the trigonometric functions $\cos n \theta$
by setting $x = \cos \theta$, the lowest degree nontrivial trigonometric function. In analogy, we make a change
of variables based on the f\/irst two nontrivial generalized cosine functions:
\begin{gather}
x = x(\tb):=\TCC_{1,0,-1}(\tb) = \frac13\left(\cos\tfrac{2\pi (t_1-t_2)}{3}+\cos\tfrac{2\pi (2t_1+t_2)}{3}+\cos\tfrac{2\pi (2t_2+t_1)}{3}  \right),\nonumber\\
y = y(\tb):=\TCC_{1,1,-2}(\tb) =\frac13\left(\cos 2\pi t_1+\cos 2\pi t_2+\cos 2\pi (t_1+t_2)  \right).\label{xy}
\end{gather}
If we change variables $(t_1, t_2) \mapsto (x, y)$, then the region $\triangle$ is mapped onto the region
$\triangle^*$ bounded by two hypocycloids,
\begin{gather} \label{Delta*}
\triangle^*
   = \left\{(x,y): \, \big(1+2y-3x^2\big)\big(24x^3-y^2-12xy-6x-4y-1\big)\ge 0\right\}.
\end{gather}

\begin{figure}[htb]
\hfill%
\begin{minipage}{0.4\textwidth}\includegraphics[width=1\textwidth]{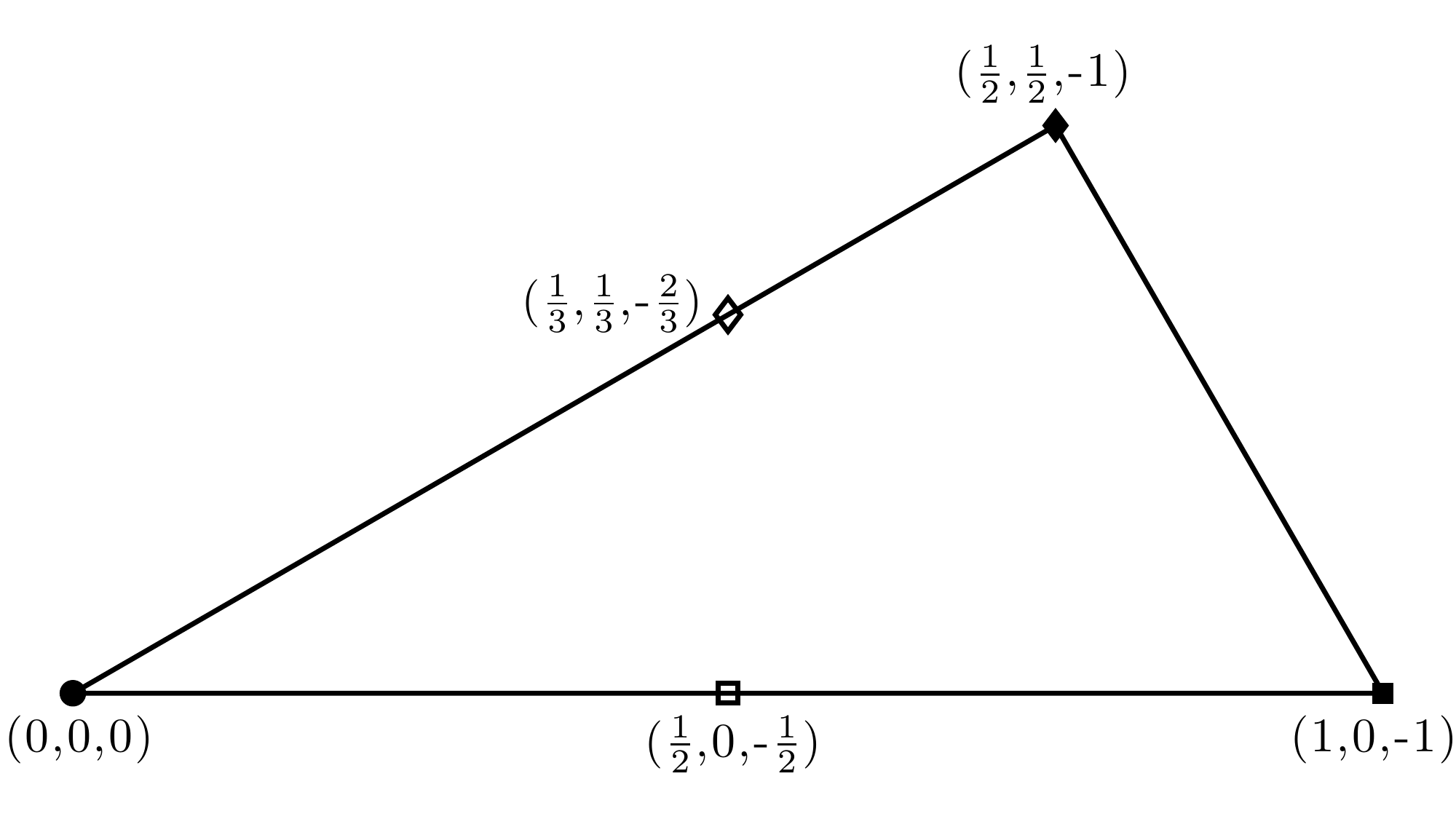}\end{minipage}%
\hfill%
\begin{minipage}{0.4\textwidth}\includegraphics[width=1\textwidth]{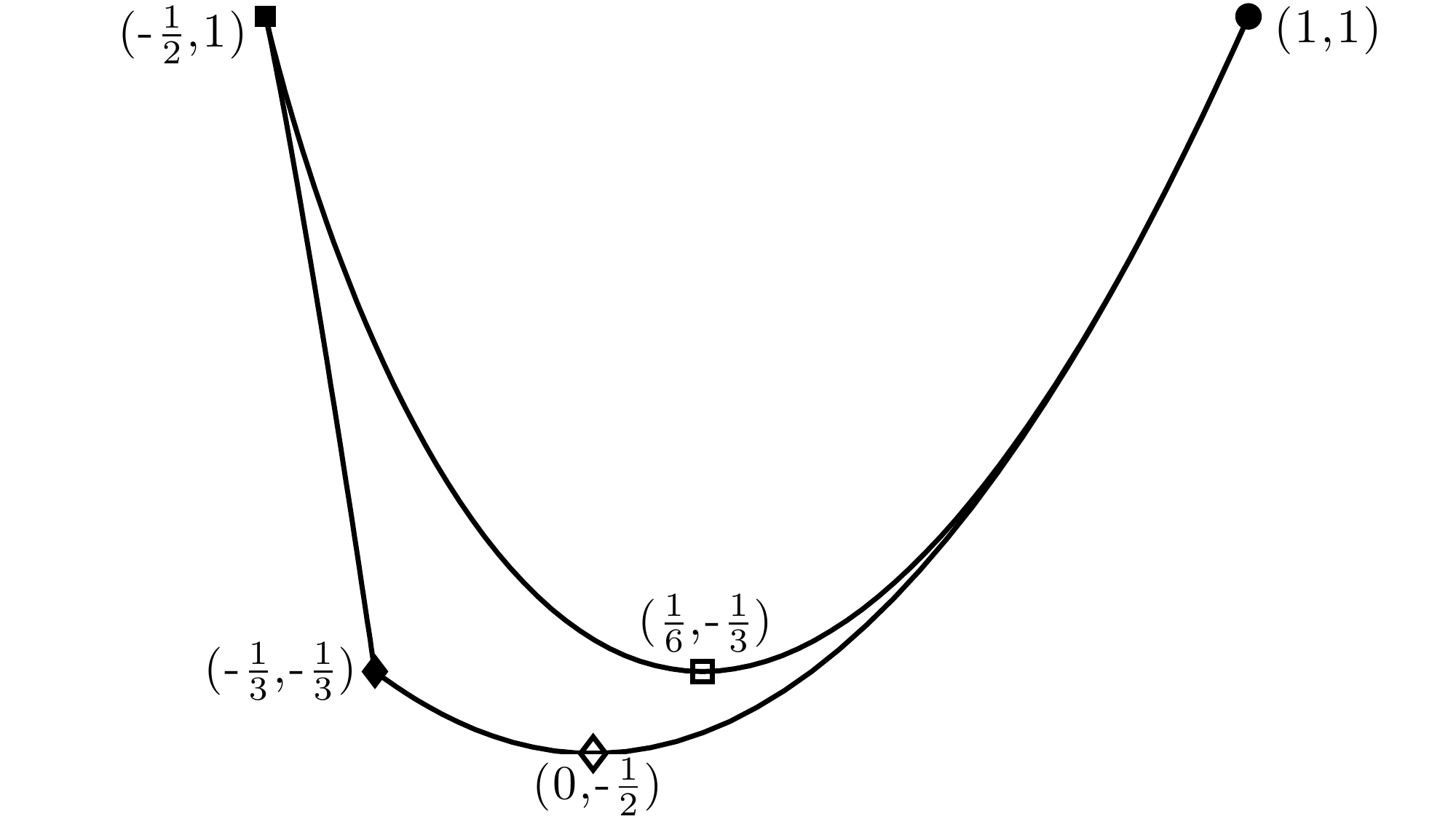}\end{minipage}%
\hspace*{\fill}
\caption{The region $\Delta^*$  (right) bounded by two hypocycloids, which is  mapped from the triangle $\Delta$ (left).}
\end{figure}

The curve that def\/ined the boundary of the domain $\Delta^*$ satisf\/ies the following relation:

\begin{lem} \label{lem:Jacobian}
Let $F(x,y):=  (1+2y-3x^2)(24x^3-y^2-12xy-6x-4y-1)$. Then, in homogeneous coordinates,
\begin{gather} \label{hypocycloid}
    F(x,y)  = 3 \left[\TSC_{1,0,-1}(\tb) \right]^2 \left[\TCS_{1,1,-2}(\tb) \right]^2
               =  \frac{1}{3} \left[\TSS_{2,1,-3}(\tb) \right]^2.
\end{gather}
Furthermore, let $J(x,y)$ be the Jacobian of the changing of variable \eqref{xy}; then
\begin{gather}
    J(x,y) =    \frac{64 \pi^2 }{27} \sin \pi t_1\sin \pi t_2 \sin \pi (t_1+t_2) \sin \frac{ \pi(t_1-t_2)}{3}
                         \sin \frac{ \pi(t_1+2t_2)}{3} \sin \frac{ \pi(2t_1+t_2)}{3} \notag\\
\hphantom{J(x,y)}{}  =   \frac{4  \pi^2}{3}  \TSC_{1,0,-1}(\tb)  \TCS_{1,1,-2}(\tb). \label{Jacobian}
\end{gather}
\end{lem}

\begin{proof}
Under the change of variables \eqref{xy}, by \eqref{WT1},  \eqref{WT2} and \eqref{WT3}, it follows that
\begin{gather}
      \left[\TSC_{1,0,-1}(\tb)\right]^2 = \frac{1}{3}\left(1+2y-3x^2\right), \nonumber\\
\left[\TCS_{1,1,-2}(\tb) \right]^2 = 24x^3-y^2-12xy-6x-4y-1,\label{2equality}
\end{gather}
from which the f\/irst equality in \eqref{hypocycloid} follows, whereas the second one follows from \eqref{WT}.

Taking derivatives and simplifying, we derive the formula of $J(x,y)$ in terms of the product of sine
functions. Furthermore, under the change of variables \eqref{xy}, it is not hard to verify that
\begin{gather*}
24x^3-y^2-12xy-6x-4y-1   =  \frac{16}{9}  \sin^2 \pi t_1\sin^2 \pi t_2 \sin^2 \pi (t_1+t_2),\\
  1+2y-3x^2   = \frac{16}{3} \sin^2 \frac{ \pi(t_1-t_2)}{3}
                         \sin^2\frac{ \pi(t_1+2t_2)}{3} \sin^2 \frac{ \pi(2t_1+t_2)}{3},
\end{gather*}
from which the second equality of \eqref{Jacobian} follows readily.
\end{proof}

\begin{defn} \label{Cheb*}
Under the change of variables \eqref{xy}, def\/ine for $k_1, k_2 \ge 0$,
\begin{gather*}
  P_{k_1,k_2}^{-\frac12, -\frac12}(x,y) :  = \TCC_{k_1+k_2, k_2, -k_1-2 k_2}(\tb),  \\
P_{k_1,k_2}^{\frac12, -\frac12}(x,y) :  = \frac{\TSC_{k_1+k_2+1, k_2, -k_1-2 k_2-1}(\tb)}{\TSC_{1,0,-1}(\tb)}, \\
 P_{k_1,k_2}^{-\frac12, \frac12}(x,y) :  = \frac{\TCS_{k_1+k_2+1, k_2+1, -k_1-2 k_2-2}(\tb)}{\TCS_{1,1,-2}(\tb)},\\
P_{k_1,k_2}^{\frac12, \frac12}(x,y) :  = \frac{\TSS_{k_1+k_2+2, k_2+1, -k_1-2 k_2-3}(\tb)}{\TSS_{2,1,-3}(\tb)}.
\end{gather*}
We call these functions generalized Chebyshev polynomials and, in particular, call
$P_{k}^{-\frac12,-\frac12}(x,y)$ and $P_{k}^{\frac12,\frac12}(x,y)$ the f\/irst kind and
the second kind, respectively.
\end{defn}

That these functions are indeed algebraic polynomials in $x$ and $y$ variables can be seen from the following
recursive relations, which can be derived from \eqref{TCCTCC}--\eqref{TCCTSS}.

\begin{prop}
\label{recurP}
For $\alpha,\beta=\pm \frac12$, $P^{\alpha,\beta}_{k_1,k_2}$
satisfy the recursion relation
\begin{gather}
      P^{\alpha,\beta}_{k_1+1,k_2}(x,y) =  6 x P^{\alpha,\beta}_{k_1,k_2}(x,y) -
  P^{\alpha,\beta}_{k_1+2,k_2-1}(x,y) - P^{\alpha,\beta}_{k_1-1,k_2+1}(x,y) \notag\\
 \hphantom{P^{\alpha,\beta}_{k_1+1,k_2}(x,y) =}{}
   - P^{\alpha,\beta}_{k_1+1,k_2-1}(x,y)-  P^{\alpha,\beta}_{k_1-2,k_2+1}(x,y) -  P^{\alpha,\beta}_{k_1-1,k_2}(x,y),
\label{eq:recurPx} \\
  P^{\alpha,\beta}_{k_1,k_2+1}(x,y) =    6 y P^{\alpha,\beta}_{k_1,k_2}(x,y) -
  P^{\alpha,\beta}_{k_1+3,k_2-2}(x,y) - P^{\alpha,\beta}_{k_1+3,k_2-1}(x,y) \notag\\
\hphantom{}{P^{\alpha,\beta}_{k_1,k_2+1}(x,y) =}
   - P^{\alpha,\beta}_{k_1-3,k_2+1}(x,y)-  P^{\alpha,\beta}_{k_1-3,k_2+2}(x,y) -  P^{\alpha,\beta}_{k_1,k_2-1}(x,y) \label{eq:recurPy}
\end{gather}
for $ k_1,k_2\ge 0$. Furthermore, the following symmetric relations hold,
\begin{gather}
\label{eq:symmPa}
 P^{\alpha,-\frac12}_{\mu,-\nu}(x,y) = P^{\alpha,-\frac12}_{\mu-3\nu,\nu}(x,y), \qquad P^{\alpha,\frac12}_{\mu,-\nu-1}(x,y) = -P^{\alpha,\frac12}_{\mu-3\nu,\nu-1}(x,y), \qquad \mu\ge 3\nu \ge 0, \\
\label{eq:symmPb}
 P^{-\frac12,\beta}_{-\mu,\nu}(x,y) = P^{-\frac12,\beta}_{\mu,\nu-\mu}(x,y), \qquad P^{\frac12,\beta}_{-\mu-1,\nu}(x,y) = -P^{\frac12,\beta}_{\mu-1,\nu-\mu}(x,y),\qquad \nu\ge \mu \ge 0.
\end{gather}
\end{prop}

\begin{proof}
The recursive relations \eqref{eq:recurPx} and \eqref{eq:recurPy} follow directly from
\eqref{TCCTCC} and~\eqref{TCCTSS}.
As for~\eqref{eq:symmPa} and~\eqref{eq:symmPb}, we resort to the following identities
of the trigonometric functions,
\begin{gather*}
\TCC_{\mu-\nu,-\nu,2\nu-\mu}(x,y) =   \TCC_{(\mu-3\nu)+\nu,\nu,\nu-\mu}(x,y),
\\
\TSC_{\mu-\nu+1,-\nu,2\nu-\mu-1}(x,y) =    \TSC_{(\mu-3\nu)+\nu+1,\nu,\nu-\mu-1}(x,y),
\\
\TCS_{\mu-(\nu+1)+1,-(\nu+1)+1,2\nu-\mu}(x,y) =   -\TCS_{(\mu-3\nu)+(\nu-1)+1, (\nu-1)+1,\nu-\mu}(x,y),
\\
  \TSS_{\mu-(\nu+1)+2,-(\nu+1)+1,2\nu-\mu-1}(x,y) =   -\TSS_{(\mu-3\nu)+ (\nu-1) +2 ,(\nu-1)+1,\nu-\mu-1}(x,y) ,
\\
\TCC_{-\mu+\nu,\nu,\mu-2\nu}(x,y) =    \TCC_{\mu+(\nu-\mu),\nu-\mu,\mu-2\nu}(x,y),
\\
\TCS_{-\mu+\nu+1,\nu+1,\mu-2\nu-2}(x,y) =    \TCS_{\mu+(\nu-\mu)+1,(\nu-\mu)+1,\mu-2\nu-2}(x,y),
\\
\TSC_{-(\mu+1)+\nu+1,\nu,\mu-2\nu}(x,y) =    -\TSC_{(\mu-1)+(\nu-\mu)+1,\nu-\mu,\mu-2\nu}(x,y),
\\
 \TSS_{-(\mu+1)+\nu+2,\nu+1,\mu-2\nu-2}(x,y) =    -\TSS_{(\mu-1)+ (\nu-\mu) +2 ,(\nu-\mu)+1,\mu-2\nu-2}(x,y) ,
\end{gather*}
which are derived from \eqref{SG2}--\eqref{SA2*-2}.
\end{proof}

The recursive relations \eqref{eq:recurPx} and \eqref{eq:recurPy} can be used to generate all polynomials
$P_{k_1,k_2}^{\a,\b}$ recursively. The task, however, is non-trivial. Below we describe an algorithm for the
recursion. Our starting point is
\begin{alignat*}{5}
& P^{-\frac12,-\frac12}_{0,0}(x,y)=1,\qquad  &&P^{-\frac12,-\frac12}_{1,0}(x,y)=x,\qquad   &&P^{-\frac12,-\frac12}_{0,1}(x,y)=y,& \\
& P^{\frac12,-\frac12}_{0,0}(x,y)=1,\qquad  &&P^{\frac12,-\frac12}_{1,0}(x,y)=6x+2,\qquad &&P^{\frac12,-\frac12}_{0,1}(x,y)=6x+3y+1,&\\
& P^{-\frac12,\frac12}_{0,0}(x,y)=1,\qquad  &&P^{-\frac12,\frac12}_{1,0}(x,y)=3x,\qquad   &&P^{-\frac12,\frac12}_{0,1}(x,y)=6y+2,&\\
& P^{\frac12,\frac12}_{0,0}(x,y)=1,\qquad   &&P^{\frac12,\frac12}_{1,0}(x,y)=6x+1,\qquad  &&P^{\frac12,\frac12}_{0,1}(x,y)=6x+6y+2.&
\end{alignat*}
The f\/irst few cases are complicated as the right side of the \eqref{eq:recurPx} and \eqref{eq:recurPy} involve
negative indexes, for which we need to use \eqref{eq:symmPa} and \eqref{eq:symmPb}. We give these cases
explicitly below
\begin{alignat*}{3}
& P^{-\frac12,-\frac12}_{2,0}(x,y)= 6x^2-2x-2y-1, \qquad && P^{-\frac12,-\frac12}_{1,1}(x,y)= 3xy-6x^2+x+2y+1, &
\\
&  P^{-\frac12,\frac12}_{2,0}(x,y)= 18x^2-3x-6y-3,\qquad &&P^{-\frac12,\frac12}_{1,1}(x,y)= 18xy+6x-18x^2+6y+3,&
\\
&  P^{\frac12,-\frac12}_{2,0}(x,y)= 36x^2-6y-3,\qquad && P^{\frac12,-\frac12}_{1,1}(x,y)= 18xy+6x+9y+2,&
\\
&  P^{\frac12,\frac12}_{2,0}(x,y)=36x^2-6y-3;\qquad &&P^{\frac12,\frac12}_{1,1}(x,y)=36xy+12x+12y+4;&
\end{alignat*}
\begin{gather*}
  P^{-\frac12,-\frac12}_{3,0}(x,y)= 36x^3-18xy-9x-6y-2,
      \\
  P^{-\frac12,\frac12}_{3,0}(x,y)= 108x^3-54xy-27x-12y-5,
      \\
  P^{\frac12,-\frac12}_{3,0}(x,y)= 216x^3-72xy-48x-24y-8,
       \\
  P^{\frac12,\frac12}_{3,0}(x,y)=216x^3-72xy-42x-18y-7;
 \\
  P^{-\frac12,-\frac12}_{0,2}(x,y)= 6y^2+10y-72x^3+36xy+18x+3,
      \\
  P^{-\frac12,\frac12}_{0,2}(x,y)= 36y^2+36y-216x^3+108xy+54x+9,
      \\
  P^{\frac12,-\frac12}_{0,2}(x,y)= 126xy+18y^2+36y+54x+10-216x^3,
       \\
  P^{\frac12,\frac12}_{0,2}(x,y)=144xy+36y^2+42y-216x^3+60x+11.
 \end{gather*}
The above formulas are derived from the recursive relations in the order of $(2,0)$, $(1,1)$, $(3,0)$, $(0,2)$, that is,
we need to deduce $(3,0)$ before proceeding to $(0,2)$. It should be pointed out that our polynomial
$P_{0,2}^{\a,\b}$ is of degree $3$, rather than degree $2$, which shows that our polynomials do not
satisfy the property of $\mathrm{span} \{P_{k_1,k_2}^{\a,\b}: k_1+k_2 \le  n\} = \Pi_n^2$. In particular,
they cannot be ordered naturally in the graded lexicographical order.

We shall show in the following section that our polynomials are best ordered in another
graded order for which the order is def\/ined by $2 k_1+ 3 k_2 = n$. We have displayed the
polynomials $P^{\alpha,\beta}_{k_1,k_2}(x,y)$ for all $2k_1+3k_2 \le 6$. In Algorithm~1 below
we give an algorithm for the evaluation of all $P^{\alpha,\beta}_{k_1,k_2}(x,y)$  with
$2k_1+3k_2=n$ and $n \ge 7$.

\begin{algorithm}[t]
\centerline{{\bf Algorithm~1.}~A recursive algorithm for the evaluation of $P^{\alpha,\beta}_{k_1,k_2}(x,y)$.}
\label{algo1}

\begin{itemize}\itemsep=0pt
\item[Step 1] if $n=2m$
\begin{gather*}
 P^{\alpha,\beta}_{m,0}(x,y) =  6x P^{\alpha,\beta}_{m-1,0}(x,y)
  -c_{\beta} P^{\alpha,\beta}_{m-2,1}(x,y)  - P^{\alpha,\beta}_{m-2,0}(x,y)
  - c_{\beta}P^{\alpha,\beta}_{m-3,1}(x,y),
   \end{gather*}
where $c_{\beta}=2$ if $\beta=-\frac12$, and $c_{\beta}=1$ if  $\beta=\frac12$;

\item[Step 2] for $k_2$ from $2-\bmod{(n,2)}$ with increment $2$ up to $\lfloor\frac{n}{3}\rfloor-2$ do
\begin{gather*}
  k_1 = \tfrac{n-3k_2}{2}  , \\
   P^{\alpha,\beta}_{k_1,k_2}(x,y)  =  6 x P^{\alpha,\beta}_{k_1-1,k_2}(x,y) -
  P^{\alpha,\beta}_{k_1+1,k_2-1}(x,y) - P^{\alpha,\beta}_{k_1-2,k_2+1}(x,y)
  \\
 \hphantom{P^{\alpha,\beta}_{k_1,k_2}(x,y)  =}{} - P^{\alpha,\beta}_{k_1,k_2-1}(x,y)-  P^{\alpha,\beta}_{k_1-3,k_2+1}(x,y) -  P^{\alpha,\beta}_{k_1-2,k_2}(x,y);
\end{gather*}
\item[Step 3] if $n=3m$
\begin{gather*}
  P^{\alpha,\beta}_{0,m}(x,y) =   6 y P^{\alpha,\beta}_{0,m-1}(x,y) -
  P^{\alpha,\beta}_{3,m-3}(x,y) - P^{\alpha,\beta}_{3,m-2}(x,y)-  P^{\alpha,\beta}_{0,m-2}(x,y)\\
\hphantom{P^{\alpha,\beta}_{0,m}(x,y) =}{}  + \begin{cases} -P^{\alpha,\beta}_{3,m-3}(x,y)-  P^{\alpha,\beta}_{3,m-2}(x,y), & \alpha= -\frac12,\\
    P^{\alpha,\beta}_{1,m-2}(x,y)+  P^{\alpha,\beta}_{1,m-1}(x,y), & \alpha= \frac12;
 \end{cases}
\end{gather*}
if $n=3m+1$
\begin{gather*}
  P^{\alpha,\beta}_{2,m-1}(x,y)  =  6 x P^{\alpha,\beta}_{1,m-1}(x,y) -
  P^{\alpha,\beta}_{3,m-2}(x,y) - P^{\alpha,\beta}_{0,m}(x,y) - P^{\alpha,\beta}_{2,m-2}(x,y)
  \\
\hphantom{P^{\alpha,\beta}_{2,m-1}(x,y)  =}{}
-  P^{\alpha,\beta}_{0,m-1}(x,y)
  - \begin{cases}  P^{\alpha,\beta}_{1,m-1}(x,y), & \alpha=-\frac12,\\
  0,  & \alpha=\frac12;
  \end{cases}
\end{gather*}
if $n=3m+2$
\begin{gather*}
  P^{\alpha,\beta}_{1,m}(x,y)  =
\begin{cases}
    3 x P^{\alpha,\beta}_{0,m}(x,y)-  P^{\alpha,\beta}_{2,m-1}(x,y) -  P^{\alpha,\beta}_{1,m-1}(x,y), &\alpha=-\frac12,
  \\   (6x+1)P^{\alpha,\beta}_{0,m}(x,y) -  P^{\alpha,\beta}_{2,m-1}(x,y) -  P^{\alpha,\beta}_{1,m-1}(x,y) , & \alpha=\frac12.
  \end{cases}
\end{gather*}
 \end{itemize}
\end{algorithm}

The polynomials $P_k^{\pm \frac12, \pm \frac12}$ def\/ined in the Def\/inition \ref{Cheb*}
satisfy an orthogonality relation. Let us def\/ine a weight function
$w_{\alpha,\beta}$ on the domain $\triangle^*$,
\begin{gather*}
  w_{\alpha,\beta}(x,y): =   \frac{(4 \pi^2)^ {\alpha+\beta}}{3^{2\alpha+\beta}} \left(1+2y-3x^2\right)^{\alpha}
      \left(24x^3-y^2-12xy-6x-4y-1\right)^{\beta}    \\
 \hphantom{w_{\alpha,\beta}(x,y)}{} \;     =   \left(\frac{4\pi^2}{3}\right)^{\alpha+\beta}  \left(\TSC_{1,0,-1}(\tb)\right)^{2\alpha} \left(\TCS_{1,1,-2}(\tb)\right)^{2\beta}
\end{gather*}
where the second equality follows from \eqref{2equality}. This weight function is closely related to
the Jacobian of the changing variables \eqref{xy}, as seen in Lemma~\ref{lem:Jacobian}. With
respect to this weight function, we def\/ine
\[
 \la f, g \ra_{w_{\alpha,\beta}}  : =  c_{\alpha,\beta}  \int_{\Delta^*} f(x,y)\overline{g(x,y)}
      w_{\alpha,\beta}(x,y) dxdy,
\]
where $c_{\alpha,\beta}:=  1/\int_{\triangle^*} w_{\alpha,\beta}(x,y) dxdy$ is a normalization constant; in
particular, $c_{-\frac{1}{2},-\frac{1}{2}}=4$,  $c_{\frac{1}{2},-\frac{1}{2}}=c_{-\frac{1}{2},\frac{1}{2}}=18/\pi^2$
and $c_{\frac{1}{2},\frac{1}{2}} =243/\pi^4$.  Since the change of variables
\eqref{xy} implies immediately that
\begin{gather} \label{int-int*}
c_{\alpha,\beta}
    \int_{\triangle^*} f(x,y) w_{\alpha,\beta}(x,y) dxdy = \frac{1}{|\triangle|} \int_{\triangle} f(\tb)
      \big(\TSC_{1,0,1}(\tb)\big)^{2\alpha+1}   \big(\TCS_{1,1,-2}(\tb)\big)^{2\beta+1} d\tb ,
\end{gather}
we can translate the orthogonality of $\TCC_\jb$, $\TSC_\jb$, $\TCS_\jb$ and $\TSS_\jb$ to that of
$P^{\alpha,\beta}_{k_1,k_2}$ for $\alpha,\beta=\pm \frac12$. Indeed, from Proposition~\ref{prop:trig-ortho*}
we can deduce the following theorem.

\begin{thm} \label{th:orthp}
For $\alpha = \pm \frac12, \beta = \pm \frac12$,
\begin{gather}\label{OPorthogonal}
\la P^{\alpha,\beta}_{k_1,k_2}, P^{\alpha,\beta}_{j_1,j_2} \ra_{ w_{\alpha,\beta}} =d_{k_1,k_2}^{\alpha,\beta}  \delta_{k_1,j_1}\delta_{k_2,j_2},
\end{gather}
where
\begin{alignat*}{3}
&    d_{k_1,k_2}^{-\frac12,-\frac12} :=     \begin{cases}  1, & k_1=k_2=0, \\
        \frac{1}{6}, &    k_1k_2=0 , \   k_1+k_2>0,\\
        \frac{1}{12}, &   k_1> 0, \  k_2>0,   \end{cases} \qquad
   &&   d_{k_1,k_2}^{\frac12,-\frac12} :=
       \begin{cases}
        \frac{1}{6}, &    k_1\ge 0, \ k_2=0,\\
        \frac{1}{12}, &   k_1 \ge 0,  \ k_2>0,   \end{cases} & \\
 &   d_{k_1,k_2}^{-\frac12,\frac12} :=
       \begin{cases}
        \frac{1}{6}, &    k_1= 0, \  k_2 \ge 0,\\
        \frac{1}{12}, &   k_1> 0, \ k_2 \ge 0 ,  \end{cases}  \qquad
    && d_{k_1,k_2}^{\frac12,\frac12} :=     \frac{1}{12},\ k_1\ge0, \  k_2\ge 0. &
\end{alignat*}
\end{thm}

\begin{proof}
All four cases follow from Proposition \ref{prop:trig-ortho*}. For $\alpha = \beta = -\frac12$, this is immediate.
For the other three cases, we observe that the weight function cancels the denominator in the def\/inition of $P_{k_1,k_2}^{\alpha,\beta}P_{j_1,j_2}^{\alpha,\beta}$ (see Def\/inition~\ref{Cheb*}),
which requires~\eqref{WT} in the case of $\alpha = \beta = \frac12$.
\end{proof}

Although the polynomials $P_{k_1,k_2}^{\pm \frac12,\pm \frac12}$ are mutually orthogonal,
they are not quite the usual orthogonal polynomials as we have seen from the recursive relations.
In fact, there are only two such polynomials with the total degree $2$, which is one less than
the number of monomials of degree $2$. As we have seen from the recursive relations, the structure
of these polynomials is much more complicated. To understand their structure, we study them
as solutions of the corresponding Sturm--Liouville problem in the following section.

\section[Sturm--Liouville eigenvalue problem and generalized Jacobi polynomials]{Sturm--Liouville eigenvalue problem\\ and generalized Jacobi polynomials}\label{section5}

Recall that our generalized trigonometric polynomials are solutions of the Sturm--Liouville eigenvalue
problems with corresponding boundary conditions. The Laplace operator becomes a~second-order
linear dif\/ferential operator in~$x$,~$y$ variables under the change of variables~\eqref{xy}. Using the fact
that $t_3 = - t_1-t_2$, we rewrite the change of variables~\eqref{xy} as
\begin{gather*}
 x =  \frac13\left(\cos\tfrac{2\pi (t_1-t_2)}{3}+\cos\tfrac{2\pi (t_2-t_3)}{3}+\cos\tfrac{2\pi (t_3-t_1)}{3}  \right),\\
 y =  \frac13(\cos 2\pi t_1+\cos 2\pi t_2+\cos 2\pi t_3  ).
\end{gather*}
A tedious but straightforward computation shows that
\begin{gather*}
   (\partial_{t_1}-\partial_{t_2})^2 + (\partial_{t_2}-\partial_{t_3})^2 + (\partial_{t_3}-\partial_{t_1})^2  \\
  \qquad  =  \frac{4\pi^2}{9}\big[ A_{1,1}(x,y) \partial_x^2  + 2 A_{1,2}(x,y)\partial_x \partial_y + A_{2,2}(x,y)\partial_y^2   + 6x\partial_x  + 18y\partial_y\big]
 =:  \frac{4\pi^2}{9} \CL_{-\frac12, -\frac12},
\end{gather*}
where we def\/ine
\begin{gather}
   A_{11}:=-6x^2+y+3x+2,  \qquad A_{12}=A_{21}:=-9xy+18x^2-6y-3, \nonumber\\
   A_{22}:=-18y^2+108x^3-54xy-27x-9y.\label{Aij}
\end{gather}
Consequently, we can translate the Laplace equation satisf\/ied by $\TCC_\kb$ into the equation
in $\CL_{-\frac12, -\frac12}$ for the polynomials $P_{k_1,k_2}^{- \frac12, - \frac12}(x,y)$. It is easy to
verify that the operator can be rewritten as
\begin{gather*}
 \CL_{-\frac12, -\frac12}    = - w_{\frac12, \frac12}\big[\partial_x w_{\alpha,\beta}\big(A_{11}\partial_x        +A_{12}\partial_y\big)  + \partial_y\omega^{-\frac12, -\frac12}\big(A_{21}\partial_x  + A_{22}\partial_y\big)\big] \\
 \hphantom{\CL_{-\frac12, -\frac12}    =}{} = -w_{\frac12, \frac12} \nabla^{\tr} w_{-\frac12, -\frac12} \Lambda \nabla,
\end{gather*}
where in the second line we have used
\[
   \nabla := (\partial_x, \partial_y)^\tr \qquad \hbox{and} \qquad  \Lambda := \begin{pmatrix} A_{11} & A_{12}\\ A_{21} & A_{22} \end{pmatrix}.
\]
It is not dif\/f\/icult to verify that the matrix $\Lambda$ is positive def\/inite in the interior of the do\-main~$\triangle^*$.
Indeed, $\det \Lambda = 3 F(x,y)$, where $F$ is def\/ined in Lemma~\ref{lem:Jacobian}, and $A_{1,1}(x,y) = 3(x-y) + 2(1+2y-3x^2)$ is positive if $x >y$ and it attains its minimal on the left most boundary, as seen by taking partial derivatives, in the rest of the domain, from which it is easy to verify that $A_{1,1} > 0$ in the interior of $\triangle^*$.
The expression of $\CL_{-\frac12,-\frac12}$ prompts the following def\/inition.

\begin{defn}
For $\alpha, \beta > -1$, def\/ine a second-order dif\/ferential operator
\begin{gather*}
\CL_{\alpha,\beta}:=  -w_{-\alpha,-\beta}\nabla^{\tr} w_{\alpha,\beta}\Lambda \nabla \\
\hphantom{\CL_{\alpha,\beta}}{}\;  =   -w_{-\alpha,-\beta} \left[\partial_x w_{\alpha,\beta}\big(A_{11}\partial_x +A_{12}\partial_y\big)  + \partial_y w_{\alpha,\beta}\big(A_{21}\partial_x  + A_{22}\partial_y\big)\right].
\end{gather*}
\end{defn}

The explicit formula of this dif\/ferential operator is given by
\begin{gather} \label{Lab}
\CL_{\alpha,\beta}=  -A_{11} \partial_x^2 -2A_{12} \partial_x \partial_y  - A_{22}\partial_y^2
   + B_1 \partial_x + B_2  \partial_y.
\end{gather}
where we def\/ine
\begin{gather*}
   B_1(x,y)   = 21x+12\alpha x+18\beta x+6\alpha+3, \\
   B_2(x,y)   =18x+36\alpha x+18\beta+45y+36\beta y+18\alpha y+9.
\end{gather*}
\begin{thm}\label{th:SA}
Let $\alpha,\beta>-1$. Then, the differential operator
\begin{gather*}
\CL_{\alpha,\beta}=-w_{-\alpha,-\beta}\nabla^{\tr} w^{\alpha,\beta}\Lambda \nabla
\end{gather*}
is self-adjoint and positive definite with respect to the inner product $\la\cdot,\cdot\ra_{w_{\alpha,\beta}}$.
\end{thm}

\begin{proof}
By Green's formula,
\begin{gather*}
  \iint_{\triangle^*}f\CL_{\alpha,\beta}g  w_{\alpha,\beta}dxdy
   = \iint_{\triangle^*} f\nabla^{\tr} w_{\alpha,\beta}\Lambda \nabla g dx dy
  =   -\iint_{\triangle^*} (\nabla f)^{\tr} \Lambda (\nabla g) w_{\alpha,\beta}  dx dy\\
\qquad\quad{}  + \oint_{\partial \triangle^*} w_{\alpha,\beta}  f \left[ (A_{11}\partial_x g + A_{12}\partial_y g )     dy
 - f (A_{22}\partial_y g + A_{21}\partial_x g )    dx\right] \\
 \qquad{} =   -\iint_{\triangle^*} (\nabla f)^{\tr} \Lambda (\nabla g) w_{\alpha,\beta}  dx dy\\
 \qquad\quad{}
 + \oint_{\partial \triangle^*}   w_{\alpha,\beta} f \left[(\partial_x g) (A_{11}dy - A_{21}dx)
 - (\partial_y g)  (A_{22}dx - A_{12}dy)\right],
 \end{gather*}
where $\partial \triangle^*$ denotes the boundary of the triangle. Recall that $\partial \triangle^*$ is
def\/ined by $F(x,y) =0$, where $F$ is def\/ined in Lemma \ref{lem:Jacobian}. It follows then
\begin{gather}  \label{dF}
          d F =  F_1 dx + F_2 d y =0, \qquad \hbox{where} \qquad F_1 = \frac{\partial F} {\partial x}, \qquad
            F_2 = \frac{\partial F} {\partial y}.
\end{gather}
On the other other hand, a quick computation shows that
 \begin{gather}
   F_1 A_{11}+ F_2 A_{21}   = - 6 (5x+1)F(x,y) =0, \label{FA1}\\
   F_1 A_{12}+F_2 A_{22}     = - 6 (3y+2x+1) F(x,y) =0\label{FA2}
 \end{gather}
on $\partial \triangle^*$. Solving \eqref{dF} and \eqref{FA1} shows that
$A_{11}dy - A_{21}dx =0$, whereas solving \eqref{dF} and \eqref{FA2} shows that
$A_{22}dx - A_{12}dy =0$ on $\partial \triangle^*$. Consequently, the integral over $\partial \triangle^*$
is zero and we conclude that
 \begin{gather*}
 -\iint_{\triangle^*} f \CL_{\alpha,\beta}g  w_{\alpha,\beta} dxdy
 =\iint_{\triangle^*} (\nabla f)^{\tr} \Lambda (\nabla g) w_{\alpha,\beta}  dx dy
 =  -\iint_{\triangle^*}g\CL_{\alpha,\beta}f w_{\alpha,\beta} dxdy ,
 \end{gather*}
which shows that $\CL_{\alpha,\beta}$ is self-adjoint and positive def\/inite.
\end{proof}

We consider polynomial solutions for the eigenvalue problem
$
    \CL_{\alpha,\beta}  u = \lambda u.
$
Dif\/ferential operators in the form of \eqref{Lab} have long been associated with orthogonal polynomials of two variables (see, for example, \cite{KS,Su}). However, in most of the studies, the coef\/f\/icients $A_{i,j}$ are chosen
to be polynomials of degree 2, which is necessary if, for each positive integer $n$, the solution of the
eigenvalue problem is required to consist of $n+1$ linearly independent polynomials of degree $n$, since
such choices ensure that the dif\/ferential operator preserves the degree of polynomials.
In our case, however, the coef\/f\/icient $A_{2,2}$ in \eqref{Aij} is of degree 3, which causes a number
of complications. In particular, our dif\/ferential operator does not preserve the polynomial degree; in other
words, it does not map $\Pi_n^2$ to $\Pi_n^2$, the space of polynomials of degree at most $n$ in two variables.

\begin{defn}
For $k_1, k_2 \ge 0$, the $m$-degree of the monomial $x^{k_1} y^{k_2}$ is def\/ined as
$|k|_*: = 2 k_1 + 3 k_2$. A polynomial $p$ in two variables is said to have $m$-degree $n$ if one monomial in $p$ has
$m$-degree of exactly $n$ and all other monomials in $p$ have $m$-degree at most $n$.
For $n \in \NN_0$, let~$\Pi_n^*$ denote the space of polynomials of $m$-degree at most $n$; that is,
\[
    \Pi_n^*: = \mathrm{span} \big\{x^{k_1}y^{k_2}: 0\le k_1,k_2;\  2k_1+3k_2\le n \big\}.
\]
\end{defn}

The dimension of the space $\Pi_n^*$ is the same as that of $\CH_n^{\cc}$, by \eqref{dimCT},
\begin{equation} \label{dimPin*}
   \dim \Pi_n^*  = \tfrac12 \big(3\lfloor\tfrac{n}{3}\rfloor-2n\big)  \big(\lfloor\tfrac{n}{3}\rfloor+1\big)
  -\big(\lfloor\tfrac{n}{2} \rfloor-n-1\big) \big(\lfloor\tfrac{n}{2} \rfloor+1\big).
\end{equation}
Here is a list of the dimension for small $n$:
\begin{table}[ht]
\centering
\begin{tabular}{|c|c|c|c|c|c|c|c|c|c|c|c|c|}
\hline
$n$  &  1  &  2 & 3 & 4 & 5 & 6 & 7 & 8 & 9 & 10 & 11 & 12\\
\hline
$\dim \Pi_n^*$ & 1 &  2  &  3 & 4 & 5 & 7 & 8 & 10 & 12 & 14 & 16 & 19\\
\hline
\end{tabular}
\end{table}

The name $m$-degree is coined in~\cite{Patera} after the marks, or co-marks, in the root system for the simple compact Lie group, where the case of the group $G_2$ is used as an example. For polynomials graded by the
$m$-degree, we introduce an ordering among monomials.

\begin{defn}
For any $k, j\in \NN_0^2$, we def\/ine an order $\prec$ by $k\prec j $ if $2j_1+3j_2>2k_1+3k_2$ or $2(k_1-j_1)
=3(j_2-k_2)>0$, and  $k \preceq j$ if $k\prec j$ or $k=j$. We call $\prec$ the $*$-order. If $p (x,y) =
\sum\limits_{(k_1,k_2)\preceq (m,n)} c_{k_1,k_2} x^{k_1} y^{k_2}$ with $c_{m,n} \ne 0$, we
call $c_{m,n}x^my^n$ the leading term of $p$ in the $*$-order.
\end{defn}

For $m,n \ge 0$,  def\/ine
\begin{gather*}
\Pi_{m,n}^* = \mathrm{span}\big\{ x^{j}y^k: (j,k) \preceq (m,n)   \big\}.
\end{gather*}
It is easy to see that $\Pi_n^*=\Pi_{2n-3 \lfloor \frac{2n}{3} \rfloor,  2 \lfloor \frac{2n}{3} \rfloor -n }^*$.
The $*$-order is well-def\/ined. The following lemma justif\/ies our def\/initions.

\begin{lem} \label{PitoPi}
For $m,n \ge 0$, the operator $\CL_{\alpha,\beta}$ maps $\Pi_{m,n}^*$ onto $\Pi_{m,n}^*$.
\end{lem}

\begin{proof}
We apply the operator $\CL_{\a,\b}$ on the monomial $x^j y^k$. The result is
\begin{gather*}
\CL_{\alpha,\beta} x^j y^k  =  -A_{11} \partial_x^2x^j y^k - A_{22}\partial_y^2x^jy^k
 - 2A_{12} \partial_x \partial_yx^jy^k  + B_1\partial_xx^jy^k + B_2\partial_yx^jy^k \\
\hphantom{\CL_{\alpha,\beta} x^j y^k}{}
 =   \left[ 6(j^2+3k^2+3jk) + 3(5+4\alpha +6\beta )j + 3(9+6\alpha+12\beta)k\right] x^j y^k \\
\hphantom{\CL_{\alpha,\beta} x^j y^k=}{}
     - 108 k(k-1)x^{j+3}y^{k-2} - j(j-1) x^{j-2} y^{k+1} + 18k(3k-2-2j+2\alpha) x^{j+1} y^{k-1}   \\
\hphantom{\CL_{\alpha,\beta} x^j y^k=}{}
  +3j(-j+2+4k+2\alpha) x^{j-1} y^{k} + 9k(k+2\beta) x^jy^{k-1} \\
\hphantom{\CL_{\alpha,\beta} x^j y^k=}{}
    - 2j(j-1) x^{j-2} y^k + 27k(k-1) x^{j+1} y^{k-2} + 6jk x^{j-1} y^{k-1}.
\end{gather*}
Introducing the notation
\[
\Upsilon= \left\{(0,0), (0,1), (1,0), (1,1), (2,1), (3,2), (4,2), (4,3), (5,3)    \right\},
\]
we write the expression as
\begin{equation} \label{L-monomial}
\CL_{\alpha,\beta} x^j y^k  =    \sum_{(\mu,\nu)\in \Upsilon}a_{\mu,\nu}^{j,k}
    x^{j-2\mu+3\nu} y^{k+\mu-2\nu},
\end{equation}
where
\begin{gather*}
    a^{j,k}_{0,0} = 6\big(j^2+3k^2+3jk\big) + 3(5+4\alpha +6\beta )j + 3(9+6\alpha+12\beta)k, \\
     a^{j,k}_{0,1} = - 108 k(k-1), \quad a^{j,k}_{1,0}=- j(j-1), \qquad a^{j,k}_{1,1}= 18k(3k-2-2j+2\alpha)  ,
 \\
  a^{j,k}_{2,1}=3j(-j+2+4k+2\alpha) , \qquad  a^{j,k}_{3,2} =9k(k+2\beta) ,
 \\
   a^{j,k}_{4,2}= - 2j(j-1),\qquad a^{j,k}_{4,3} =  27k(k-1),  \qquad a^{j,k}_{5,3}= 6jk.
\end{gather*}
From this computation, it follows readily that $\CL_{\a,\b}$ maps $\Pi_{m,n}^*$ into $\Pi_{m,n}^*$.
Furthermore, with respect to  the $*$-order, it is easy to see that $a^{m,n}_{0,0} x^m y^n$
is the leading term of $\CL_{\a,\b}$ by~\eqref{L-monomial}, which shows that $\CL_{\a,\b}$ maps
$\Pi_{m,n}^*$ onto $\Pi_{m,n}^*$.
\end{proof}

The identity \eqref{L-monomial} also shows that $\CL_{\a,\b}$ has a complete set of eigenfunctions
in $\Pi_{m,n}^*$.

\begin{thm}
For $\a,\b \ge  -1/2$ and $k_1, k_2  \ge 0$, there exists a polynomial $P_{k_1,k_2}^{\a,\b} \in \Pi_{k_1,k_2}^*$ with the leading  term
$x^{k_1}y^{k_2}$
such that
\begin{gather}\label{JacobiP}
    \CL_{\a,\b} P_{k_1,k_2}^{\a,\b} = \lambda_{k_1,k_2}^{\a,\b} P_{k_1,k_2}^{\a,\b},
\end{gather}
where
\begin{gather}\label{eigenvalue}
  \lambda_{k_1,k_2}^{\alpha,\beta} :=  \frac{3}{2}|k|_* (|k|_*+5+4\alpha+6\beta)+\frac{9}{2}k_2 (k_2+1+2\beta).
\end{gather}
Furthermore, if we require all the polynomials are orthogonal to each other with respect to the inner product $\langle\cdot,\cdot\rangle_{w_{\a,\b}}$,  then $P_{k_1,k_2}^{\a,\b}$ is uniquely determined by its leading term
in the $*$-order.
\end{thm}

\begin{proof}
We f\/irst apply the Gram--Schmidt orthogonality process to monomials $\left\{ x^{k_1}y^{k_2} \right\}$
in the $*$-order, which uniquely determines a complete system of orthogonal polynomials with
leading term  $x^{k_1}y^{k_2} $ with respect to $\la \cdot, \cdot\ra_{w_{\a,\b}}$; that is,
$P_{0,0}^{\a,\b}(x,y)=1$ and
\begin{gather*}
 P_{k_1,k_2}^{\a,\b}(x,y)= x^{k_1}y^{k_2} - \sum_{(j_1,j_2)\prec (k_1,k_2)}
  \frac{ \la x^{k_1}y^{k_2} , P_{j_1,j_2}^{\a,\b}\ra_{w_{\a,\b}}}{ \la P_{j_1,j_2}^{\a,\b}, P_{j_1,j_2}^{\a,\b}\ra_{w_{\a,\b}} }  P_{j_1,j_2}^{\a,\b}(x,y),  \qquad (0,0) \prec (k_1,k_2).
\end{gather*}
The Gram--Schmidt orthogonality and  Lemma \ref{PitoPi} show that
\begin{gather}\label{CLSpace}
\CL_{\a,\b} P_{k_1,k_2}^{\a,\b}(x,y)  \in \mathrm{span}\big\{P_{j_1,j_2}^{\a,\b}(x,y):   (j_1,j_2) \preceq (k_1,k_2) \big\} = \Pi_{k_1,k_2}^*.
\end{gather}
Evidently, $\CL_{\a,\b}P_{0,0}^{\a,\b} = 0=  a^{0,0}_{0,0} P_{0,0}^{\a,\b}$. We apply induction. Assume
that
\begin{gather*}
 \CL_{\a,\b}P_{j_1,j_2}^{\a,\b}  = a^{0,0}_{j_1,j_2} P_{j_1,j_2}^{\a,\b}, \qquad  (j_1,j_2) \prec (k_1,k_2).
\end{gather*}
It then follows from Theorem \ref{th:SA} and the orthogonality of $P_{k_1,k_2}^{\a,\b}$ that
\begin{gather*}
 \la  \CL_{\a,\b}  P_{k_1,k_2}^{\a,\b} ,   P_{j_1,j_2}^{\a,\b} \ra_{w_{\a,\b}}=
 \la  P_{k_1,k_2}^{\a,\b} ,   \CL_{\a,\b} P_{j_1,j_2}^{\a,\b} \ra_{w_{\a,\b}} =  a_{0,0}^{j_1,j_2}\la  P_{k_1,k_2}^{\a,\b} ,   P_{j_1,j_2}^{\a,\b} \ra_{w_{\a,\b}}=0,
 \end{gather*}
so that, as a consequence of \eqref{CLSpace},
 \begin{gather}
 \label{Constc}
\CL_{\a,\b}  P_{k_1,k_2}^{\a,\b}  =  c   P_{k_1,k_2}^{\a,\b}.
\end{gather}
Comparing the leading term of the above identity, we obtain from \eqref{L-monomial} that
\begin{gather*}
     a_{0,0}^{k_1,k_2} x^{k_1,k_2} =  c x^{k_1,k_2},
\end{gather*}
which gives $c= a_{0,0}^{k_1,k_2}$. Ultimately, this inductive process  shows that
\begin{gather*}
\CL_{\a,\b}  P_{k_1,k_2}^{\a,\b}   =  \lambda_{k_1,k_2}^{\a,\b} P_{k_1,k_2}^{\a,\b},
\qquad \hbox{with} \quad
 \lambda_{k_1,k_2}^{\a,\b} =a_{0,0}^{k_1,k_2}.
\end{gather*}
As shown in the proof of Lemma \ref{PitoPi},
\begin{gather*}
\lambda_{k_1,k_2}^{\a,\b}=
6\big(k_1^2+3k_2^2+3k_1k_2\big) + 3(5+4\alpha +6\beta )k_1 + 3(9+6\alpha+12\beta)k_2\\
\hphantom{\lambda_{k_1,k_2}^{\a,\b}}{} = \frac{3}{2}(2k_1+3k_2)((2k_1+3k_2)+4+4\alpha+6\beta)+\frac{9}{2}k_2 (k_2+2+2\beta) + 3k_1,
\end{gather*}
which is \eqref{eigenvalue} since $|k|_* = 2k_1+3k_2$ by def\/inition.

Moreover,  suppose  $\wt P_{k_1,k_2}^{\a,\b} (x,y) \in \Pi_{k_1,k_2}^*$ is another polynomial with the leading term $x^{k_1}y^{k_2}$ such that
  \begin{gather*}
 \CL_{\a,\b} \wt P_{k_1,k_2}^{\a,\b} (x,y) =  \lambda  \wt  P_{k_1,k_2}^{\a,\b} (x,y),\\
  \la \wt P_{k_1,k_2}^{\a,\b}, p \ra_{w_{\a,\b}} =0,\qquad \forall \,  p \in \mathrm{span} \left\{ x^{j_1}y^{j_2}: (j_1,j_2) \prec (k_1,k_2)  \right\}.
 \end{gather*}
Using the same argument that determines $c$ in \eqref{Constc}, we see that $\lambda=
\lambda_{k_1,k_2}^{\a,\b} =a_{0,0}^{k_1,k_2} $. Moreover, it is easy to see that
\begin{gather*}
  P_{k_1,k_2}^{\a,\b} -\wt P_{k_1,k_2}^{\a,\b}   \in \mathrm{span} \left\{ x^{j_1}y^{j_2}: (j_1,j_2) \prec (k_1,k_2)     \right\},
\\
  \big\langle P_{k_1,k_2}^{\a,\b}-\wt P_{k_1,k_2}^{\a,\b}, P_{j_1,j_2}^{\a,\b} \big\rangle_{w_{\a,\b}} =0,
   \qquad \forall \,  (j_1,j_2) \prec (k_1,k_2) .
\end{gather*}
This f\/inally leads to $P_{k_1,k_2}^{\a,\b}-\wt P_{k_1,k_2}^{\a,\b} = 0$, which  shows that
$P_{k_1,k_2}^{\a,\b}$ is  uniquely determined by its leading term in the $*$-order and the
orthogonality $\langle P_{k_1,k_2}^{\a,\b}, x^{j_1}y^{j_2}\rangle_{w_{\a,\b}}=0$ for all
$(j_1,j_2) \prec (k_1,k_2) $. This completes the proof.
\end{proof}

Let $P_{k_1,k_2}^{\alpha,\beta}$ be orthogonal  to each other with respect
to the inner product $\la\cdot,\cdot\ra_{w_{\alpha,\beta}}$. The f\/irst few polynomials
and the eigenvalues can be readily checked to be
\begin{gather*}
  P_{0,0}^{\a,\b}(x,y)  = 1,\qquad     \lambda_{0,0}^{\a,\b}= 0; \\
  P_{1,0}^{\a,\b}(x,y)  = x+\frac{1+2\a}{7+4\a+6\b},\qquad     \lambda_{1,0}^{\a,\b} = 3(7+4\a+6\b);
  \\
    P_{0,1}^{\a,\b}(x,y)  =      y+ \frac{3(1+2\a)}{4+\a+3\b}x+\frac{5+5\a+11\b+2\a\b+6\b^2+4\a^2}{(4+\a+3\b)(5+2\a+4\b)},\\
\hphantom{P_{0,1}^{\a,\b}(x,y)  =}{}~
         \lambda_{0,1}^{\a,\b} =9(5+2\a+4\b);
  \\
    P_{2,0}^{\a,\b}(x,y)  =x^2-  \frac{2y}{3(3+2\a )}y+\frac{  4(\a +1)(2\a -1)}{(3+2\a )(4\a +11+6\b )} x
  \\
\hphantom{P_{2,0}^{\a,\b}(x,y)  =}{}
 + \frac{-105-86\a -120\b -36\b ^2-48\b \a +8\a ^2+24\a ^3}{3(3+2\a )(4\a +11+6\b )(4\a +6\b +9)},\\
\hphantom{P_{2,0}^{\a,\b}(x,y)  =}{}~
  \lambda_{2,0}^{\a,\b}
  =6(9+4\b+6\a),
  \\
    P_{1,1}^{\a,\b}(x,y)  = xy+\frac{3(2\a -1)}{5+\a +3\b }x^2+ \frac{6\b \a +11\a +15\b +2\a ^2+27}{(5+\a +3\b )(4\a +6\b +13)}y
  \\
  \hphantom{P_{1,1}^{\a,\b}(x,y)  =}{}
  +\frac{119\!+\!229\b\! +\!40\a ^3\!+\!36\b ^3\!+\!111\a\! +\!80\b \a ^2\!+\!156\b ^2\!+\!132\b \a\! +\!140\a ^2\!+\!36\b ^2\a }{(2\a +7+4\b )(4\a +6\b +13)(5+\a +3\b )}x
    \\
\hphantom{P_{1,1}^{\a,\b}(x,y)  =}{}
  +\frac{8\a ^3+6\a ^2+4\b \a ^2+12\b ^2\a +9\a +28\b \a +70+93\b +30\b ^2}{(2\a +7+4\b )(4\a +6\b +13)(5+\a +3\b )},\\
\hphantom{P_{1,1}^{\a,\b}(x,y)  =}{}~
 \lambda_{1,1}^{\a,\b}
   = 6(14+9\b+5\a).
  \end{gather*}

For each $P^{\a,\b}_{m,n}$, \eqref{CLSpace} shows that the $\CL_{\a,\b}P^{\a,\b}_{m,n}$
involves only $P_{j_1,j_2}^{\a,\b}$ with $(j_1,j_2)$ in
\[
\Gamma_{m,n} := \left\{(j_1,j_2) \in \NN^2:  (j_1,j_2) \preceq (m,n)\right\}.
\]
This set of dependence of the polynomial solution is determined by the $*$-ordering. Indeed,
it is easy to see that
\begin{gather*}
\Gamma_{m,n}= \Gamma_{m,n}^{+} \cup\Gamma_{m,n}^{-},\\
\Gamma_{m,n}^{+} :=\left\{  (m-2p+3q,n+p-2q))\in \ZZ^2:  0 \le q \le \lfloor \tfrac{p+n}{2}\rfloor ,
0\le p \le 2m+3n \right\} ,
\\
 \Gamma_{m,n}^{-}  := \left\{ (m-2p+3q,n+p-2q) \in \ZZ^2:    \lceil \tfrac{2p-m}{3} \rceil \le q \le -1 ,    1\le p   \le \lfloor \tfrac{m-3}{2}\rfloor \right\}.
\end{gather*}
For $p$, $q$ as in $\Gamma_{k_1,k_2}$ but not both $0$, we have that for $\alpha,\beta\geq -\frac12$,
\begin{gather*}
 \lambda_{k_1,k_2}^{\a,\b} -  \lambda_{{k_1-2p+3q,k_2+p-2q}}^{\a,\b}
  = 3(2k_1-2p+3q+2\a+1)p + 9(2k_2+p-2q+2\b+1)q >0,
\end{gather*}
which shows that $\lambda_{k_1,k_2}^{\a,\b}\neq  \lambda_{j_1,j_2}^{\a,\b} $ for any $(j_1,j_2)\in  \Gamma_{k_1,k_2}^{+}$. This implies that polynomial solutions of the same $m$-degree below
to dif\/ferent eigenvalues. Moreover, if $ \lambda_{k_1,k_2}^{\a,\b}=  \lambda_{j_1,j_2}^{\a,\b} $ for
$(j_1,j_2)\prec (k_1,k_2)$, then $(j_1,j_2)\in  \Gamma_{k_1,k_2}^{-}$.

In the case of $(\a,\b) = (-\frac12, -\frac12)$, our polynomials $P_{k_1,k_2}^{\a,\b}$ agree
with the generalized Chebyshev polynomial that we def\/ined in the last section. For the other
three cases of $(\a,\b) = (\pm \frac12, \pm \frac12)$, this requires proof. Let us denote the
Chebyshev polynomials temporarily by $Q_{k_1,k_2}^{\alpha,\beta}$, $(\a,\b) = (-\frac12, -\frac12)$.
It is not hard to see,  from Algorithm~1, that the leading term of
$Q^{\alpha,\beta}_{k_1,k_2}$ is $cx^{k_1}y^{k_2}$ with certain $c>0$, which implies that
$\mathrm{span}\big\{Q_{j_1,j_2}^{\a,\b}(x,y):   (j_1,j_2) \preceq (k_1,k_2) \big\} = \Pi_{k_1,k_2}^*.$
Thus, we can write
\begin{gather} \label{LQ}
   \CL_{\alpha,\beta} Q^{\alpha,\beta}_{k_1,k_2}(x,y) = \lambda_{k_1,k_2}^{\alpha,\beta}
    Q^{\alpha,\beta}_{k_1,k_2}(x,y) + \sum_{(j_1,j2)\prec (k_1,k_2)}  c_{j_1,j2}^{k_1,k_2}
    Q^{\alpha,\beta}_{j_1,j2}(x,y).
\end{gather}
On the other hand, by the orthogonality and the self-adjointness of $\CL_{\a,\b}$, for
any $(l_1,l_2)\prec (k_1,k_2)$,
\begin{gather*}
      \big(\CL_{\alpha,\beta} Q^{\alpha,\beta}_{k_1,k_2}, Q^{\alpha,\beta}_{l_1,l_2}\big)_{w_{\alpha,\beta}}
   =    \big( Q^{\alpha,\beta}_{k_1,k_2}, \CL_{\alpha,\beta} Q^{\alpha,\beta}_{l_1,l_2}\big)_{w_{\alpha,\beta}}   \\
   \hphantom{\big(\CL_{\alpha,\beta} Q^{\alpha,\beta}_{k_1,k_2}, Q^{\alpha,\beta}_{l_1,l_2}\big)_{w_{\alpha,\beta}}}{}
   = \left( Q^{\alpha,\beta}_{k_1,k_2}, \lambda_{l_1,l_2}^{\alpha,\beta} Q^{\alpha,\beta}_{l_1,l_2}
       + \sum_{(j_1,j_2)\prec (l_1,l_2)}  c_{j_1,j_2}^{l_1,l_2} Q^{\alpha,\beta}_{j_1,j_2}\right)_{w_{\alpha,\beta}} =  0.
\end{gather*}
As a result, we deduce from \eqref{LQ} that
 \begin{gather*}
 \CL_{\alpha,\beta} Q^{\alpha,\beta}_{k_1,k_2}(x,y) =
    \lambda_{k_1,k_2}^{\alpha,\beta} Q^{\alpha,\beta}_{k_1,k_2}(x,y).
 \end{gather*}
Consequently, up to a constant multiple, we see that $Q_{k_1,k_2}^{\a,\b}$ coincides with the
Jacobi polynomials.

\begin{cor}
The Chebyshev polynomials defined in Definition~{\rm  \ref{Cheb*}} satisfy the equation \eqref{JacobiP}.
\end{cor}

In particular, this shows that the Chebyshev polynomials are elements in $\Pi_{|k|_*}^*$ and they are
determined, as eigenfunctions of $\CL_{\a,\b}$, uniquely by the leading term in the $*$-order.

\section{Cubature rules for polynomials}\label{section6} 

In the case of the equilateral triangle, the cubature rules for the trigonometric functions
are transformed into cubature rules of high quality for polynomials on the region bounded
by the Steiner's hypocycloid. In this section we discuss analogous results for the cubature
rules in the Section~\ref{section3}. To put the results in perspective, let us f\/irst recall the relevant
background.

Let $w$ be a nonnegative weight function def\/ined on a compact set $\Omega$ in $\RR^2$.
A cubature rule of degree $2n-1$ for the integral with respect to $w$ is a sum of point
evaluations that satisf\/ies
\begin{gather*}
  \int_\Omega f(x) w(x) dx = \sum_{j=1}^N \lambda_j f(x_j),
   \qquad \lambda_j \in \RR
\end{gather*}
for every $f \in \Pi_{2n-1}^2$. It is well-known that a cubature rule of degree $2n-1$
exists only if $N \ge \dim \Pi_{n-1}^2 = n(n+1)/2$. A cubature that attains such a lower
bound is called Gaussian. Unlike one variable, the Gaussian cubature rule exists
rarely and it exists if and only if the corresponding orthogonal polynomials of degree
$n$, all $n+1$ linearly independent ones,  have $n(n+1)/2$ real distinct common
zeros. We refer to~\cite{DX,St} for these results and further discussions. At the
moment there are only two regions with weight functions that admit the Gaussian
cubature rule. One is the region bounded by the Steiner's hypocycloid and the Gaussian
cubature rule is obtained by transformation from one cubature rule for trigonometric
functions on the equilateral triangle.

\subsection[Gaussian cubature rule of $m$-degree]{Gaussian cubature rule of $\boldsymbol{m}$-degree}
We f\/irst consider the case of $w_{\frac12,\frac12}$, which turns out to admit the
Gaussian cubature rule in the sense of $m$-degree.

\begin{thm}
For $w_{\frac{1}{2},\frac{1}{2}}$ on $\triangle^*$, the cubature rule
\begin{gather} \label{GaussCuba}
  c_{\frac12,\frac12}\iint_{\Delta^*} f(x,y) w^{\frac{1}{2},\frac12}(x,y) dxdy
       = \frac{12}{(n+5)^2} \sum_{\jb \in \Upsilon_{n+5}^{{\degree}}}
     \big|\TSS_{2,1,-3}\big(\tfrac{\jb}{n+5}\big)\big|^2
          f\big(x\big(\tfrac{\jb}{n+5}\big),y\big(\tfrac{\jb}{n+5}\big)\big),
\end{gather}
is exact for all polynomials $f \in P_{2n-1}^*$.
\end{thm}

\begin{proof}
Using \eqref{int-int*} with $\alpha = \beta = \frac12$ and \eqref{WT}, we see that
\begin{gather}\label{int-int-1/2}
c_{\frac12,\frac12} \int_{\triangle^*} f(x,y) w_{\frac12,\frac12}(x,y) dxdy  =
   \frac{1}{|\triangle|}  \int_{\triangle} f (x(\tb),y(\tb))  \left[ \TSS_{2,1,-3}(\tb) \right]^2 d\tb.
\end{gather}
By \eqref{hypocycloid} and \eqref{2equality}, $[\TSS_{2,1,-3}(\tb)]^2$ has $m$-degree 10, so that
$f (x(\tb),y(\tb))  \left[\TSS_{2,1,-3}(\tb) \right]^2 \in \CH_{2n+9}^{cc}$ if $f \in \Pi_{2n-1}^*$.
Since $\TSS_{2,1,-3}(\tb)$ vanishes on the boundary of $\triangle$, applying the cubature
rule \eqref{cuba-HHT2} of degree $2n +9$ to the right hand side of \eqref{int-int-1/2} gives
the stated result.
\end{proof}

What makes this result interesting is the fact that, by \eqref{dimCT},
\[
    |\Upsilon_{n+5}^{{\degree}} | = |\Gamma_{n+5}^{\ss} |= |\Gamma_{n-1}^{\cc} | = \dim \Pi_{n-1}^*,
\]
which shows that the cubature rule \eqref{GaussCuba} resembles the Gaussian cubature rule
under the $m$-degree. Furthermore, it turns out that it is again characterized by the common
zeros of ortho\-go\-nal polynomials.  Let~$Y_n^{{\degree}}$ be the image of
$\big\{\tfrac{\jb}{n}:  \jb \in\Upsilon_n^{\ss} \big\}$ under the mapping $\tb \mapsto x$,
\begin{gather*}
Y_n^{{\degree}} :  =  \big\{\big(x\big(\tfrac{\jb}{n} \big), y\big(\tfrac{\jb}{n} \big) \big)  :\; \jb\in \Upsilon_n^{{\degree}} \big\},
\end{gather*}
which is the set of nodes for \eqref{GaussCuba}. Then all polynomials $P_{k_1,k_2}^{\frac12,\frac12}$ with
$m$-degree~$n$ vanish on~$Y_n^{\degree}$.

\begin{thm} \label{thm:zeroU}
The set $Y_{n+5}^{\degree}$ is the variety of the polynomial ideal
$\big\langle P_{k_1,k_2}^{\frac12,\frac12}(x): 2 k_1 + 3 k_2  = n \big\rangle$.
\end{thm}

\begin{proof}
By the def\/inition of $P_{k_1,k_2}^{\frac12,\frac12}$, it suf\/f\/ices to show that
\begin{gather}
\label{Yon+5}
  \TSS_{\kb} \left(\tfrac{\jb}{n+5}\right) =0 \qquad
   \hbox{for}  \quad \jb\in \Upsilon, \  \kb\in \Gamma \quad  \hbox{and}  \quad k_1-k_{3}=n+5.
\end{gather}
Directly form its def\/inition,
\begin{gather*}
\TSS_{\kb} \big(\tfrac{\jb}{n+5}\big) =
  \frac{1}{3} \Big[ \sin \tfrac{\pi(k_1-k_3)(j_1-j_3)}{3(n+5)} \sin \tfrac{\pi k_2 j_2}{n+5}
  + \sin \tfrac{\pi(k_1-k_3)(j_2-j_1)}{3(n+5)} \sin \tfrac{\pi k_2 j_3}{n+5}
\\
\hphantom{\TSS_{\kb} \big(\tfrac{\jb}{n+5}\big) = }{}
 + \sin \tfrac{\pi(k_1-k_3)(j_3-j_2)}{3(n+5)} \sin \tfrac{\pi k_2 j_1}{n+5}  \Big]
\\
\hphantom{\TSS_{\kb} \big(\tfrac{\jb}{n+5}\big)}{}
=  \frac{1}{3} \Big[ \sin \tfrac{\pi(j_1-j_3)}{3} \sin \tfrac{\pi k_2 j_2}{n+5}
+ \sin \tfrac{\pi(j_2-j_1)}{3} \sin \tfrac{\pi k_2 j_3}{n+5}
+ \sin \tfrac{\pi(j_3-j_2)}{3} \sin \tfrac{\pi k_2 j_1}{n+5}  \Big],
\end{gather*}
Since $j_1\equiv j_2 \equiv j_3$ $(\bmod{3})$, we conclude then $\TSS_{\kb} \big(\tfrac{\jb}{n+5}\big)=0$.
The proof is completed.
\end{proof}

In \cite{Patera}, the existence of the Gaussian cubature rule in the sense of $m$-degree and the connection to
orthogonal polynomials were established in the context of compact simple Lie groups. The case of the group
$G_2$ was used as an example, where a numerical example was given. The domain~$\triangle^*$ and
the one in~\cite{Patera} dif\/fer by an af\/f\/ine change of variables.

 Our results give explicit nodes and weights of the cubature rule and provide
further explanation for the result.

\begin{figure}[htb]
\hfill%
\begin{minipage}{0.4\textwidth}\includegraphics[width=1\textwidth]{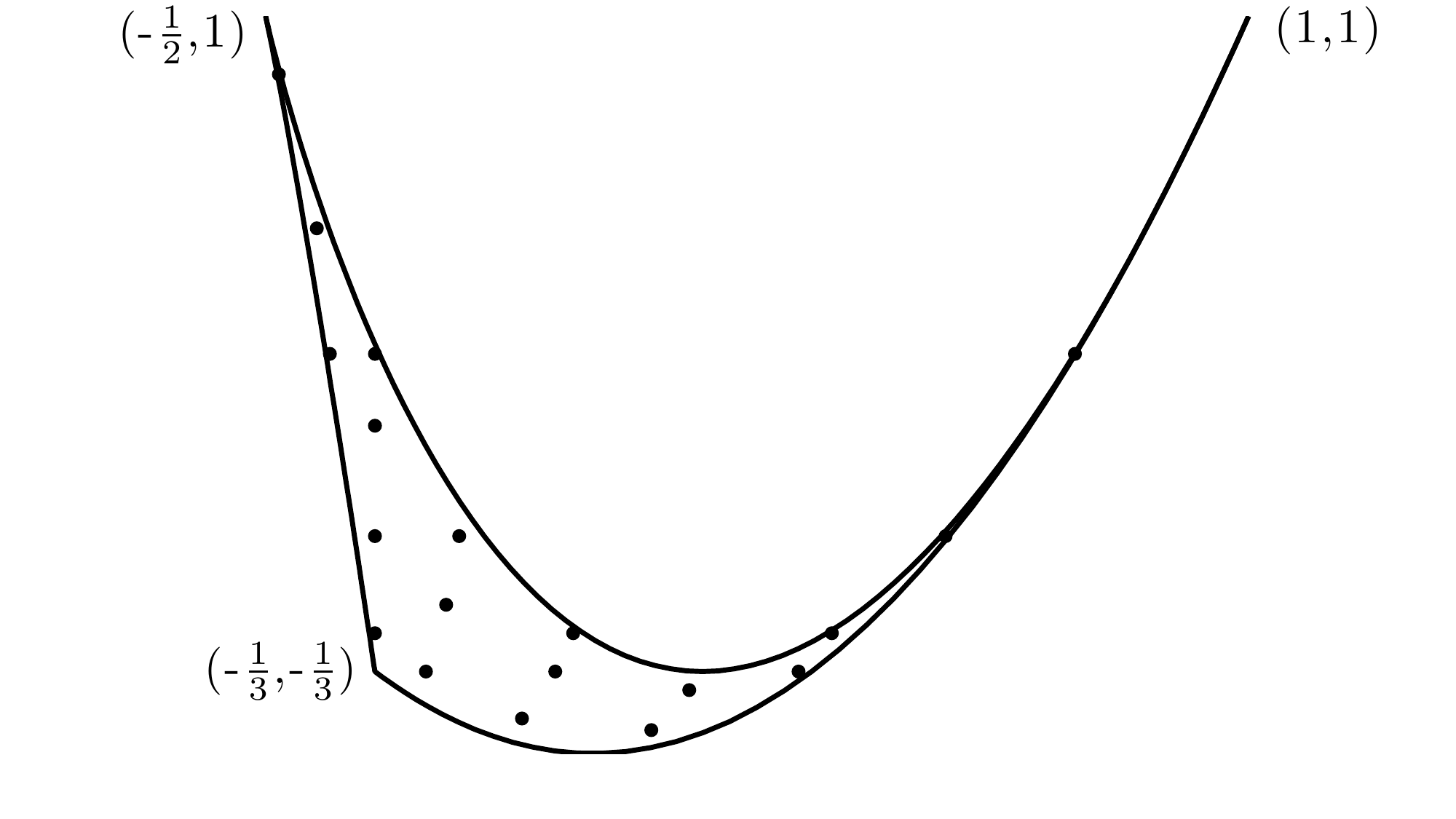}\\
\centering Chebyshev--Guass\end{minipage}%
\hfill%
\begin{minipage}{0.4\textwidth}\includegraphics[width=1\textwidth]{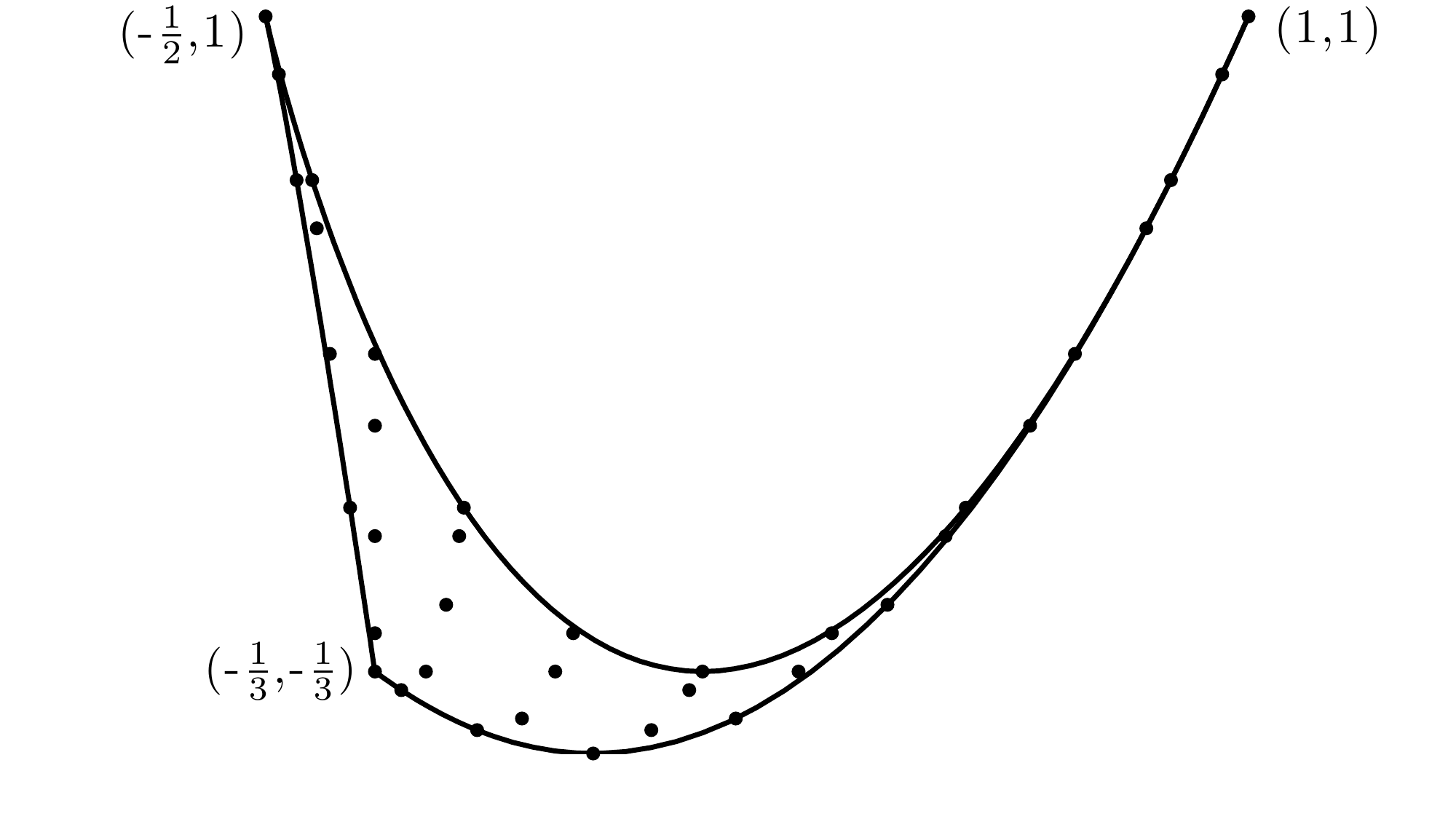}\\
\centering Chebyshev--Guass--Lobatto\end{minipage}%
\hspace*{\fill}

 \vspace*{0.5em}

\hfill%
\begin{minipage}{0.4\textwidth}\includegraphics[width=1\textwidth]{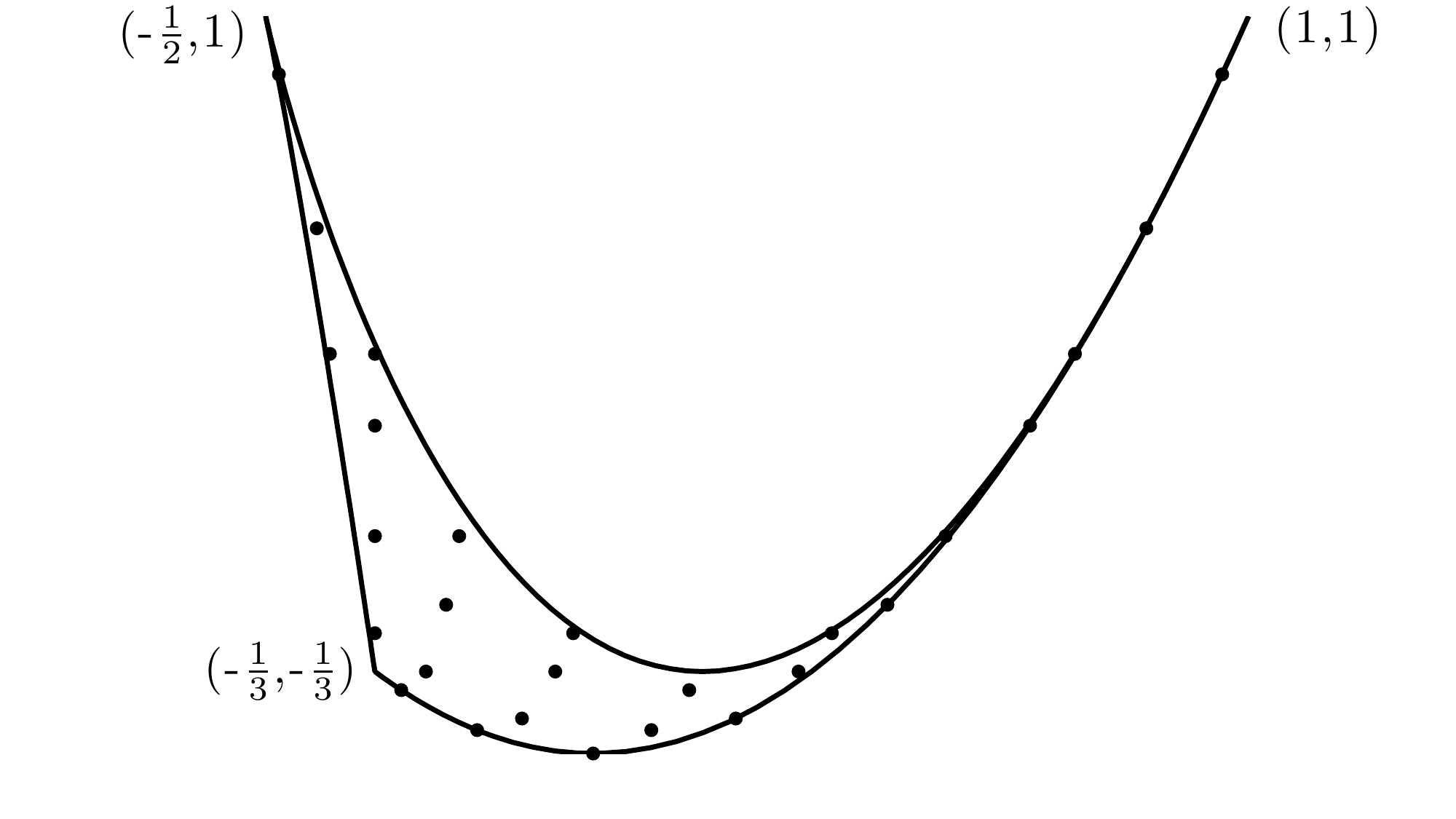}\\
\centering Chebyshev--Guass--Radau I\end{minipage}%
\hfill%
\begin{minipage}{0.4\textwidth}\includegraphics[width=1\textwidth]{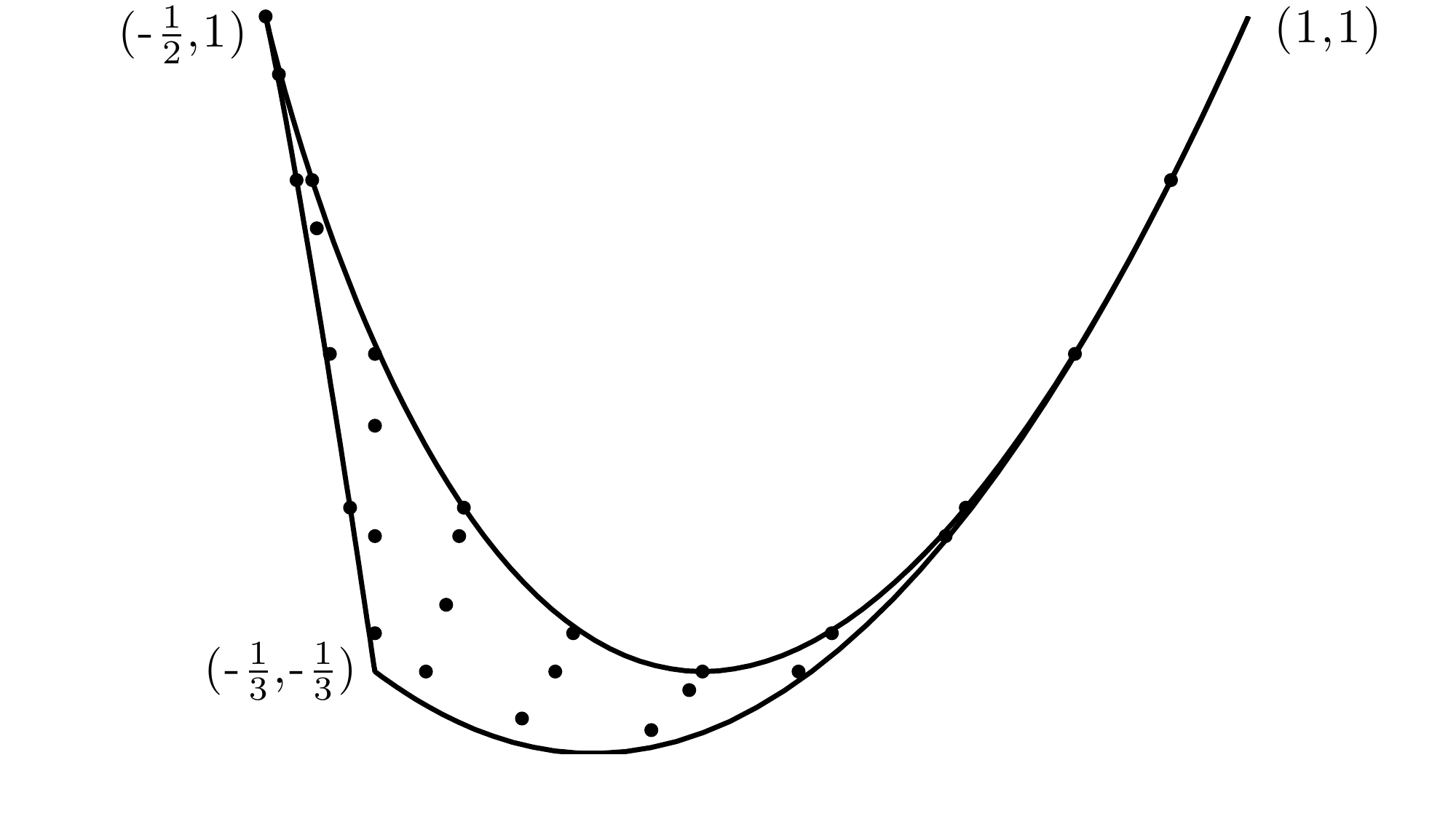}\\
\centering Chebyshev--Guass--Radau II\end{minipage}%
\hspace*{\fill}
\caption{The cubature nodes on the region $\Delta^*$.}
\end{figure}

\subsection[Gauss-Lobatto cubature and Chebyshev polynomials of the first kind]{Gauss--Lobatto cubature and Chebyshev polynomials of the f\/irst kind}

In the case of $w_{-\frac12, -\frac12}$,  the change of variables $\tb \mapsto x$ shows that
\eqref{cuba-HHT2} leads to a cubature of $m$-degree $2n-1$ based on the nodes of $Y_n$.

\begin{thm}
For the weight function $w_{-\frac{1}{2},-\frac{1}{2}}$ on $\triangle^*$ the
cubature rule
\begin{equation} \label{GaussCuba2}
  c_{- \frac{1}{2},-\frac12} \iint_{\triangle^*} f(x,y) w_{-\frac{1}{2},-\frac{1}{2}}(x,y) dxdy =
     \frac{1}{n^2}  \sum_{\jb\in \Upsilon_n}\omega^{(n)}_{\jb}
            f\big(x\big(\tfrac{\jb}{n}\big),y\big(\tfrac{\jb}{n}\big)\big)
\end{equation}
holds for $f \in \Pi_{2n-1}^*$.
\end{thm}

The set $Y_n$ includes points on the boundary of $\triangle^*$, hence, the cubature rule in
\eqref{GaussCuba2} is an analogue of the Gauss--Lobatto type cubature for $w_{-\frac12,-\frac12}$ on
$\triangle^*$. The number of nodes of this cubature is $\dim \Pi_n^*$, instead of $\dim \Pi_{n-1}^*$.
In this case, the corresponding orthogonal polynomials are the generalized Chebyshev polynomials
of the f\/irst kind, $T_{k_1,k_2}(x,y) : = P_n^{-\frac12,-\frac12}(x,y)$. The polynomials in
$\{T_\a: |\a|_* = n\}$ do not have enough common zeros in general. In fact, the two orthogonal polynomials
of $m$-degree $6$,
\begin{gather*}
T_{3,0}(x,y)= 36x^3-18xy-9x-6y-2,\\  T_{2,2}(x,y)=6y^2+10y-72x^3+36xy+18x+3.
\end{gather*}
only have three common zeros on $\triangle^*$,
\begin{gather*}
(x,y)= \Big(\tfrac{\sqrt{2}}{\sqrt{7}+1}\cos(\tfrac{2\pi \mu}{3} + \tfrac13\arccos \tfrac{3\sqrt{2}}{2\sqrt{7}+1}),-\tfrac{1}{\sqrt{7}+1} \Big), \qquad \mu=0,1,2,
\end{gather*}
whereas $\dim \Pi_5^*=5$. For cubature rules in the ordinary sense, that is, with $\Pi_n^2$ in place of~$\Pi_n^*$, the nodes of a cubature rule of degree $2n-1$ with $\dim \Pi_n^2$ nodes must be
the variety of a~polynomial ideal generated by $\dim \Pi^{*}_{n+1}$ linearly independent polynomials of degree
$n+1$, and these polynomials are necessarily quasi-orthogonal in the sense that they are orthogonal to all polynomials of degree $n-2$~\cite{X00}. Our next theorem shows that this characterization of such a cubature
carries over to the case of $m$-degree.

\begin{thm} 
Denote $\a^*=(\a_1  -1,\a_2)$, $a_1> a_2 $, and $\a^*=(\a_1, \a_1-1)$  if $\a_1=\a_2$. Then
$Y_n$ is the variety of the polynomial ideal
\begin{gather} \label{ideal}
  \left \langle T_{\a}(x) - T_{\a^*}(x): \  |\a|_* = n+1 \right \rangle.
\end{gather}
Furthermore, the polynomial $ T_{\a}(x) - T_{\a^*}(x)$ is of $m$-degree
$n+1$ and orthogonal to all polynomials in $\Pi_{n-2}^*$ with respect to $w_{-\frac12,-\frac12}$.
\end{thm}

\begin{proof}
A direct computation shows that, for any $\kb\in \Gamma$ with $k_1-k_3=n+1$,
\begin{gather*}
   \TCC_{k_1-1,k_2,k_3+1}(\tb)-\TCC_{\kb}(\tb)
 = \frac{1}{3} \Big[ \cos {\tfrac{\pi (n-1)(t_1-t_3)}{3}}\cos \pi k_2 t_2
    + \cos {\tfrac{\pi (n-1)(t_2-t_1)}{3}}\cos \pi k_2 t_3\\
\qquad\quad{}
    + \cos {\tfrac{\pi (n-1)(t_3-t_2)}{3}}\cos \pi k_2 t_1\Big]
-\frac{1}{3} \Big[ \cos {\tfrac{\pi (n+1)(t_1-t_3)}{3}}\cos \pi k_2 t_2\\
\qquad\quad{}
              + \cos {\tfrac{\pi (n+1)(t_2-t_1)}{3}}\cos \pi k_2 t_3
            + \cos {\tfrac{\pi (n+1)(t_3-t_2)}{3}}\cos \pi k_2 t_1\Big]
\\
  \qquad {}=  \frac{2}{3} \Big[ \sin {\tfrac{\pi n(t_1-t_3)}{3}} \sin {\tfrac{\pi(t_1-t_3)}{3}}   \cos \pi k_2 t_2
      + \sin {\tfrac{\pi n(t_2-t_1)}{3}}\sin {\tfrac{\pi(t_1-t_3)}{3}} \cos \pi k_2 t_3
 \\
   \qquad\quad{}            + \sin {\tfrac{\pi n(t_3-t_2)}{3}}\sin {\tfrac{\pi(t_1-t_3)}{3}} \cos \pi k_2 t_1\Big],
\end{gather*}
where we have used the def\/inition of $\TCC_{\kb}$ for the f\/irst equality sign.
Hence, for any $\jb \in \Upsilon_{n}$,
\begin{gather*}
   \TCC_{k_1-1,k_2,k_3+1}\left(\tfrac{\jb}{n}\right)-\TCC_{\kb}\left(\tfrac{\jb}{n}\right)
 =  \frac{2}{3} \Big[ \sin {\tfrac{\pi (j_1-j_3)}{3}} \sin {\tfrac{\pi(j_1-j_3)}{3n}}
  \cos \tfrac{\pi k_2 j_2}{n}\\
   \quad\qquad     {}  + \sin {\tfrac{\pi (j_2-j_1)}{3}} \sin {\tfrac{\pi(j_2-j_1)}{3n}}
  \cos \tfrac{\pi k_2 j_3}{n} + \sin {\tfrac{\pi(j_3-j_2)}{3}} \sin {\tfrac{\pi(j_3-j_2)}{3n}}
  \cos \tfrac{\pi k_2 j_1}{n}\Big] =0,
\end{gather*}
where the last equality sign uses the fact $j_1\equiv j_2\equiv j_3 \pmod{3}$.
With $\a_1 = k_1+k_2$, this shows that $T_\a - T_{\a^*}$ vanishes on $Y_n$. Finally, we
note that $|\a^*|_*=|\a|_*-2$ or $|\a|_*-1$,  so that $T_{\a^*}$ is a Chebyshev polynomial
of degree at least $n-1$ and $T_\a - T_{\a^*}$ is orthogonal to all polynomials in
$\Pi_{n-2}^*$.
\end{proof}

\subsection[Gauss-Radau cubature and Chebyshev polynomials of mixed kinds]{Gauss--Radau cubature and Chebyshev polynomials of mixed kinds}
Under the change of variables $\tb \mapsto x$ def\/ined in \eqref{xy}, we can also transform
\eqref{cuba-HHT2}  into cubature rules with respect to $w_{-\frac12, \frac12}$ and $w_{-\frac12, \frac12}$,
which have nodes on part of the boundary and are analogue of Gauss--Radau cubature rule. They are
associated with Chebyshev polynomials of the mixed types. We state the result without proof.

\begin{thm}
The following cubature rules hold,
\begin{gather}
  c_{- \frac{1}{2},\frac12}  \iint_{\Delta^*} f(x,y) w_{-\frac{1}{2},\frac{1}{2}}(x,y) dxdy \nonumber\\
 \qquad{} =
     \frac{4\pi^2}{9(n+2)^2}  \sum_{\jb\in \Upsilon_{n+2}}\omega^{(n+2)}_{\jb}\left|\TSC_{1,0,-1} \big(\tfrac{\jb}{n+2}\big)\right|^2
            f\big(x\big(\tfrac{\jb}{n+2}\big),y\big(\tfrac{\jb}{n+2}\big)\big), \qquad \forall \, f \in \Pi_{2n-1}^*,\!\!\label{GR1Cuba}
\\
  c_{\frac{1}{2},-\frac12}  \iint_{\Delta^*} f(x,y) w_{\frac{1}{2},-\frac{1}{2}}(x,y) dxdy \nonumber
 \\
 \qquad{} =
     \frac{4\pi^2}{9(n+3)^2}  \sum_{\jb\in \Upsilon_{n+3}}\omega^{(n+3)}_{\jb}\left|\TCS_{1,1,-2} \big(\tfrac{\jb}{n+3}\big)\right|^2
            f\big(x\big(\tfrac{\jb}{n+3}\big),y\big(\tfrac{\jb}{n+3}\big)\big), \qquad \forall\, f \in \Pi_{2n-1}^*.\!\! \label{GR2Cuba}
\end{gather}
\end{thm}

Since by \eqref{2equality}, $\TSC_{1,0,-1}$ and $\TCS_{1,1,-2}$ vanish on part of the boundary of
$\triangle$, the summation is not over the entire $\Upsilon_{n+2}$ or $\Upsilon_{n+3}$ but over a
subset that exclude points on the respective boun\-da\-ry. Let $Y_{n+1}^{\sc}$ and
$Y_{n+3}^{\cs}$ denote the set of nodes for the above two cubature rules, respectively.

\begin{thm} \label{prop2}
$Y_{n+2}^{\sc}$ is the variety of the polynomial ideal
\begin{gather} \label{ideal2}
  \big \langle P^{-\frac12,\frac12}_{\a}(x): \  |\a|_* = n \big \rangle.
\end{gather}
And
$Y_{n+3}^{\cs}$ is the variety of the polynomial ideal
\begin{gather} \label{ideal3}
  \big\langle P^{\frac12,-\frac12}_{\a}(x)-P^{\frac12,-\frac12}_{\a^*}(x): \  |\a|_* = n+1 \big \rangle.
\end{gather}
\end{thm}

It is of some interests to notice that, in terms of the number of nodes vs the degree,  \eqref{GR1Cuba}
is an analogue of the Gauss cubature rule in $m$-degree.

\subsection*{Acknowledgements}
The work of the f\/irst author was partially supported by
NSFC  Grants 10971212 and 91130014. The work of the second author was partially supported by
NSFC  Grant 60970089.  The work  of the third author was supported in part by NSF Grant DMS-110
6113 and  a grant from the Simons Foundation (\#~209057 to Yuan Xu).

\pdfbookmark[1]{References}{ref}
\LastPageEnding

\end{document}